\documentclass[final,3p,11pt]{elsarticle}
\usepackage[T1]{fontenc}
\usepackage{lmodern,amsmath,amssymb,amsthm,mathrsfs,booktabs,listings}
\usepackage{mathtools,microtype}
\usepackage{xr-hyper}
\usepackage[hidelinks]{hyperref}
\newtheorem{theorem}{Theorem}
\newtheorem{lemma}[theorem]{Lemma}

\newcommand{\QQ}{\mathbb Q}
\newcommand{\cGc}{\widetilde G_C}
\newcommand{\cGb}{\widetilde G_B}
\newcommand{\cLam}{\mathcal L}

\hypersetup{pdftitle={Counting formulas and refined generating functions for inversion-sequence Classes 152 and 207},
  pdfauthor={Haobo Ma, Toufik Mansour, Wenlin Zhang}}
\newtheorem{proposition}{Proposition}
\newcommand{\In}{\mathcal I}
\newcommand{\U}{\mathsf U}
\newcommand{\Ssite}{\mathsf S}
\newcommand{\ind}{\mathbf1}
\begin{document}
\begin{frontmatter}
\title{Counting formulas and refined generating functions for inversion-sequence Classes 152 and 207}
\author[chrono]{Haobo Ma}
\ead{auric@aelf.io}
\author[haifa]{Toufik Mansour}
\ead{tmansour@univ.haifa.ac.il}
\author[nus,chrono]{Wenlin Zhang}
\ead{e1327962@u.nus.edu}
\affiliation[chrono]{organization={ChronoAI Pte Ltd},city={Singapore},country={Singapore}}
\affiliation[haifa]{organization={Department of Mathematics, University of Haifa},city={Haifa},country={Israel}}
\affiliation[nus]{organization={National University of Singapore},country={Singapore}}
\begin{abstract}
For the two Class~152 pattern sets, we prove that the common number of
inversion sequences of length $n\ge1$ is
\[
 a_n=1+\frac12\displaystyle\sum_{j=1}^{n-1}\binom{2j}{j}.
\]
The Class~152 proof uses first-descent and last-zero decompositions,
with weighted Dyck paths evaluating the resulting generating functions.
For Class~207, we prove equality of the three total-count series and give
their common generating function as an explicit Lambert-series expression
after the substitution $z=q/(1+q)^2$. Complete legal-site and catalytic
recurrences lead to a Lambert series for one class and $q$-product orbit
sums for another. A cubic limit of the balanced nonterminating
${}_3\phi_2$ $q$-Saalsch\"utz summation identifies these expressions.
Fixed-maximum expansions, finite verification data, refined analytic
obstructions, and signed specializations are recorded in the supplementary
material.
\end{abstract}
\end{frontmatter}

\section{Introduction}
An inversion sequence of length $n$ is a word $e_0\cdots e_{n-1}$ with
$0\le e_i\le i$. Patterns occur in subsequences, with all equalities and
strict inequalities retained. Write $\In_n(X)$ for the sequences avoiding
the pattern set $X$, and $G_X(z)=\sum_{n\ge0}|\In_n(X)|z^n$.
The index $n$ is the number of entries. In the notation of
Callan--Mansour~\cite{CM}, a word $e_0\cdots e_m$ belongs to $I_m(X)$;
thus their index is one less than its length. We use the length convention
throughout.

The systematic study of length-three patterns in inversion sequences
was initiated by Corteel, Martinez, Savage and Weselcouch~\cite{CMSW}
and Mansour and Shattuck~\cite{MS}. Martinez and Savage~\cite{MSII}
organized a broader classification by triples of relations, and
Lin~\cite{Lin} studied such restrictions in inversion and ascent sequences.
For sets of ordinary patterns, Yan and Lin~\cite{YL} investigated pairs,
Callan, Jel\'inek and Mansour~\cite{CJM} triples, and Callan and
Mansour~\cite{CM} quadruples; the last classification is the source of the
two conjectures addressed here. Related work develops statistic refinements
for restricted inversion sequences~\cite{LK}, consecutive-pattern
enumeration~\cite{AEC}, and vincular Wilf equivalences~\cite{AEV}.
These adjacency-restricted occurrences use different conventions from the
subsequence occurrences considered here.

Generating trees give both structural descriptions and catalytic equations:
Kotsireas, Mansour and Y{\i}ld{\i}r{\i}m~\cite{KMY} give an algorithmic
approach for inversion sequences, while Testart~\cite{Testart} uses trees
growing on the left for further pattern classes. Britt and
Beaton~\cite{BrittBeaton} complete the enumeration of the triples-of-relations
classification using left- and right-growing trees. Within the quadruple
classification, Asinowski and Polley~\cite[Section~4.1]{AP} resolve
Class~166 through rectangulations. Here we resolve Classes~152 and~207;
the Class~207 proof evaluates the published B tree and our C legal-site
recurrence and identifies their total-count series by a classical
$q$-summation.

Callan--Mansour's numbering~\cite[Table~1, p.~28]{CM} is
\[
\begin{aligned}
207(1)=A&=\{100,101,110,201\},\\
207(2)=C&=\{100,101,110,210\},\\
207(3)=B&=\{101,110,120,210\}.
\end{aligned}
\]
Their bijection~\cite[Section~4.10]{CM} gives $G_A=G_C$: in the strict
left-to-right-maximum decomposition, leave the maximum blocks fixed and
refill the intervening boxes from left to right with the largest available
entry below the corresponding maximum.

For Class~152, cutting at the first descent describes one class by a
nondecreasing prefix and a suffix whose nonmaximum entries lie in one gap
between prefix values. Cutting at the last zero describes the other by a prefix with
increasing positive entries and a suffix of rotated increasing blocks.
Dyck paths encode the prefixes, and weights on their down-runs evaluate
both decompositions as the same generating function.

The more involved Class~207 step is $G_C(z)=G_B(z)$. The C legal-site
system gives a kernel equation whose two roots eliminate the diagonal
boundary; symmetry and a second-derivative condition then determine the
total as a Lambert series. The published B tree gives two $q$-product
orbit sums. To identify their difference, we specialize the balanced
nonterminating ${}_3\phi_2$ $q$-Saalsch\"utz summation
\cite[Eq.~17.7.5]{DLMF}. Its multiplier has a triple pole in the limiting
parameter. The linear and quadratic terms of the logarithmic product
expansion cancel, leaving the cubic term that produces the Lambert
series. The transformed tail retains the initial B orbit term; a uniform
geometric bound justifies termwise passage to the limit.
Section~\ref{sec:class207-equality} gives the
specialization and the passage back to formal length series.

The paper has two focused enumeration results. Section~\ref{sec:C} derives
the complete C legal-site recurrence and its catalytic equations. Sections
\ref{sec:B}--\ref{sec:Ckernel} give the corresponding B system and the two
kernel reductions. Section~\ref{sec:class207-equality} proves the common
Class~207 generating function, and Section~\ref{sec:class152} proves the
Class~152 formula by two direct decompositions. The expanded supplementary
material records the fixed-maximum rational expansion, numerical checks,
the auxiliary B-boundary calculation, refined analytic obstructions, signed
specializations, and the reproducibility code.

\section{An explicit recurrence for 207(2)}\label{sec:C}
For a nonempty $C$-avoiding prefix $w$ of length $n$, let $M=\max w$
and write $s=n-M\ge1$ for the slack, called the bounce statistic in
Asinowski--Polley~\cite[Section~2.2]{AP}. List its admissible next values at most $M$ in increasing
order as a word $p$: write $\Ssite$ for a value already seen and $\U$ for
an unused value. Let $h=1$ if $M$ is admissible, and $h=0$ otherwise.
All values $M+1,\ldots,n$ are admissible. The initial word $0$ has
state $(\Ssite,1,1)$.

We use the labelled generating-tree framework of West~\cite{West}.
Its translation into generating functions follows the general approach of
Banderier et al.~\cite{Banderier}. The complete succession rule is as follows;
each choice is counted once:
\begin{align}
(p,s,h)&\longrightarrow(p\U^{j-1}\Ssite,s-j+1,1),
                   &&1\le j\le s;\label{eq:newmax}\\
(p,s,1)&\longrightarrow(\Ssite,s+1,1),
                   &&\text{repeat }M;\label{eq:repeat}\\
(p,s,h)&\longrightarrow(\U^k,s+1,0),
                   &&\text{select a site below }M.\label{eq:lower}
\end{align}
In the last line, $k$ is the number of unused sites strictly to the right
of the selected site. An equal pair $a,a$ forbids a later entry below $a$
by $110$. A descent $a>b$ forbids entries at most $b$ by $100,210$, and
forbids $a$ by $101$. Thus selecting $x<M$ removes all sites at most $x$
and all seen sites above $x$, leaving precisely the $k$ unused sites.
Repeating $M$ leaves only its own site; a new maximum introduces no new
exclusion. This proves the rule and shows why gaps between listed sites
do not affect it.

For $r\ge0$, $s\ge1$ and $k\ge0$, let $a_r(s,k)$ count continuations
of $r$ further entries from the formal state $(\U^k,s,0)$, and let $b_r(s)$
count them from $(\Ssite,s,1)$. The rules above define these numbers
even when a state is not reached from the root.

\begin{proposition}\label{prop:recurrence}
The initial conditions are $a_0(s,k)=b_0(s)=1$. For $r\ge0$,
\begin{align}
a_{r+1}(s,k)
 &=\sum_{i=0}^{k-1}a_r(s+1,i)\notag\\
 &\quad+\sum_{j=1}^{s}\left(a_r(s-j+1,k+j-1)
                         +b_r(s-j+1)-a_r(s-j+1,0)\right),\label{eq:a}\\
b_{r+1}(s)
 &=b_r(s+1)
   +\sum_{j=1}^{s}\left(a_r(s-j+1,j)
                         +b_r(s-j+1)-a_r(s-j+1,0)\right).\label{eq:b}
\end{align}
Empty sums are zero. In particular,
\begin{equation}\label{eq:counts}
|\In_0(C)|=1,\qquad |\In_n(C)|=b_{n-1}(1)\quad(n\ge1).
\end{equation}
These are recurrences for nonnegative integers; the differences in the
formulas come from an affine reduction of the site word.
\end{proposition}

\begin{proof}
For a site $a$ in $p$, let $k(a)$ count unused sites strictly to its right,
and put $T_p(y)=\sum_{a\in p}y^{k(a)}$. Write $\phi_r(p,s,h)$ for the
number of continuations of length $r$. We first prove that it is affine
in $h$ and the coefficients of $T_p$, for fixed $r,s$.
At $r=0$ it is $1$. For the induction step, lower selections and a
possible repetition of $M$ contribute
\[
\sum_{a\in p}a_r(s+1,k(a))
       +h\bigl(b_r(s+1)-a_r(s+1,0)\bigr).
\]
For a new maximum, the child's site polynomial is
\begin{equation}\label{eq:Ttransform}
T_{p\U^{j-1}\Ssite}(y)
 =y^{j-1}T_p(y)+\frac{1-y^{j-1}}{1-y}+1,
\end{equation}
so its contribution is also affine by induction. Evaluating this affine
function on the empty word, on $\U^k$, and on $\Ssite$ gives
\begin{equation}\label{eq:affine}
\phi_r(p,s,h)=a_r(s,0)+h\bigl(b_r(s)-a_r(s,1)\bigr)
 +\sum_{a\in p}\bigl(a_r(s,k(a)+1)-a_r(s,k(a))\bigr).
\end{equation}
Indeed, $T_{\U^{k+1}}-T_{\U^k}=y^k$ and
$T_{\Ssite}=T_{\U}=1$. Apply \eqref{eq:affine} to the new-maximum
children $\U^{k+j-1}\Ssite$ and $\Ssite\U^{j-1}\Ssite$; their
continuation counts are respectively
$a_r(s',k+j-1)+b_r(s')-a_r(s',0)$ and
$a_r(s',j)+b_r(s')-a_r(s',0)$, where $s'=s-j+1$.
The lower children of $\U^k$ are $\U^0,\ldots,\U^{k-1}$;
the repetition child of $\Ssite$ is $\Ssite$ with slack $s+1$.
This proves \eqref{eq:a}--\eqref{eq:b}. The root has length one,
which proves \eqref{eq:counts}.
\end{proof}

For example, $a_1(s,k)=s+k$ and $b_1(s)=s+1$.
To calculate all counts through length $N$, at depth $r$ it suffices to
retain $a_r(s,k)$ for $s+k\le N-r$ and $b_r(s)$ for $s\le N-r$.
Row and diagonal partial sums evaluate \eqref{eq:a}--\eqref{eq:b} in
$O(N^3)$ integer operations and $O(N^2)$ stored integers;
\ref{supp-sec:program} of the supplement gives the complete program.

\subsection*{Catalytic equations and their proof}
For $C$, let the sums below run over all nonempty $C$-avoiding inversion
prefixes. Define
\[
F_C(x)=\sum_w z^{|w|}x^{s(w)},\qquad
P_C(x)=\sum_w z^{|w|}x^{s(w)}h(w),\qquad
H_C(x,y)=\sum_w z^{|w|}x^{s(w)}T_{p(w)}(y).
\]
Throughout this paper, $F_X(x)$ and $P_X(x)$ are univariate in the
catalyst $x$, whereas $H_X(x,y)$ and $Z_X(x,y)$ are bivariate in the
catalysts $x,y$. The length variable $z$ is suppressed in these
arguments, for example $F_X(x)=F_X(z;x)$. A specialization such as
$H_X(x,1)$ or $H_X(y,y)$ always retains both argument positions.
Put
\[
J_C(x)=zx\frac{F_C(x)-F_C(1)}{x-1}.
\]
The succession rule gives
\begin{align}
F_C(x)&=zx+zxH_C(x,1)+J_C(x),\label{eq:CF}\\
P_C(x)&=zx+zxP_C(x)+J_C(x),\label{eq:CP}\\
H_C(x,y)&=zx+zxP_C(x)
 +\frac{zx}{1-y}\bigl(H_C(x,1)-H_C(x,y)\bigr)\notag\\
&\quad+\frac{zx}{x-y}\bigl(H_C(x,y)-H_C(y,y)\bigr)
 +\frac{2-y}{1-y}J_C(x)
 -\frac{zx\bigl(F_C(x)-F_C(y)\bigr)}{(1-y)(x-y)}.\label{eq:CH}
\end{align}
The coefficient ring is $\mathbb Q[x,y][[z]]$; for a univariate
series use $\mathbb Q[x][[z]]$. Apparent singularities in these equations
are removable coefficientwise. In particular, terms with denominator
$1-y$ must be combined before setting $y=1$. All constant terms in $z$
are zero, and
\[
G_C(z)=1+F_C(1)=1+\sum_{r\ge0}b_r(1)z^{r+1}.
\]

\begin{proof}[Proof of \eqref{eq:CF}--\eqref{eq:CH}]
The root contributes $zx$ to each series. For a parent with slack $s$,
the new maxima have slacks $s-j+1$, $1\le j\le s$, and
\begin{equation}\label{eq:sumidentities}
\sum_{j=1}^s x^{s-j+1}=\frac{x(x^s-1)}{x-1},\qquad
\sum_{j=1}^s x^{s-j+1}y^{j-1}
 =\frac{x(x^s-y^s)}{x-y}.
\end{equation}
The first identity gives $J_C(x)$ in $F_C(x)$ and $P_C(x)$.
Every admissible site at most $M$ gives a child with slack $s+1$,
so their total contribution to $F_C(x)$ is $zxH_C(x,1)$.
Among them, only repetition of an admissible $M$ has an admissible
maximum afterward; it contributes $zxP_C(x)$ to $P_C(x)$.
These observations prove \eqref{eq:CF} and \eqref{eq:CP}.

For $H_C(x,y)$, repetition of $M$ leaves the one-site word $\Ssite$,
so contributes $zxP_C(x)$. A lower site $a$ leaves $\U^{k(a)}$,
whose polynomial is $(1-y^{k(a)})/(1-y)$. Summing over all sites gives
$zx(H_C(x,1)-H_C(x,y))/(1-y)$: the maximum site, when present,
has $k(M)=0$ and contributes zero to this expression.
For a new maximum, \eqref{eq:Ttransform} separates the child's
polynomial into the old-site term and the term independent of the
old sites. The old-site term, summed using the second identity in
\eqref{eq:sumidentities}, is
$zx(H_C(x,y)-H_C(y,y))/(x-y)$.
The independent term is $(2-y-y^{j-1})/(1-y)$; its sum is
\[
\frac{2-y}{1-y}J_C(x)
 -\frac{zx(F_C(x)-F_C(y))}{(1-y)(x-y)}.
\]
This accounts for every child exactly once and proves \eqref{eq:CH}.
Both identities in \eqref{eq:sumidentities} are finite geometric sums;
thus each displayed quotient, after its indicated combination, is a
polynomial at each length. Finally the factor $z$ on each child
contribution determines the coefficients recursively from the root.
\end{proof}

\section{207(3) and a common catalytic form}\label{sec:B}

Here is a direct four-type generating tree for $B$. Let $M,L$ be the
largest and second largest entries, counting multiplicity; for the
initial word $0$ put $L=0$. Define
\[
D(w)=\{w_i:w_i>w_j\text{ for some }j>i\},\qquad
s=n-M,\quad d=M-L.
\]
The types are $R_s$ if $M=L$, and, when $M>L$,
\[
\begin{array}{c|cc}
\text{type}&M\in D(w)&L\in D(w)\\\hline
U_{s,d}&\text{no}&\text{no}\\
V_{s,d}&\text{yes}&\text{no}\\
W_{s,d}&\text{no}&\text{yes}
\end{array}
\]
The root is $R_1$, and the complete succession rule is
\begin{align}
R_s&\longrightarrow R_{s+1},\quad U_{s-j+1,j}\ (1\le j\le s),\notag\\
U_{s,d}&\longrightarrow R_{s+1},\quad V_{s+1,k}\ (1\le k\le d),
                         \quad U_{s-j+1,j}\ (1\le j\le s),\notag\\
V_{s,d}&\longrightarrow V_{s+1,k}\ (1\le k\le d),
                         \quad W_{s-j+1,j}\ (1\le j\le s),\label{eq:Btree}\\
W_{s,d}&\longrightarrow R_{s+1},\quad V_{s+1,k}\ (1\le k<d),
                         \quad U_{s-j+1,j}\ (1\le j\le s).\notag
\end{align}
To prove the rule, an earlier pair $a\le b$ forbids an appended $x<a$
by $110$ or $120$; a descent $a>b$ forbids $x<b$ by $210$ and $x=a$
by $101$. To identify their intersection, the pair consisting of the two largest
entries forbids every $x<L$, while every earlier pair has its smaller
entry at most $L$ and forbids nothing larger by its inequality condition.
The equality exclusions are exactly $D(w)$. Thus the legal next entries
are $[L,n]\setminus D(w)$ (also for the single-entry root).
If $M>L$, the maximum occurs once. The two endpoints cannot both
belong to $D(w)$: if every $L$ occurs after $M$, a descent from $L$
gives $M,L,v$ with $v<L$, a $210$. If some $L$ occurs before $M$,
a descent from $M$ gives a later $v\le L$; then $L,M,v$ is $120$
when $v<L$, and $101$ when $v=L$.
If $M=L$ and the maximum is repeated, a descent after its second
occurrence gives $110$, and a descent before that occurrence followed
by the second maximum gives $101$. Thus the maximum remains admissible.
No earlier value lies strictly between $L$ and $M$, so only the two
endpoints can be excluded from $[L,M]$. Appending $M+j$ has parameters
$(s-j+1,j)$ and type $W$ precisely when the old $M$ was excluded.
Repeating a legal $M$ gives $R_{s+1}$. Appending $M-k<M$ gives
$V_{s+1,k}$ in the ranges displayed above. This also proves inductively
that the four types exhaust all prefixes. The tree in
\cite[Table~1, p.~28]{CM} gives another enumeration of $B$.

\subsection{The two explicit systems}
For $X=B,C$, use the admissible sites of $X$ and the definitions
of $p,h,T_p$ in Section~\ref{sec:C}. Define
\[
\begin{aligned}
F_X(x)&=\sum_{w\in\bigcup_{n\ge1}\In_n(X)}z^{|w|}x^{s(w)},&
P_X(x)&=\sum_{w\in\bigcup_{n\ge1}\In_n(X)}z^{|w|}x^{s(w)}h(w),\\
H_X(x,y)&=\sum_{w\in\bigcup_{n\ge1}\In_n(X)}
 z^{|w|}x^{s(w)}T_{p(w)}(y),&
Z_X(x,y)&=H_X(x,y)-P_X(x).
\end{aligned}
\]
Thus $Z_X(x,y)$ sums the site polynomials with the maximum site removed:
that site, if admissible, has weight $y^0=1$.
In both classes $x$ marks the slack $s=n-M$, and a site's exponent
of $y$ counts unused admissible sites strictly to its right. In
particular, $y$ is a catalytic variable recording the unused sites above a selected site; a prefix contributes a polynomial in $y$, rather than necessarily a single monomial.
The unknowns in each system are $P_X(x),F_X(x),Z_X(x,y)$;
$H_X(x,y)=P_X(x)+Z_X(x,y)$ is then recovered.

\begin{proposition}\label{prop:common}
For $X=B$, the complete system in ordinary functional-equation form is
\begin{subequations}\label{eq:Bsystem}
\begin{align}
P_B(x)&=zx+zxP_B(x)+zx\frac{F_B(x)-F_B(1)}{x-1},\label{eq:BP}\\
F_B(x)&=P_B(x)+zxZ_B(x,1),\label{eq:BF}\\
Z_B(x,y)&=zx\left[\frac{Z_B(x,1)-Z_B(x,y)}{1-y}+Z_B(x,y)
 +\frac{P_B(x)-P_B(y)}{x-y}\right.\notag\\
&\hspace{16mm}\left.
 +\frac{1}{1-y}\left(\frac{F_B(x)-F_B(1)}{x-1}
                    -\frac{F_B(x)-F_B(y)}{x-y}\right)\right].\label{eq:BZ}
\end{align}
\end{subequations}
For $X=C$, it is
\begin{subequations}\label{eq:Csystem}
\begin{align}
P_C(x)&=zx+zxP_C(x)+zx\frac{F_C(x)-F_C(1)}{x-1},\label{eq:CPexplicit}\\
F_C(x)&=P_C(x)+zxZ_C(x,1),\label{eq:CFexplicit}\\
Z_C(x,y)&=zx\left[\frac{Z_C(x,1)-Z_C(x,y)}{1-y}
 +\frac{Z_C(x,y)-Z_C(y,y)}{x-y}
 +\frac{P_C(x)-P_C(y)}{x-y}\right.\notag\\
&\hspace{16mm}\left.
 +\frac{1}{1-y}\left(\frac{F_C(x)-F_C(1)}{x-1}
                    -\frac{F_C(x)-F_C(y)}{x-y}\right)\right].\label{eq:CZ}
\end{align}
\end{subequations}
All constant terms in $z$ are zero. The root conditions are
\[
[z]P_X(x)=[z]F_X(x)=x,\qquad [z]Z_X(x,y)=0,\qquad [z]H_X(x,y)=x.
\]
Each system has a unique solution in the coefficient rings stated
above, and $G_X(z)=1+F_X(1)$.
\end{proposition}

\begin{proof}
For either class, a new maximum is admissible and contributes
$J_X(x)=zx(F_X(x)-F_X(1))/(x-1)$.
Repetition of an admissible maximum leaves a single seen site and
contributes $zxP_X(x)$. A lower choice makes the old maximum
inadmissible: it creates a descent ending below that maximum and
$101$ forbids repeating it. Including the root proves
\eqref{eq:BP} and \eqref{eq:CPexplicit}.
Every admissible site gives a child of slack $s+1$, and every new
maximum is already counted in $J_X(x)$, so
\[
F_X(x)=zx+zxH_X(x,1)+J_X(x).
\]
Subtracting the equation for $P_X(x)$ proves \eqref{eq:BF} and
\eqref{eq:CFexplicit}.

For $C$, use the proved equation \eqref{eq:CH} and substitute
$H_C(x,y)=P_C(x)+Z_C(x,y)$. Its old-site new-maximum term becomes
\[
\frac{zx(P_C(x)-P_C(y))}{x-y}
 +\frac{zx(Z_C(x,y)-Z_C(y,y))}{x-y}.
\]
The lower-site difference loses its $P_C(x)$ terms. Subtracting
\eqref{eq:CPexplicit} leaves the independent source
\[
\frac{zx}{1-y}\left(\frac{F_C(x)-F_C(1)}{x-1}
                         -\frac{F_C(x)-F_C(y)}{x-y}\right),
\]
which proves \eqref{eq:CZ}.

For $B$, the legal interval is $[L,n]\setminus D(w)$, as proved
above. A lower choice $a$ leaves one seen site at $a$ and the $k(a)$
unused sites to its right; the old maximum is forbidden. Its site
polynomial is therefore $(1-y^{k(a)+1})/(1-y)$.
Repetition of $M$ leaves one site, the same expression for $k(M)=0$.
A new maximum $M+j$ discards every old site below $M$: the old
maximum becomes the second largest entry. The old maximum itself
remains admissible exactly when $h=1$. Consequently the child's
polynomial is
\[
hy^{j-1}+\frac{1-y^{j-1}}{1-y}+1.
\]
The finite sums \eqref{eq:sumidentities} now give the full equation
\begin{align}
H_B(x,y)&=zx+\frac{zx}{1-y}\bigl(H_B(x,1)-yH_B(x,y)\bigr)
 +\frac{zx(P_B(x)-P_B(y))}{x-y}\notag\\
&\quad+\frac{2-y}{1-y}J_B(x)
 -\frac{zx(F_B(x)-F_B(y))}{(1-y)(x-y)}.\label{eq:BH}
\end{align}
Inserting $H_B(x,y)=P_B(x)+Z_B(x,y)$ and subtracting
\eqref{eq:BP}, the site term left is
\[
\frac{Z_B(x,1)-yZ_B(x,y)}{1-y}
 =\frac{Z_B(x,1)-Z_B(x,y)}{1-y}+Z_B(x,y).
\]
The independent source is the same fraction involving $F_B$ as in
\eqref{eq:BZ}. This proves that equation too.

For completeness, the last fraction in either system is polynomial:
for $f(x)=x^a$ with $a\ge2$ it is
\[
\frac{1}{1-y}\left(\frac{x^a-1}{x-1}
                      -\frac{x^a-y^a}{x-y}\right)
 =\sum_{\substack{i,j\ge0\\i+j\le a-2}}x^iy^j,
\]
and it is zero for $a=0,1$. The other divided differences are also
polynomial finite sums. Starting with the constant terms zero,
compute the coefficient of $P_X(x)$ and $Z_X(x,y)$ at length $n$
from length $n-1$, then obtain that of $F_X(x)$ from its displayed
equation. This proves uniqueness and the root conditions.
The definition of $F_X(x)$ proves $G_X(z)=1+F_X(1)$.
\end{proof}

The identity proved below concerns the specialization $x=1$.

\section{Kernel solution for 207(3)}\label{sec:Bkernel}
The elimination below uses catalytic variables and iterated kernel roots,
as in the generating-tree applications of the kernel method surveyed and
developed by Bousquet-M\'elou~\cite{BM}. We give the specialized elimination
for the present succession rule explicitly.

Write $b_B=F_B(1)$, $\beta_B=b_B+(z-1)/z$, and put
\[
\phi(x)=\frac1{1-zx},\qquad
D(x)=zx^2-x+1,\quad A(x)=D(x)+zx,\quad M(x)=zx+z-1.
\]
Here $D,A,M$ are auxiliary polynomials, distinct from the forbidden
set $D(w)$ in the generating tree.
Clearing denominators in \eqref{eq:BZ} gives
\begin{align}
(1-y+zxy)Z_B(x,y)
={}&zxZ_B(x,1)+zx(1-y)\frac{P_B(x)-P_B(y)}{x-y}\notag\\
&+zx\left(\frac{F_B(x)-b_B}{x-1}
                    -\frac{F_B(x)-F_B(y)}{x-y}\right).\label{eq:Bcleared}
\end{align}
The kernel vanishes at $y=\phi(x)$; this substitution is valid
coefficientwise because the catalytic coefficients are polynomials.
Solving \eqref{eq:BP} first gives
\begin{equation}\label{eq:PXsolution}
P_X(x)=\frac{zx}{1-zx}
 \left(1+\frac{F_X(x)-F_X(1)}{x-1}\right),\qquad X=B,C.
\end{equation}
With $Y=\phi(x)$, the vanishing-kernel condition, multiplied by $x-Y$,
is
\begin{align*}
0={}&(x-Y)\frac{F_B(x)-P_B(x)}{zx}
 +(1-Y)(P_B(x)-P_B(Y))\\
&+(x-Y)\frac{F_B(x)-b_B}{x-1}-F_B(x)+F_B(Y).
\end{align*}
Inserting \eqref{eq:PXsolution} and $Y=1/(1-zx)$ and clearing
rational denominators yields the scalar equation
\begin{equation}\label{eq:Bscalar}
M(x)A(x)F_B(x)+zx(1-zx)^2F_B(\phi(x))
 =xz^2D(x)\beta_B.
\end{equation}
These operations are identities over $\mathbb Q(z,x)$ with the three
series values $F_B(x),F_B(Y),b_B$ treated as independent symbols.
More explicitly, after the substitutions the right-hand side of the
vanishing-kernel display is $-1/(zx(1-zx)M(x))$ times the difference
between the left and right sides of \eqref{eq:Bscalar}. This verifies
all factors and signs in the reduction.

Define
\[
\lambda=\frac{1-\sqrt{1-4z}}{2z},\qquad
\alpha=\frac{1-z-\sqrt{1-6z+z^2}}{2z},\qquad \rho=\phi(\alpha),
\]
where the square roots have constant term $1$. Then $D(\lambda)=0$,
$A(\alpha)=0$, and both $\lambda$ and $\alpha$ have constant term $1$.
For the rational functions
\[
f(x)=\frac{xz^2D(x)}{M(x)A(x)},\qquad
 g(x)=-\frac{zx(1-zx)^2}{M(x)A(x)},
\]
put $a_i=\phi^i(a)$ and define
\begin{equation}\label{eq:orbits}
\Pi(a)=\prod_{i\ge0}g(a_i),\qquad
\mathscr S(a)=\sum_{i\ge0}f(a_i)\prod_{0\le j<i}g(a_j).
\end{equation}
The empty product is $1$. Below these are coefficientwise products
and sums, for $a=1$, $a=\rho$, or an indeterminate $a=x$.

\begin{proposition}\label{prop:Bproduct}
Set
\[
\tau=-z^2\alpha\rho^2,\qquad
\mathscr Q=\mathscr S(1)+\frac{\Pi(1)}{\Pi(\rho)}
                     \bigl(\tau-\mathscr S(\rho)\bigr).
\]
The orbit expressions are well-defined formal series,
$\mathscr Q=-z^2+O(z^3)$, and
\begin{equation}\label{eq:GBproduct}
G_B(z)=1+\frac{(1-z)\mathscr Q}{z(\mathscr Q-1)}.
\end{equation}
With $b_B=G_B(z)-1$ and $\beta_B=b_B+(z-1)/z$, the full slack series is
\begin{equation}\label{eq:FBproduct}
F_B(x)=\beta_B\left[\mathscr S(x)+\frac{\Pi(x)}{\Pi(\rho)}
                          \bigl(\tau-\mathscr S(\rho)\bigr)\right].
\end{equation}
The series $P_B(x)$ is given by \eqref{eq:PXsolution}. The remaining
series is recovered explicitly as
\begin{align}
Z_B(x,y)=\frac{1}{1-y+zxy}\left[F_B(x)-P_B(x)
 +zx(1-y)\frac{P_B(x)-P_B(y)}{x-y}\right.\notag\\
\left.\hspace{17mm}
 +zx\left(\frac{F_B(x)-b_B}{x-1}
                    -\frac{F_B(x)-F_B(y)}{x-y}\right)\right],\label{eq:ZBsolution}
\end{align}
and $H_B(x,y)=P_B(x)+Z_B(x,y)$.
\end{proposition}

\begin{proof}
The fixed point $\phi(\lambda)=\lambda$ and the identity
\[
\phi(u)-\lambda=\frac{z(u-\lambda)}{(1-zu)(1-z\lambda)}
\]
show that each orbit step raises the $z$-order of $u-\lambda$ by one.
The initial orders of $1-\lambda$ and $\rho-\lambda$ are $1$ and $2$,
respectively. Along either orbit $M(a_i)$ is a unit and $A(a_i)$ has
order $1$. Also $g(\lambda)=1$ and $f(\lambda)=0$.
To see the orders explicitly, factor
\[
D(u)=(u-\lambda)\bigl(z(u+\lambda)-1\bigr).
\]
Then the factor $z^2/A(u)$ in $f(u)$ raises the order by one.
For $g(u)-1$, its numerator is
$-zu(1-zu)^2-M(u)A(u)$, which vanishes at $u=\lambda$;
division by $A(u)$ can lower the order by at most one.
Consequently, on the $1$ orbit,
\[
g(a_i)-1=O(z^i)\ (i\ge1),\qquad f(a_i)=O(z^{i+2}),
\]
and on the $\rho$ orbit the respective orders are at least $i+1$
and $i+3$. Direct expansion gives $g(1)=1/2+O(z)$ and
$g(\rho)=1+O(z)$. The products therefore converge and are units,
with constant terms $1/2$ and $1$, and the sums converge.
For an indeterminate $x$, the first iterate has constant term $1$;
$A(\phi(x))=z(2-x)+O(z^2)$ is nonzero over $\mathbb Q(x)$.
After the second iterate the same order estimates apply, proving
convergence over $\mathbb Q(x)[[z]]$ as well.

Equation \eqref{eq:Bscalar} is
\[
F_B(u)=\beta_B f(u)+g(u)F_B(\phi(u)).
\]
After $k$ iterations this is a finite sum plus
$\prod_{i<k}g(a_i)F_B(a_k)$.
Since each coefficient of $F_B(u)$ is polynomial in $u$ and
$a_k\longrightarrow\lambda$, its coefficientwise limit is
\begin{equation}\label{eq:orbitlimit}
F_B(a)=\beta_B\mathscr S(a)+\Pi(a)F_B(\lambda).
\end{equation}
At $u=\alpha$, the scalar equation is valid before dividing by
$A(u)$. As $D(\alpha)=-z\alpha$, it gives
\[
F_B(\rho)=-\frac{z^2\alpha}{(1-z\alpha)^2}\beta_B=\tau\beta_B.
\]
Use \eqref{eq:orbitlimit} at $a=\rho$ to eliminate $F_B(\lambda)$.
At $a=1$ the result is $b_B=\beta_B\mathscr Q$. Solving this
linear equation using $\beta_B=b_B-(1-z)/z$ proves
\eqref{eq:GBproduct}; using $a=x$ proves \eqref{eq:FBproduct}.
The first terms are $\mathscr S(1)=-z^2/2+O(z^3)$,
$\mathscr S(\rho)=O(z^3)$ and $\tau=-z^2+O(z^3)$, so
$\mathscr Q=-z^2+O(z^3)$ and $\mathscr Q-1$ is a unit.
Finally \eqref{eq:PXsolution} and \eqref{eq:Bcleared} prove
\eqref{eq:ZBsolution}. Its apparent catalytic denominators are
interpreted after cancellation; the recovered coefficients are the
polynomials of the unique solution to Proposition~\ref{prop:common}.
\end{proof}

\section{The kernel reduction for 207(2)}\label{sec:Ckernel}
Clearing denominators in \eqref{eq:CZ} gives
\begin{align}
\kappa(x,y)Z_C(x,y)
={}&zx\left[(x-y)Z_C(x,1)-(1-y)Z_C(y,y)
 +(1-y)(P_C(x)-P_C(y))\right.\notag\\
&\left.\hspace{12mm}
 +(x-y)\frac{F_C(x)-F_C(1)}{x-1}-F_C(x)+F_C(y)\right],\label{eq:Ccleared}\\
\kappa(x,y)&=(1-y)(x-y)+zx(x-1).\notag
\end{align}
Let $b_C=F_C(1)$ and $\beta_C=b_C+(z-1)/z$. Define
\[
W_C(x)=\frac{zx(1-x)}{1-zx}F_C(x),\qquad
R(y)=\frac{z(1-y)^2}{(2-y)\{2-y+z(y^2-y-1)\}}.
\]
The following gives a scalar reduction with both kernel roots kept
inside legitimate formal compositions.

\begin{proposition}\label{prop:Ckernel}
For an indeterminate $t$, put
\[
T(t)=\frac{t(zt-t+1)}{zt^2-t+1},\qquad
x_0=T(t),\quad x_1=-\frac{t-1}{zt^2-t+1}=T(\phi(t)),\quad y_0=T(t)/t.
\]
Then $x_0,x_1$ are the two roots of $\kappa(x,y_0)=0$ and
\begin{equation}\label{eq:Cadditive}
W_C(x_0)-W_C(x_1)=\beta_C(x_0-x_1)R(y_0).
\end{equation}
\end{proposition}

\begin{proof}
Substitution into the quadratic verifies $\kappa(x_i,y_0)=0$.
The constant terms of $x_0,x_1,y_0$ are $t,1,1$, respectively;
thus compositions with the original prefix series are defined over
$\mathbb Q(t)[[z]]$. Rational divisions below may be taken in its
Laurent-series field.
In \eqref{eq:CH}, set $\kappa(x,y)=0$ and use
\eqref{eq:CF}, \eqref{eq:CP}, and \eqref{eq:PXsolution}.
The common diagonal boundary is
\begin{equation}\label{eq:Cboundary}
H_C(y,y)-\frac{F_C(y)}{1-y}
 =a(x,y)F_C(x)+q(x,y)b_C+c(x,y),
\end{equation}
where
\begin{align*}
a(x,y)&=\frac{x-y}{(x-1)(1-zx)}
                      +\frac{x-y-zx}{zx(1-y)},\\
q(x,y)&=-\frac{x-y}{(x-1)(1-zx)},\qquad
c(x,y)=\frac{x-y}{1-zx}-\frac{x-y}{1-y}.
\end{align*}
For an explicit derivation, the cleared equation for $H_C(x,y)$
on the kernel is
\begin{align*}
H_C(y,y)-\frac{F_C(y)}{1-y}
={}&(x-y)(1+P_C(x))
 +\frac{x-y}{1-y}H_C(x,1)\\
&+\frac{(2-y)(x-y)}{zx(1-y)}J_C(x)-\frac{F_C(x)}{1-y}.
\end{align*}
Insert $zxH_C(x,1)=F_C(x)-zx-J_C(x)$ and
$1+P_C(x)=(1+J_C(x))/(1-zx)$; the displayed $a,q,c$ result.

In the following root calculations set $y=y_0$. On the kernel, $(x-y)/(x-1)=-zx/(1-y)$, which gives
\[
a(x,y)\frac{1-zx}{zx(1-x)}=-\frac{2-y}{(1-y)^3}.
\]
The root sum and product of
$zx^2+(1-y-z)x-y(1-y)$ give
\[
(1-zx_0)(1-zx_1)=2-y+z(y^2-y-1).
\]
Also, writing $\delta q=q(x_0,y)-q(x_1,y)$, the kernel relation gives
\[
\delta q=\frac{z(x_0-x_1)}{(1-y)(1-zx_0)(1-zx_1)},\qquad
c(x_0,y)-c(x_1,y)=\frac{z-1}{z}\delta q.
\]
Subtract \eqref{eq:Cboundary} for the two roots; its left sides
are identical. The result is
\[
\frac{2-y}{(1-y)^3}\bigl(W_C(x_0)-W_C(x_1)\bigr)
 =\beta_C\delta q.
\]
The preceding product formula proves \eqref{eq:Cadditive} at $y=y_0$.
\end{proof}

The $C$ scalar equation has an additive forcing term, whereas the
 $B$ scalar equation has the convergent product solution above. The formal-iteration
warning and the coefficientwise fixed-maximum expansion are recorded in the
supplementary material; the Lambert-series calculation below supplies the
all-length boundary value used in the equality proof.

\section{The Class 207 equality and its common generating function}
\label{sec:class207-equality}
The C kernel equation gives a Lambert series quotient, whereas the
published B tree gives a difference of two product orbit sums. We evaluate
that difference by a cubic limit of the nonterminating
$q$-Saalsch\"utz summation. Together with the A--C bijection of
Callan--Mansour~\cite[Section~4.10]{CM}, this proves the Class~207 case
of their Conjecture~1.

Let $q=\lambda-1$, where $\lambda=(1-\sqrt{1-4z})/(2z)$ is the Catalan series, so $z=q/(1+q)^2$, and set
\[
 L=1-2q-2q^2+q^3,\qquad q+Lr+q^2r^2=0,\qquad r=-q+O(q^2).
\]
This specifies $r\in q\QQ[[q]]$ uniquely; explicitly,
$r=-2q/(L+\sqrt{L^2-4q^3})$, with square root $1+O(q)$.
Define
\begin{equation}
 T(q,r)=\sum_{k\ge1}\frac{k^2(qr)^k}{1-q^k},\qquad
 \theta=\frac{(1-q)^3(1+q)}{q^2},\qquad H=1+\theta T(q,r(q)).
 \label{count:eq:H}
\end{equation}
\begin{theorem}\label{count:thm:main}
For all $n\ge0$, the three classes $A,B,C$ have the same number of
inversion sequences of length $n$. Their common generating function satisfies
\begin{equation}
 G_A\!\left(\frac q{(1+q)^2}\right)
 =G_B\!\left(\frac q{(1+q)^2}\right)
 =G_C\!\left(\frac q{(1+q)^2}\right)
 =\frac{(1+q)^2}{q}+\frac1H.
 \label{count:eq:main}
\end{equation}
Here $H=-q+O(q^2)$; the two Laurent expressions on the right have
canceling negative powers and their sum belongs to $\QQ[[q]]$.
\end{theorem}
The proof first obtains the two counting formulas and then evaluates the
orbit difference that connects them. The summation used in the last step is
classical~\cite[Eq.~(II.24)]{GR}; see also~\cite[Eq.~17.7.5]{DLMF}.
We give the specialization and cubic limit explicitly.

\subsection{The C boundary and its Lambert series}
Use the C statistics of Section~\ref{sec:C} and, in this subsection,
write $P=P_C$, $F=F_C$ and $Z=Z_C$. In particular,
$G_C=1+F(1)$ and $\deg_x[z^n]P\le n$, since the slack is at most
the prefix length. Equations~\eqref{eq:CPexplicit}--\eqref{eq:CZ}
provide the complete system used below.

Write $\kappa(x,y)=(1-y)(x-y)+zx(x-1)$. At a kernel root, clearing
\eqref{eq:CZ} and substituting \eqref{eq:CPexplicit}--\eqref{eq:CFexplicit} gives
\begin{equation}
 (1-y)^2\left(Z(y,y)+\frac{P(y)}{zy}\right)
 =(1-x)(G_C-x)+\frac{(2-y)(1-x)^2}{1-y}P(x).
 \label{count:eq:kernel}
\end{equation}
For two roots $x,x'$, subtraction and $x+x'=1+(y-1)/z$ give
\begin{equation}
 (2-y)\bigl((1-x)^2P(x)-(1-x')^2P(x')\bigr)
 =(1-y)(x-x')\left(G_C+\frac{1-y}{z}\right).
 \label{count:eq:two}
\end{equation}
In particular, for an indeterminate $t$, put
\[
 d=t^2-t+z,\quad X=\frac{t(t-1)}d,\quad
 X'=\frac{z-t}d,\quad Y=\frac{(t-1)(t-z)}d.
\]
They satisfy $\kappa(X,Y)=\kappa(X',Y)=0$ and have nonnegative $z$-orders.
Substitution in \eqref{count:eq:two} gives
\begin{equation}
 (d+tz)\bigl(z^2P(X)-t^4P(X')\bigr)
 =zt(t^2-z)(dG_C+t).
 \label{count:eq:input}
\end{equation}
All substitutions so far are coefficientwise legitimate over $\QQ(t)[[z]]$.

Put $D(u)=qu^2+Lu+q^2$ and define
\begin{align*}
 j_G(u)&=\frac{(u-1)(u-q)(u^2-q)(1-q)^2}{uD(u)},&
 j_0(u)&=\frac{(u-1)^2(u-q)^2(u^2-q)(1+q)}{u^2D(u)},\\
 \cLam(j)&=\sum_{k\ge1}\frac{k^2q^k}{1-q^k}[u^k]j(u),&
 \Lambda_G&=\cLam(j_G),\qquad\Lambda_0=\cLam(j_0).
\end{align*}
Laurent expansions here are $q$-first: expand in $q$ with Laurent-polynomial
coefficients in $u$, and then extract $[u^k]$.
Take
\[
 t=\frac{u-q}{(1+q)(u-1)},\quad
 d=\frac{(1-q)^2u}{(1+q)^2(u-1)^2},\quad
 X=1-\frac q{(1-q)^2}(u+u^{-1}-2),\quad Q(u)=P(X)/d^2.
\]
Equation~\eqref{count:eq:input} becomes, with $\cGc=G_C(q/(1+q)^2)$ and $\cGb=G_B(q/(1+q)^2)$,
\begin{equation}
 Q(u)-Q(q/u)=\frac{(1+q)^4}{q(1-q)^3}
                 (j_G(u)\cGc+j_0(u)).\label{count:eq:add}
\end{equation}
The substitution deserves an order check. Writing $w=u+u^{-1}-2$,
$Q$ is a unit in $\QQ[[q]]$ times $w^2P(z,1-qw/(1-q)^2)$.
Since $\deg_x[z^n]P\le n$, a term $w^{j+2}$ from $P$ has $q$-order at
least $2j$. Consequently, if $Q_k=[u^k]Q$, then
$\operatorname{ord}_q Q_k\ge2|k|-4$ for $|k|\ge3$.
This proves that $Q(q/u)$ and the coefficient sums below are well defined
in the $q$-adic Laurent field. Also $Q(u)=Q(1/u)$ and $Q''(1)=0$, since
$d^{-2}$ supplies a fourth-order zero at $u=1$.
Taking positive coefficients in \eqref{count:eq:add} gives
\[
 Q_k=-\frac{(1+q)^4}{q(1-q)^3}\frac{q^k}{1-q^k}
           [u^k](j_G\cGc+j_0).
\]
Thus $0=Q''(1)=2\sum_{k>0}k^2Q_k$ implies
\begin{equation}
 \cGc=-\Lambda_0/\Lambda_G.\label{count:eq:Cquot}
\end{equation}
To evaluate it, use the two rational identities
\begin{align}
 j_0+\frac{(1+q)^2}{q}j_G
 &=\frac{1+q}{q}\frac{(u-1)(u-q)(u^2-q)}{u^2},\label{count:eq:red}\\
 j_G&=\frac{(1-q)^2}{q}\left(u-\frac q u\right)
 +\frac{(1-q)^4(1+q)}{q^2}
       \left(\frac1{1-ru}-\frac1{1-qr/u}\right).\label{count:eq:partial}
\end{align}
The second follows by $D(u)=-(q/r)(1-ru)(u-qr)$ and the quadratic equation
for $r$. Applying $\cLam$ yields
\[
 \Lambda_0+\frac{(1+q)^2}{q}\Lambda_G=-(1-q),\qquad
 \Lambda_G=(1-q)H.
\]
The first coefficients are $H=-q+q^2-q^3+O(q^4)$, so division is
legitimate in $\QQ((q))$. Therefore
\begin{equation}
 \cGc=\frac{(1+q)^2}{q}+\frac1H.\label{count:eq:Cfinal}
\end{equation}

\subsection{The published B tree and its orbit sums}
We use the tree for $B$ in Callan--Mansour~\cite[Table~1, p.~28]{CM}.
Let $a_m(z),b_{m,j}(z),c_{m,j}(z)$ count nodes of its three types by length,
with root $a_1$ at length one and $1\le j\le m$.
Writing $m'=m+1-j$, the rules are
\begin{align*}
 a_m&\longrightarrow a_{m+1},\ b_{m,1},\ldots,b_{m,m},\\
 b_{m,j}&\longrightarrow c_{m,j},c_{m-1,j-1},\ldots,c_{m',1},\
                    a_{m'+1},b_{m',1},\ldots,b_{m',m'},\\
 c_{m,j}&\longrightarrow c_{m+1,j},c_{m,j-1},\ldots,c_{m'+1,1},\
                    a_{m'+1},b_{m',1},\ldots,b_{m',m'}.
\end{align*}
There is exactly one $a$-child per node, and the empty object precedes
the root. Hence $\sum_m a_m=zG_B$. Define a separate boundary $f$ by
\[
 a(x)=\sum_{m\ge1}a_m x^{m-1}=z(1+f(x)),\quad
 \mathsf B(x,y)=\sum_{m,j}b_{m,j}x^{m-j}y^{j-1},\quad
 \mathsf C(x,y)=\sum_{m,j}c_{m,j}x^{m-j}y^{j-1}.
\]
In particular $f(1)=G_B-1$. This $f$ is a boundary of the published tree;
it is different from both slack series $F_B$ and $F_C$ of
Proposition~\ref{prop:common}.
The rules give $b_{m,j}=a_{m+1}$, since their parent sets coincide, and
therefore
\begin{align}
 \mathsf B(x,y)&=z\frac{f(x)-f(y)}{x-y},\notag\\
 (1-y+zxy)\mathsf C(x,y)
 &=z(\mathsf B(x,1)-y\mathsf B(x,y))+zx\mathsf C(x,1),\notag\\
 a(x)&=z+zx a(x)+zx(\mathsf B(x,1)+\mathsf C(x,1)).\label{count:eq:Btree}
\end{align}
At $y=(1-zx)^{-1}$, eliminate $\mathsf C(x,1)$ to obtain
$(1-zx)f(x)-zx=(x-1)\mathsf B(x,1)+y\mathsf B(x,y)$.
Inserting $\mathsf B$ and clearing denominators proves
\begin{equation}
 zf\!\left(\frac1{1-zx}\right)
 =(1-zx)(1-x+zx+zx^2)f(x)
       +z(1-x+zx^2)(f(1)-x).\label{count:eq:Bfun}
\end{equation}
Thus the scalar equation determines the total of the B class through
$f(1)=G_B-1$, without identifying $f$ with the slack series $F_B$.

Put
\[
 x(v)=\frac{(1+q)(1-v)}{1-qv},\quad E(v)=q+Lv+q^2v^2,\quad
 \mathfrak A(v)=\frac{(1-q^2v)E(v)}{q(1-qv)^3},\quad
 \mathfrak b(v)=\frac{(1-q)^2v}{(1-qv)^2}.
\]
For $\mathcal F(v)=f(x(v))$ and $g=G_B-1$, equation~\eqref{count:eq:Bfun} is
\begin{equation}
 \mathcal F(qv)=\mathfrak A(v)\mathcal F(v)+\mathfrak b(v)(g-x(v)).\label{count:eq:Bdiff}
\end{equation}
With $(a;q)_\infty=\prod_{j\ge0}(1-aq^j)$ and products of several
arguments denoting products of these factors, define
\begin{equation}
 U(v)=\frac{(qv/r,q^2rv;q)_\infty}{(qv;q)_\infty^2},\quad
 W(v)=\frac{(1-q)^2v}{1-q^2v}U(v),\quad
 J(v)=\frac{(v/r,qrv;q)_\infty}{(1-qv)(qv;q)_\infty^2}.\label{count:eq:products}
\end{equation}
Set
\[
 M_a=\sum_{k\ge0}W(aq^k),\quad N_a=\sum_{k\ge0}W(aq^k)x(aq^k),\quad
 \mathcal D=M_1-M_r,\quad\mathcal N=N_1-N_r.
\]
Cancellation of product factors, using $E(r)=0$, gives
\[
 J(qv)\mathfrak A(v)=J(v),\qquad J(qv)\mathfrak b(v)=W(v).
\]
Multiplying \eqref{count:eq:Bdiff} by $J(qv)$ and summing on an orbit therefore
telescopes. On the orbit from $1$, use $\mathcal F(1)=f(0)=0$;
on the orbit from $r$, use $J(r)=0$ without dividing by $\mathfrak A(r)=0$.
Since type $a_m$ first occurs at length $m$, $\deg_x[z^n]f\le n$.
Thus $\mathcal F(aq^k)\to\mathcal F(0)$ and $J(aq^k)\to1$ $q$-adically.
Since $J(0)=1$, the result in either case is
$\mathcal F(0)=gM_a-N_a$. Thus
\begin{equation}
 \cGb=1+\frac{\mathcal N}{\mathcal D}.\label{count:eq:Bquot}
\end{equation}
The formal limits exist: $W(q^k)=O(q^k)$,
$W(rq^k)=O(q^{k+1})$, and $\mathcal D=2+O(q)$.

For $K=(1+q+q^2)/q$, finite product cancellation also gives
\begin{equation}
 J(qv)-J(v)=W(v)(x(v)-K).\label{count:eq:tel}
\end{equation}
Indeed, after division by $U(v)$, the required polynomial identity is
\[
 (1-qv)^3-(1-q^2v)(1-v/r)(1-qrv)
 =-\frac{(1-q)^2}{q}v(1-q^3v).
\]
Since $J(0)=1$ and $J(r)=0$, summing \eqref{count:eq:tel} yields
$\mathcal N-K\mathcal D=-J(1)$. Equation~\eqref{count:eq:Bquot} becomes
\begin{equation}
 \cGb=\frac{(1+q)^2}{q}-\frac{J(1)}{\mathcal D}.\label{count:eq:Bfinal}
\end{equation}
It remains to prove $\mathcal D=-J(1)H$.

\subsection{A two-parameter summation}
\begin{lemma}\label{count:lem:sum}
For independent real parameters $0<q<1$ and $-1<r<0$, define
$U,W,J,M$ by \eqref{count:eq:products} and the orbit sums above. Then
\begin{equation}
 M_1-M_r=-J(1)(1-\beta T(q,r)),\qquad
 \beta=\frac{(1-q)(1-r)(1-qr)}{qr}.\label{count:eq:sum}
\end{equation}
The quadratic relation between $q$ and $r$ is not assumed here.
\end{lemma}
\begin{proof}
For $(a;q)_n=\prod_{j=0}^{n-1}(1-aq^j)$, use the convention
\[
 {}_3\phi_2\!\left(\begin{matrix}a,b,c\\ e,f\end{matrix};q,q\right)
 =\sum_{n\ge0}\frac{(a,b,c;q)_n}{(q,e,f;q)_n}q^n.
\]
The nonterminating $q$-Saalsch\"utz identity~\cite[Eq.~(II.24)]{GR},
in the form of~\cite[Eq.~17.7.5]{DLMF}, holds for $0<|q|<1$ and
$ef=abcq$, with denominator factors nonzero. It states
\begin{align}
 {}_3\phi_2\!\left(\begin{matrix}a,b,c\\ e,f\end{matrix};q,q\right)
 &+\frac{(q/e,a,b,c,qf/e;q)_\infty}
              {(e/q,aq/e,bq/e,cq/e,f;q)_\infty}
 {}_3\phi_2\!\left(\begin{matrix}aq/e,bq/e,cq/e\\ q^2/e,qf/e\end{matrix};q,q\right)\notag\\
 &=\frac{(q/e,f/a,f/b,f/c;q)_\infty}
               {(aq/e,bq/e,cq/e,f;q)_\infty}.\label{count:eq:saal}
\end{align}
Set
\[
 a=b=qr,\quad c=q^2r,\quad e=q^3rs,\quad f=q^2r^2/s,
 \qquad h=1-1/s,\quad t=1-h.
\]
The balance is exact. Denote the first series by $\Phi(s)$, the multiplier
of the second by $\mathcal A(s)$, and its terms by $b_n(s)$; explicitly,
\begin{align*}
 b_n(s)&=\frac{(t/q,t/q,t;q)_n}{(q,t/(qr),rt^2;q)_n}q^n,\\
 \mathcal A(s)&=
 \frac{(t/(q^2r),qr,qr,q^2r,rt^2;q)_\infty}
      {(q^2r/t,t/q,t/q,t,q^2r^2t;q)_\infty}.
\end{align*}
At $s=1$, the first series is $\Phi(1)=M_r/W(r)$, as follows by comparing
the $n=0$ terms and successive ratios. The multiplier has a triple pole:
\[
 (t/q;q)_\infty^2(t;q)_\infty=(1-t/q)^2h^3(qt;q)_\infty^3.
\]
Write $A_{-3}=\lim_{h\to0}h^3\mathcal A(s)$. Finite shifts give
\begin{align}
 A_{-3}&=\frac{q^2(q^{-2}/r;q)_\infty(qr;q)_\infty^2(r;q)_\infty}
 {(1-q)^2(q^2r^2;q)_\infty(q;q)_\infty^3},&
 W(r)A_{-3}&=\beta J(1).\label{count:eq:Aminus}
\end{align}
For clarity, the last cancellation reduces to
\[
 \frac{W(r)A_{-3}}{J(1)}
 =\frac{q^2r(1-q)(1-r)}{1-q^2r}
        (1-q^{-2}/r)(1-q^{-1}/r)=\beta.
\]
For every $n\ge2$,
\begin{equation}
 \lim_{h\to0}\mathcal A(s)b_n(s)
 =-\frac{W(q^{n-2})}{W(r)}.\label{count:eq:tail}
\end{equation}
A base case and ratio prove this for all indices. At $n=2$,
\[
 \lim_{h\to0}\frac{b_2}{h^3}
 =\frac{(1-q^{-1})^2q^2}
 {(1-q^2)(1-(qr)^{-1})(1-r^{-1})(1-r)(1-qr)};
\]
multiplication by $W(r)A_{-3}/W(1)$ gives $-1$ by finite product shifts.
The common successive ratio for \eqref{count:eq:tail} is
\[
 q\frac{(1-q^{n-1})^2(1-q^n)}
 {(1-q^{n+1})(1-q^{n-1}/r)(1-rq^n)}.
\]
In particular the $n=2$ term retains the orbit endpoint $W(1)$.

Divide the right side of \eqref{count:eq:saal} by $\mathcal A(s)$ and call it
\[
 P(h)=\frac{(qrt;q)_\infty^2(rt;q)_\infty(q^2r/t;q)_\infty}
 {(qr;q)_\infty^2(q^2r;q)_\infty(rt^2;q)_\infty}.
\]
The two exceptional terms are $b_0=1$ and
\[
 b_1=-\frac{rh(1-q-h)^2}
 {(1-q)(1-qr-h)(1-r(1-h)^2)}.
\]
Let
\[
 R(h)=\frac{(1-r+rh)(1-qr)(1-h)}
 {(1-qr-h)(1-r(1-h)^2)}.
\]
The following two identities isolate the finite part:
\begin{align}
 R-1-b_1&=\frac{rqh^3}{(1-q)(1-qr-h)(1-r(1-h)^2)},\label{count:eq:cubic}\\
 \frac{P}{R}&=\frac{(qrt;q)_\infty^3(qr/t;q)_\infty}
 {(qr;q)_\infty^3(qrt^2;q)_\infty}.\label{count:eq:PR}
\end{align}
The first is a polynomial identity after clearing denominators; the second
uses four single-step product shifts. Put $a=qr$, $\lambda=\log t$ and
$L_j(a)=(a\partial_a)^j\log(a;q)_\infty$. Taking logarithms of
\eqref{count:eq:PR} gives
\[
 \log(P/R)=\sum_{j\ge1}\frac{3+(-1)^j-2^j}{j!}L_j(a)\lambda^j.
\]
The linear and quadratic coefficients vanish, the cubic coefficient is
$-L_3(a)$, and
\[
 L_3(a)=-\sum_{k\ge1}\frac{k^2a^k}{1-q^k}=-T(q,r),\qquad
 \lambda^3=-h^3+O(h^4).
\]
Consequently $P-R=-Th^3+O(h^4)$, and \eqref{count:eq:cubic} proves
\begin{equation}
 [h^3](P-1-b_1)=\frac1\beta-T.\label{count:eq:finitepart}
\end{equation}
Now \eqref{count:eq:saal} reads
$\Phi(s)=\mathcal A(s)(P-1-b_1-\sum_{n\ge2}b_n)$.
Using \eqref{count:eq:tail} and \eqref{count:eq:finitepart} gives
\[
 \frac{M_r}{W(r)}=A_{-3}(1/\beta-T)+\frac{M_1}{W(r)}.
\]
Equation~\eqref{count:eq:Aminus} proves \eqref{count:eq:sum}.

To justify the limit through the infinite tail, for $n\ge2$ extract the
three zeros explicitly:
\[
 \frac{b_n}{h^3}=
 \frac{(1-t/q)^2(qt;q)_{n-2}^2(qt;q)_{n-1}q^n}
 {(q;q)_n(t/(qr);q)_n(rt^2;q)_n}.
\]
For fixed $0<q<1$, $-1<r<0$, on a sufficiently small closed complex
neighborhood of $h=0$ the denominator products are bounded away from
zero uniformly in $n$, and the numerator products are bounded. Therefore
$|b_n/h^3|\le Cq^n$. The two ${}_3\phi_2$ series have argument $q$,
so they converge, and this bound justifies termwise passage to the limit.
All product shifts and logarithmic derivatives above are then valid there.
\end{proof}

\subsection{Formal specialization and completion of the proof}
On the branch $q+Lr+q^2r^2=0$, elementary algebra gives
\[
 \beta=-\theta,\qquad J(1)=\frac{1-q^2}{q}U(1).
\]
The implicit function theorem makes $r(q)$ analytic near zero; moreover,
$q/r=-L-q^2r$ removes the apparent singularity in $U$.
The functions $U(1)$, $M_1$, $M_r$ and $qJ(1)$ are analytic germs,
since their nonconstant factors and orbit tails converge uniformly on
sufficiently small disks. In $qJ(1)$ the prefactor $q$ removes its simple pole.
In particular,
\[
 \mathcal D(0)=2,\qquad (qJ(1))(0)=2.
\]
For sufficiently small positive $q$, the branch has $-1<r(q)<0$.
Lemma~\ref{count:lem:sum} thus gives
\[
 q\mathcal D+(qJ(1))H=0
\]
on a real interval. Both sides are analytic at zero, so the identity theorem
proves the corresponding formal identity. It gives $H=-q+O(q^2)$ and
$\mathcal D=-J(1)H$. Equations~\eqref{count:eq:Cfinal} and \eqref{count:eq:Bfinal}
coincide. Since $z=q+O(q^2)$ is formally invertible, $G_B=G_C$ as length
series. The published bijection $G_A=G_C$ completes Theorem~\ref{count:thm:main}.

\section{A common counting formula for Class 152}\label{sec:class152}
An inversion sequence of length $n$ is $e=(e_0,\ldots,e_{n-1})$ with
$0\le e_i\le i$. A pattern occurs on three increasing indices when both
inequalities and equalities agree with the pattern. Write $\In_n(B_{152})$ for the
sequences avoiding every pattern in $B_{152}$, and put
\[
 A_{152}=\{010,100,102,210\},\qquad B_{152}=\{011,201,210\}.
\]
\begin{theorem}\label{count:class152}
For every $n\ge0$, $|\In_n(A_{152})|=|\In_n(B_{152})|$. The common count is $1$ at $n=0$;
for every $n\ge1$ it is
\[
 a_n=1+\frac12\sum_{j=1}^{n-1}\binom{2j}{j}.
\]
\end{theorem}
For $n\ge1$, the number $a_n$ is \textup{OEIS A024718}$(n-1)$~\cite{OEIS}.
That entry already records the partial central-binomial sum and its generating
function
\[
 \sum_{k\ge0}\textup{A024718}(k)z^k
 =\frac{1+(1-4z)^{-1/2}}{2(1-z)},
\]
the latter attributed to Michael Somos (14 February 2006). The new
contribution is the enumeration of the two specified inversion-sequence
avoidance classes and the proof of their Wilf equivalence.
This is Class 152 of Conjecture 1 in Callan--Mansour, \emph{Inversion sequences avoiding quadruple length-3 patterns}~\cite{CM}, Section 1.
We use ordinary length $n$ throughout, including the empty sequence. The triple-avoidance class $B_{152}$ appears as class~73 in Table~7 of Callan, Jel\'inek and Mansour~\cite{CJM}, where its initial counting sequence is recorded. The theorem above gives its closed counting formula and its equality with the quadruple-avoidance class $A_{152}$.

\subsection{The first-descent decomposition}
If a value $a$ follows a larger value, any second occurrence of $a$ gives
$010$ or $100$. Thus $a$ occurs exactly once. If an adjacent descent $M>a$
has a larger value elsewhere, that value forms $210$ when it is earlier and
$102$ when it is later. Consequently every adjacent descent starts at the
global maximum.

A nonmonotone avoider has a unique cut $p\,s$ at its first adjacent descent:
$p$ is a nondecreasing inversion sequence ending in the global maximum $M$,
and $s$ begins below $M$. Its entries other than $M$ are distinct and increasing
(a decrease between two of them gives $210$ with the preceding $M$).
They avoid all values of $p$. Moreover they lie in one gap $(u,v)$ between
successive distinct values of $p$: crossing a prefix value $v$ gives the triple
$v,a,b$ of type $102$. If $v<M$, no $M$ can occur in $s$, since $v,a,M$
would give $102$. Conversely these conditions exclude all four forbidden
patterns. The bounds $e_i\le i$ are automatic in $s$, since its entries are
at most $M\le |p|-1$.

Let $m=|p|$, $r=|s|\ge1$, and let $g_1,\ldots,g_k$ be the positive jumps
between consecutive distinct values of $p$ (starting from $0$).
An internal gap of size $g$ contributes $\binom{g-1}{r}$.
The top gap contributes
\begin{equation}
 \sum_{b\ge1}\binom{g_k-1}{b}\binom{r-1}{b-1}
 =\binom{g_k+r-2}{r}.
\label{c152:gap}
\end{equation}
Here $b$ is the number of smaller entries and the second factor chooses their
positions, with the first position required to be smaller. Empty sums contribute
zero; a gap of size $1$ contributes zero.

Encode $p$ as a Dyck path by requiring exactly $p_i$ down-steps in total before its $(i+1)$st
up-step and completing with down-steps. Nondecreasing inversion sequences
of length $m$ are thereby in bijection with Dyck paths of semilength $m$.
Their positive jumps are precisely the lengths of the nonfinal down-runs;
the top jump is the penultimate down-run. This translates \eqref{c152:gap} into the
weighted path series below.

\subsection{Weighted Dyck paths and the left generating function}
Work in formal power series over $\mathbb Q$. Set
\begin{equation}
 \mathcal C=\sum_{m\ge0}\frac1{m+1}\binom{2m}{m}z^m,
 \qquad t=z\mathcal C,\qquad \mathcal C=1+z\mathcal C^2=\frac1{1-t},\qquad z=t(1-t).
\label{c152:cat}
\end{equation}
For a weight $w$ with $w(0)=0$, let $D_w$ count paths weighted by their final
down-run, $H_w$ count the sum of weights of all nonfinal down-runs, and
$J_w$ count the weight of the penultimate down-run (zero if absent).
All series include the factor $z^m$ for semilength $m$. Then
\begin{equation}
 (2-\mathcal C)H_w=(\mathcal C-1)D_w,\qquad (1-z)J_w=tD_w.
\label{c152:runs}
\end{equation}
Here is a derivation that also checks the marked-run multiplicities.
Reflect a Dyck path to exchange final down-runs and initial up-runs.
A path whose initial up-run has length $k$ has a unique factorization
$U^k D P_1 D P_2\cdots D P_k$, with $k$ arbitrary Dyck paths $P_i$.
Thus paths with that initial run contribute $t^k$. Let $S_w$ count the sum
of weights of \emph{all} up-runs. The factorization gives
\[
 S_w=D_w+S_w\sum_{k\ge1}k z^k \mathcal C^{k-1}
     =D_w+t\mathcal C S_w.
\]
Since $t\mathcal C=\mathcal C-1$ and reflection gives $H_w=S_w-D_w$, the first identity in \eqref{c152:runs}
follows. For the second up-run, select the first nonempty $P_i$; factors
before it are empty, and factors after it are arbitrary. Hence
\[
 J_w=D_w\sum_{k\ge1}z^k\sum_{j=0}^{k-1}\mathcal C^j
     =\frac{z\mathcal C}{1-z}D_w.
\]
Reflection gives the asserted identity for the penultimate down-run.

For $q$ of constant term $1$, use $w_q(\ell)=q^{\ell-1}-1$ for $\ell\ge1$
and $w_q(0)=0$, and write $D_q,H_q,J_q$. The same run factorization yields
\begin{equation}
 D_q=\sum_{\ell\ge1}(q^{\ell-1}-1)t^\ell
     =\frac{\mathcal C t^2(q-1)}{1-qt}.
\label{c152:D}
\end{equation}
With $q_1=1+z$ and $q_2=(1-z)^{-1}$, binomial expansion gives
$[z^r]w_{q_1}(g)=\binom{g-1}{r}$ and
$[z^r]w_{q_2}(g)=\binom{g+r-2}{r}$ for $r\ge1$.
The first-descent decomposition therefore gives the left series
\begin{equation}
 L=\sum_{n\ge0}|\In_n(A_{152})|z^n=\mathcal C+H_{q_1}+J_{q_2}-J_{q_1}.
\label{c152:L}
\end{equation}
In \eqref{c152:L}, $H_{q_1}$ initially weights every nonfinal down-run as an internal
gap; the difference of the two $J$ terms replaces the top-gap weight.
The path's power of $z$ records prefix length, and the weight's power records
suffix length.

Equations \eqref{c152:cat}--\eqref{c152:D} give, directly,
\[
 D_{q_1}=\frac{t^3}{(1-t)(1-t^2)},\qquad
 D_{q_2}=\frac{t^3}{(1-t)^2}=(1+t)D_{q_1},
\]
\begin{equation}
 H_{q_1}+J_{q_2}-J_{q_1}
 =D_{q_1}\left(\frac{t}{1-2t}+\frac{t^2}{1-z}\right)
 =\frac{\mathcal C t^4}{(1-z)(1-2t)}.
\label{c152:simplify}
\end{equation}
All denominators used here have nonzero constant term, so these operations
are valid in formal power series.

\subsection{The last-zero decomposition and the right generating function}
Avoiding $011$ means that all positive entries are distinct, since $e_0=0$.
Cut a nonempty avoider uniquely as $p\,0\,s$ at its last zero, with $m=|p|$
and $M=\max(p)$ (take $M=0$ for empty $p$). Positive entries of $p$ are
increasing: a positive descent followed by the last zero gives $210$.
Every entry of $s$ exceeds $M$: equality is excluded by $011$, and a positive
entry below $M$ gives $M,0,a$ of type $201$. Thus $s$ consists of distinct
entries above $M$, avoids $201,210$, and satisfies $s_j\le m+1+j$.
Conversely these conditions suffice: a forbidden triple crossing the cut
would contradict either the increasing positive prefix or the separation
of prefix and suffix values.

The running-maximum map sends $p$ bijectively to a nondecreasing inversion
sequence: the inverse retains strict increases and replaces other entries
by zero (retaining the first entry $0$). It preserves $M$. Its Dyck path has final down-run length
$h=m-M$; the empty prefix corresponds to the empty path with $h=0$.

A distinct-entry word avoiding $201,210$ has a unique decomposition into
rotated increasing blocks: sort its entries, partition into nonempty
consecutive blocks, and move the minimum of each block to its end.
Indeed, at the position of the global minimum, the preceding entries must
increase (otherwise $210$), and must be below every following entry
(otherwise $201$). Remove this first block and repeat; this also proves
uniqueness and the converse. For example, sorted blocks
$(3,4)\mid(6)$ become the suffix $4,3,6$: the entry $4$ moves one
place to the left and $3$ moves to the end of its block, while the
singleton $6$ stays in place.

For the sorted entries $a_i$, put $\lambda_i=a_i-M-(i+1)$, with $i$ starting
at $0$. These labels are weakly increasing and nonnegative. The bounds on
the rotated word are exactly $0\le\lambda_i\le h$, with labels in block
tails strictly below $h$. Here the \emph{tail} of a block consists of
all entries other than its minimum in the unrotated increasing block;
each such entry moves one place earlier, whereas the minimum moves
to the end of its block. A singleton can have label $h$;
in a longer block its minimum is bounded by its tail labels. For the example above, take $M=1$, $m=3$, so $h=2$.
The labels of $3,4,6$ are $1,1,2$: the tail label $1$ is strictly
below $h$, and the singleton label is $h$. Adding back $M+i+1$
reverses this encoding.

For each label $0,\ldots,h-1$, an occurrence may start a block or continue
one, giving two choices; the first occurrence overall must start a block.
Every occurrence of $h$ must start a singleton. Words containing a lower
label therefore have half the unrestricted binary choices, while words
containing only $h$ have one choice. The suffix series (including the empty
suffix) is consequently
\begin{equation}
 S_h=\frac1{1-z}+\frac{(1-2z)^{-h}-1}{2(1-z)}
     =\frac{1+(1-2z)^{-h}}{2(1-z)}.
\label{c152:suffix}
\end{equation}
The last zero contributes $z$. The final-run factorization in the preceding subsection,
now with weight $r^h$, gives $\sum_P z^{|P|}r^{h(P)}=(1-rt)^{-1}$.
Using $r=(1-2z)^{-1}$ in \eqref{c152:suffix}, we obtain
\begin{equation}
 R=\sum_{n\ge0}|\In_n(B_{152})|z^n
 =1+\frac{z}{2(1-z)}\left(\mathcal C+\frac1{1-t/(1-2z)}\right)
 =1+\frac{z(1-t)}{(1-z)(1-2t)}=L,
\label{c152:R}
\end{equation}
where the last two simplifications use $z=t(1-t)$ and \eqref{c152:simplify}.
Finally $1-2t=\sqrt{1-4z}$, so
\[
 L=R=1+\frac{z}{2(1-z)}\left(1+\frac1{\sqrt{1-4z}}\right).
\]
Extracting coefficients proves the theorem. For $n=0,\ldots,10$, the counts are
\[1,1,2,5,15,50,176,638,2354,8789,33099.\]

\paragraph{Formal verification record}
The formal record uses Lean~4~\cite{Lean} and mathlib~\cite{Mathlib}.
The source theorem \texttt{InversionSeqClass152.result} states the all-length
equality of the finite avoider counts. Its pinned source is
\href{https://github.com/the-omega-institute/trureturing/blob/7e22ca771362e94db46dd4394637c98db0f942f1/D5/S3/Combinatorics/InversionSeq/InversionSeqClass152.lean}{\texttt{InversionSeqClass152.lean}}.
The closed coefficient formula above follows by the displayed generating-function
calculation. The mathematical proof here is complete independently of
this verification record. The supplementary material gives the pinned
Lean and mathlib versions and a command to check this source theorem.

\section{Conclusion}
The two enumeration problems have explicit solutions.  For Class~152, the
common count is a partial central-binomial sum.  For Class~207,
Theorem~\ref{count:thm:main} proves equality of all three counting series and
gives their common Lambert-series expression.  The C kernel and the published
B tree lead to different scalar formulas; the cubic limit of the
nonterminating $q$-Saalsch\"utz summation is the identity that connects them,
while the A--C step is the published bijection of Callan--Mansour.

The expanded supplement preserves the fixed-maximum rational expansion
(Section~\ref{supp-sec:fixed}), finite computations (Section~\ref{supp-sec:checks}),
refined generating functions and analytic obstruction theorems
(Sections~\ref{supp-ana:sec:obstructions}--\ref{supp-exp:sec:exponential}),
the signed catalytic identity (Section~\ref{supp-sec:signed}), and the Python
verification program (\ref{supp-sec:program}). The pole calculations
there are related to the orbit-residue viewpoint of rational
$q$-summability in Chen and Singer~\cite[Proposition~2.10]{CS}.
Those results provide additional structure and reproducibility without
interrupting the two complete proof paths in the main text.

\section*{Declarations}
\paragraph{Funding}
This research received no specific grant from any funding agency in the public, commercial, or not-for-profit sectors.

\paragraph{Declaration of competing interest}
The authors declare that they have no competing interests.

\paragraph{Data availability}
The data and verification code supporting the finite computations are
available with reproduction instructions in the
\href{https://github.com/the-omega-institute/inversion-sequence-classes-152-207}{companion code and data repository}.
The accompanying supplementary PDF contains the expanded proofs and the
recurrence program; its Section~\ref{supp-sec:checks} describes the finite checks.
The ordinary proofs are contained in this article.

\paragraph{Declaration of generative AI and AI-assisted technologies in the manuscript preparation process}
Haobo Ma and Wenlin Zhang used AI tools, namely ChatGPT Pro (GPT-6 Pro), Codex (gpt-6.1-sol) and Claude Opus 5.5 (via Claude Code), for literature search, computation and double-checking, including formal verification of the results in Lean. All proofs can be checked by hand, and the computations can be reproduced from the verification package. H. Ma and W. Zhang take full responsibility for their contributions.

\appendix
\section{Fixed-maximum expansion for 207(2)}\label{sec:fixed-main}
Let $c_{n,m}$ count the $C$-avoiding inversion sequences of maximum
exactly $m$, and write $R_m(z)=\sum_{n\ge1}c_{n,m}z^n$.
Fix the alphabet $\{0,\ldots,m\}$. A state $q=(M,S,D)$ records the
maximum, the seen values still admissible, and the forbidden values.
Start at $q_0=(0,\{0\},\varnothing)$. For $x\notin D$, set
\begin{equation}\label{eq:finiteupdate}
D'=\begin{cases}
D,&x>M,\\
D\cup\{0,\ldots,M-1\},&x=M,\\
D\cup\{0,\ldots,x\}\cup\{a\in S:a>x\},&x<M,
\end{cases}
\qquad
q_x=\bigl(\max(M,x),(S\cup\{x\})\setminus D',D'\bigr).
\end{equation}
Let $Q_m$ be the closure of $q_0$ under these children, and define
\[
(T_m)_{q,t}=\#\{x\in\{0,\ldots,m\}\setminus D(q):q_x=t\},
\qquad r_m(q)=\ind_{\{M(q)=m\}}.
\]
Let $p_m(q)$ count legal inversion prefixes of length $m+1$ at $q$:
start at $q_0$ at length one and, at length $n\le m$, allow only
$0\le x\le n$ in \eqref{eq:finiteupdate}. Regard $p_m,r_m$ as column vectors.

\begin{proposition}\label{prop:fixed}
For every $m\ge0$,
\begin{align}
\det(I-zT_m)&=(1-z)^{m+1},\notag\\
R_m(z)&=z^{m+1}p_m^{\mathsf T}(I-zT_m)^{-1}r_m
       =\frac{P_m(z)}{(1-z)^{m+1}},\label{eq:Rm}\\
G_C(z)&=1+\sum_{m\ge0}\frac{P_m(z)}{(1-z)^{m+1}}.\label{eq:expansion}
\end{align}
Here $P_m(z)=z^{m+1}p_m^{\mathsf T}\operatorname{adj}(I-zT_m)r_m$
is an integer polynomial, $z^{m+1}\mid P_m$, $\deg P_m\le4m+2$,
and $P_m(1)=1$. For fixed $m$, $c_{n,m}$ is polynomial in $n$ of
degree $m$ and leading coefficient $1/m!$ for $n\ge3m+2$.
\end{proposition}

\begin{proof}
The pair exclusions used in Section~\ref{sec:C} give
\eqref{eq:finiteupdate}. Seen forbidden values can be forgotten because
their subsequent exclusions are already forbidden. Every non-loop step
increases $M$ or strictly enlarges $D$. The only self-loops, each of
multiplicity one, are
$L_j=(j,\{j\},\{0,\ldots,j-1\})$, $0\le j\le m$.
Each is reachable by introducing and repeating $j$ (with $L_0=q_0$).
Ordering states by $(M,|D|)$ makes the matrix upper triangular and gives
the determinant. After length $m+1$, every value in the fixed alphabet
satisfies the inversion bound, so subsequent walks give \eqref{eq:Rm}.

Remove self-loops from a walk. Its remaining path has at most $m$ steps
increasing $M$ and at most $m$ other steps, since each of the latter
forbids a new value below $m$. Its contribution has the form
$z^{m+1+d}/(1-z)^\ell$, with $d\le2m$ and $\ell\le m+1$,
including prefix multiplicity. Bringing this finite sum to the common
denominator proves the degree bound. Exactly one contribution visits
all $m+1$ loops: the prefix is $0^{m+1}$, followed between $L_j$ and
$L_{j+1}$ by $j+1,j+1$. It contributes
$z^{3m+1}/(1-z)^{m+1}$, proving $P_m(1)=1$.
If $P_m(z)=\sum_i u_i z^i$, then for $n\ge3m+2$ the coefficient is the
polynomial $\sum_i u_i\binom{n-i+m}{m}$. Terms with $i>n$ vanish
in its polynomial interpretation by the degree bound. Its leading
coefficient is $P_m(1)/m!$. Finally $z^{m+1}\mid P_m$ makes
\eqref{eq:expansion} coefficientwise convergent.
\end{proof}

The first four summands are
\begin{align*}
R_0(z)&=\frac{z}{1-z},&
R_1(z)&=\frac{z^2+z^3-z^4}{(1-z)^2},\\
R_2(z)&=\frac{2z^3+4z^4-7z^5+z^6+z^7}{(1-z)^3},\\
R_3(z)&=\frac{6z^4+12z^5-34z^6+11z^7+11z^8-4z^9-z^{10}}
                    {(1-z)^4}.
\end{align*}
This is an all-length rational expansion for $G_C$, with a finite
definition of every numerator. Section~\ref{sec:class207-equality} gives
a Lambert series expression for the same generating function.

\section{Finite checks and reproducibility record}\label{sec:checks-main}
The finite checks below are supporting evidence for the displayed formulas;
the all-length conclusions follow from the proofs in the main text.  The
recurrence in Proposition~\ref{prop:recurrence}, an independent
implementation of the four-type $B$ tree, and a literal site-word
implementation reproduce the Class~207 counts through length $12$:
\[
\begin{array}{c|rrrrrrrrrrrrr}
 n&0&1&2&3&4&5&6&7&8&9&10&11&12\\ \hline
207(2)&1&1&2&6&21&80&324&1375&6052&27425&127273&602487&2900432\\
207(3)&1&1&2&6&21&80&324&1375&6052&27425&127273&602487&2900432
\end{array}
\]
Literal pattern checks through length $8$ verify the initial counts for both
Class~152 sets and both Class~207 sets, and verify every legal-site transition
in that range.  The same enumeration checks the first four fixed-maximum
rational functions through length $8$.  Exact-integer expansion of the two
product orbits gives a zero residual for the summation identity through
internal order $136$ (hence through degree $129$ after the stated
normalization), and the common Lambert formula agrees with the C recurrence
and the four-type B tree through length $128$.  These checks are independent
of, and are not premises for, Theorem~\ref{count:thm:main}.

The checks are implemented in Python~\cite{Python}; exact symbolic
certificates use SymPy~\cite{SymPy}, and an optional numerical audit uses
mpmath~\cite{mpmath}. The complete recurrence program, independent generators, literal checker,
exact symbolic certificates, fixed-maximum expansion, and refined analytic
arguments are in the supplementary material.  Its README records the finite
ranges, software versions, commands, and expected outputs.

\end{document}


\begin{frontmatter}
\title{Supplementary material: counting formulas and refined generating functions for inversion-sequence Classes 152 and 207}
\author[chrono]{Haobo Ma}
\ead{auric@aelf.io}
\author[haifa]{Toufik Mansour}
\ead{tmansour@univ.haifa.ac.il}
\author[nus,chrono]{Wenlin Zhang}
\ead{e1327962@u.nus.edu}
\affiliation[chrono]{organization={ChronoAI Pte Ltd},city={Singapore},country={Singapore}}
\affiliation[haifa]{organization={Department of Mathematics, University of Haifa},city={Haifa},country={Israel}}
\affiliation[nus]{organization={National University of Singapore},country={Singapore}}
\begin{abstract}

For the two Class~152 pattern sets, we prove that the common
number at length $n\ge1$ is
\[
 a_n=\frac{1+\displaystyle\sum_{j=0}^{n-1}\binom{2j}{j}}2.
\]
This identifies both avoidance classes with the known sequence
\textup{OEIS A024718}, shifted by one index.
For Class~207, we prove equality of the three total-count series and
 give their common generating function as an explicit Lambert series
 expression after the substitution $z=q/(1+q)^2$. The proof connects
 the C kernel equation to the published B generating tree through a
 cubic limit of the nonterminating $q$-Saalsch\"utz summation.
 We also give complete legal-site recurrences and catalytic systems for
 $C=207(2)$ and $B=207(3)$, a fixed-maximum rational expansion for C,
 and a formal product solution for B.
For each of the pairs $(Z_B,Z_C)$ and $(Z_B,Z_C(y,y))$, the two
series are algebraically independent over $\CC(z,x,y)$.
We also exclude every finite rational linear relation among the
$y$-derivatives of $Z_B$ and the C diagonal.
For each sufficiently small fixed positive length variable, a further
obstruction excludes finite sums of rational prefactors times
polynomial exponentials in a same-argument B input, subject to a
separation condition on the exponents. These results follow from
meromorphic continuation and different pole sets of the refined series.
\end{abstract}
\begin{keyword}
Inversion sequence \sep Pattern avoidance \sep Wilf equivalence \sep
Catalytic generating function \sep Kernel method \sep Algebraic independence
\MSC[2020] 05A05 \sep 05A15 \sep 05A19 \sep 30D30
\end{keyword}
\end{frontmatter}

\section{Introduction}
This expanded supplement accompanies the main article
\emph{Counting formulas and refined generating functions for inversion-sequence
Classes 152 and 207}. The main article contains the complete enumeration proofs
in Sections~\ref{main-sec:class207-equality} and~\ref{main-sec:class152};
its \ref{main-sec:fixed-main} and~\ref{main-sec:checks-main}
give the fixed-maximum expansion and a concise verification record.
The present file retains the expanded arguments, auxiliary results and program.
Section and equation numbers refer to this supplement unless explicitly
identified as references to the main text.

An inversion sequence of length $n$ is a word $e_0\cdots e_{n-1}$ with
$0\le e_i\le i$. Patterns occur in subsequences, with all equalities and
strict inequalities retained. Write $\In_n(X)$ for the sequences avoiding
the pattern set $X$, and $G_X(z)=\sum_{n\ge0}|\In_n(X)|z^n$.
The index $n$ here is the number of entries. In the notation of
Callan--Mansour~\cite{CM}, a word $e_0\cdots e_m$ belongs to $I_m(X)$;
thus their index is one less than its length. All counts and powers
of $z$ in this paper use the length convention above.
The Class~152 sets used later are denoted by $A_{152},B_{152}$, to
keep them distinct from the Class~207 sets $A,B,C$.

The systematic study of length-three patterns in inversion sequences
was initiated by Corteel, Martinez, Savage and Weselcouch~\cite{CMSW}
and Mansour and Shattuck~\cite{MS}. Martinez and Savage~\cite{MSII}
organized a broader classification by triples of relations, and
Lin~\cite{Lin} studied such restrictions in inversion and ascent sequences.
For sets of ordinary patterns, Yan and Lin~\cite{YL} investigated pairs,
Callan, Jel\'inek and Mansour~\cite{CJM} triples, and Callan and
Mansour~\cite{CM} quadruples; the last classification is the source of the
two conjectures addressed here. Related work develops statistic refinements
for restricted inversion sequences~\cite{LK}, consecutive-pattern
enumeration~\cite{AEC}, and vincular Wilf equivalences~\cite{AEV}.
These adjacency-restricted occurrences use different conventions from the
subsequence occurrences considered here.

Generating trees give both structural descriptions and catalytic equations:
Kotsireas, Mansour and Y{\i}ld{\i}r{\i}m~\cite{KMY} give an algorithmic
approach for inversion sequences, while Testart~\cite{Testart} uses trees
growing on the left for further pattern classes. Britt and
Beaton~\cite{BrittBeaton} complete the enumeration of the triples-of-relations
classification using left- and right-growing trees. Within the quadruple
classification, Asinowski and Polley~\cite[Section~4.1]{AP} resolve
Class~166 through rectangulations. Here we resolve Classes~152 and~207;
the Class~207 proof evaluates the published B tree and our C legal-site
recurrence and identifies their total-count series by a classical
$q$-summation.

We study two enumeration questions from the quadruple classification. For Class~152,
we prove a common counting formula for its two avoidance classes and identify
the known numerical sequence that they enumerate
(Section~\ref{sec:class152}). For Class~207, we prove the equality of
 the three counting series and obtain an explicit common generating function
 (Theorem~\ref{count:thm:main}). We also derive refined generating functions
 and establish analytic obstructions to several types of relations between them.

Callan--Mansour's numbering \cite[Table~1, p.~28]{CM} is
\[
\begin{aligned}
207(1)=A&=\{100,101,110,201\},\\
207(2)=C&=\{100,101,110,210\},\\
207(3)=B&=\{101,110,120,210\}.
\end{aligned}
\]
Their bijection \cite[Section~4.10]{CM} gives $G_A=G_C$: in the strict
left-to-right-maximum decomposition, leave the maximum blocks fixed and
refill the intervening boxes from left to right with the largest available
entry below the corresponding maximum. Its inverse sorts the entries in
those boxes increasingly. We prove the remaining equality $G_C(z)=G_B(z)$
by evaluating the C kernel equation as a Lambert series and the published
B tree as two $q$-product orbit sums. A cubic limit of the nonterminating
$q$-Saalsch\"utz summation identifies these expressions
(Section~\ref{sec:class207-equality}). The summation is classical;
the specialization and limiting calculation supply the identity needed
for the counting problem. We use the form in
Gasper--Rahman~\cite[Eq.~(II.24)]{GR}, also recorded in~\cite[Eq.~17.7.5]{DLMF}.

The refined enumerations give additional information beyond the total counts.
We obtain an exact numerical recurrence and a coefficientwise rational
expansion for $G_C$, explicit catalytic systems with complete proofs for
both classes, and a formal-product solution for the full slack series $F_B$.
Although the totals agree, the refined series have different pole structures;
these yield the algebraic and differential independence results below.

\section{An explicit recurrence for 207(2)}\label{sec:C}
For a nonempty $C$-avoiding prefix $w$ of length $n$, let $M=\max w$
and write $s=n-M\ge1$ for the slack, called the bounce statistic in
Asinowski--Polley~\cite[Section~2.2]{AP}. List its admissible next values at most $M$ in increasing
order as a word $p$: write $\Ssite$ for a value already seen and $\U$ for
an unused value. Let $h=1$ if $M$ is admissible, and $h=0$ otherwise.
All values $M+1,\ldots,n$ are admissible. The initial word $0$ has
state $(\Ssite,1,1)$.

We use the labelled generating-tree framework of West~\cite{West}.
Its translation into generating functions follows the general approach of
Banderier et al.~\cite{Banderier}. The complete succession rule is as follows;
each choice is counted once:
\begin{align}
(p,s,h)&\longrightarrow(p\U^{j-1}\Ssite,s-j+1,1),
                   &&1\le j\le s;\label{eq:newmax}\\
(p,s,1)&\longrightarrow(\Ssite,s+1,1),
                   &&\text{repeat }M;\label{eq:repeat}\\
(p,s,h)&\longrightarrow(\U^k,s+1,0),
                   &&\text{select a site below }M.\label{eq:lower}
\end{align}
In the last line, $k$ is the number of unused sites strictly to the right
of the selected site. An equal pair $a,a$ forbids a later entry below $a$
by $110$. A descent $a>b$ forbids entries at most $b$ by $100,210$, and
forbids $a$ by $101$. Thus selecting $x<M$ removes all sites at most $x$
and all seen sites above $x$, leaving precisely the $k$ unused sites.
Repeating $M$ leaves only its own site; a new maximum introduces no new
exclusion. This proves the rule and shows why gaps between listed sites
do not affect it.

For $r\ge0$, $s\ge1$ and $k\ge0$, let $a_r(s,k)$ count continuations
of $r$ further entries from the formal state $(\U^k,s,0)$, and let $b_r(s)$
count them from $(\Ssite,s,1)$. The rules above define these numbers
even when a state is not reached from the root.

\begin{proposition}\label{prop:recurrence}
The initial conditions are $a_0(s,k)=b_0(s)=1$. For $r\ge0$,
\begin{align}
a_{r+1}(s,k)
 &=\sum_{i=0}^{k-1}a_r(s+1,i)\notag\\
 &\quad+\sum_{j=1}^{s}\left(a_r(s-j+1,k+j-1)
                         +b_r(s-j+1)-a_r(s-j+1,0)\right),\label{eq:a}\\
b_{r+1}(s)
 &=b_r(s+1)
   +\sum_{j=1}^{s}\left(a_r(s-j+1,j)
                         +b_r(s-j+1)-a_r(s-j+1,0)\right).\label{eq:b}
\end{align}
Empty sums are zero. In particular,
\begin{equation}\label{eq:counts}
|\In_0(C)|=1,\qquad |\In_n(C)|=b_{n-1}(1)\quad(n\ge1).
\end{equation}
These are recurrences for nonnegative integers; the differences in the
formulas come from an affine reduction of the site word.
\end{proposition}

\begin{proof}
For a site $a$ in $p$, let $k(a)$ count unused sites strictly to its right,
and put $T_p(y)=\sum_{a\in p}y^{k(a)}$. Write $\phi_r(p,s,h)$ for the
number of continuations of length $r$. We first prove that it is affine
in $h$ and the coefficients of $T_p$, for fixed $r,s$.
At $r=0$ it is $1$. For the induction step, lower selections and a
possible repetition of $M$ contribute
\[
\sum_{a\in p}a_r(s+1,k(a))
       +h\bigl(b_r(s+1)-a_r(s+1,0)\bigr).
\]
For a new maximum, the child's site polynomial is
\begin{equation}\label{eq:Ttransform}
T_{p\U^{j-1}\Ssite}(y)
 =y^{j-1}T_p(y)+\frac{1-y^{j-1}}{1-y}+1,
\end{equation}
so its contribution is also affine by induction. Evaluating this affine
function on the empty word, on $\U^k$, and on $\Ssite$ gives
\begin{equation}\label{eq:affine}
\phi_r(p,s,h)=a_r(s,0)+h\bigl(b_r(s)-a_r(s,1)\bigr)
 +\sum_{a\in p}\bigl(a_r(s,k(a)+1)-a_r(s,k(a))\bigr).
\end{equation}
Indeed, $T_{\U^{k+1}}-T_{\U^k}=y^k$ and
$T_{\Ssite}=T_{\U}=1$. Apply \eqref{eq:affine} to the new-maximum
children $\U^{k+j-1}\Ssite$ and $\Ssite\U^{j-1}\Ssite$; their
continuation counts are respectively
$a_r(s',k+j-1)+b_r(s')-a_r(s',0)$ and
$a_r(s',j)+b_r(s')-a_r(s',0)$, where $s'=s-j+1$.
The lower children of $\U^k$ are $\U^0,\ldots,\U^{k-1}$;
the repetition child of $\Ssite$ is $\Ssite$ with slack $s+1$.
This proves \eqref{eq:a}--\eqref{eq:b}. The root has length one,
which proves \eqref{eq:counts}.
\end{proof}

For example, $a_1(s,k)=s+k$ and $b_1(s)=s+1$.
To calculate all counts through length $N$, at depth $r$ it suffices to
retain $a_r(s,k)$ for $s+k\le N-r$ and $b_r(s)$ for $s\le N-r$.
Row and diagonal partial sums evaluate \eqref{eq:a}--\eqref{eq:b} in
$O(N^3)$ integer operations and $O(N^2)$ stored integers; the appendix
gives the complete program.

\subsection*{Catalytic equations and their proof}
For $C$, let the sums below run over all nonempty $C$-avoiding inversion
prefixes. Define
\[
F_C(x)=\sum_w z^{|w|}x^{s(w)},\qquad
P_C(x)=\sum_w z^{|w|}x^{s(w)}h(w),\qquad
H_C(x,y)=\sum_w z^{|w|}x^{s(w)}T_{p(w)}(y).
\]
Throughout this paper, $F_X(x)$ and $P_X(x)$ are univariate in the
catalyst $x$, whereas $H_X(x,y)$ and $Z_X(x,y)$ are bivariate in the
catalysts $x,y$. The length variable $z$ is suppressed in these
arguments, for example $F_X(x)=F_X(z;x)$. A specialization such as
$H_X(x,1)$ or $H_X(y,y)$ always retains both argument positions.
Put
\[
J_C(x)=zx\frac{F_C(x)-F_C(1)}{x-1}.
\]
The succession rule gives
\begin{align}
F_C(x)&=zx+zxH_C(x,1)+J_C(x),\label{eq:CF}\\
P_C(x)&=zx+zxP_C(x)+J_C(x),\label{eq:CP}\\
H_C(x,y)&=zx+zxP_C(x)
 +\frac{zx}{1-y}\bigl(H_C(x,1)-H_C(x,y)\bigr)\notag\\
&\quad+\frac{zx}{x-y}\bigl(H_C(x,y)-H_C(y,y)\bigr)
 +\frac{2-y}{1-y}J_C(x)
 -\frac{zx\bigl(F_C(x)-F_C(y)\bigr)}{(1-y)(x-y)}.\label{eq:CH}
\end{align}
The coefficient ring is $\mathbb Q[x,y][[z]]$; for a univariate
series use $\mathbb Q[x][[z]]$. Apparent singularities in these equations
are removable coefficientwise. In particular, terms with denominator
$1-y$ must be combined before setting $y=1$. All constant terms in $z$
are zero, and
\[
G_C(z)=1+F_C(1)=1+\sum_{r\ge0}b_r(1)z^{r+1}.
\]

\begin{proof}[Proof of \eqref{eq:CF}--\eqref{eq:CH}]
The root contributes $zx$ to each series. For a parent with slack $s$,
the new maxima have slacks $s-j+1$, $1\le j\le s$, and
\begin{equation}\label{eq:sumidentities}
\sum_{j=1}^s x^{s-j+1}=\frac{x(x^s-1)}{x-1},\qquad
\sum_{j=1}^s x^{s-j+1}y^{j-1}
 =\frac{x(x^s-y^s)}{x-y}.
\end{equation}
The first identity gives $J_C(x)$ in $F_C(x)$ and $P_C(x)$.
Every admissible site at most $M$ gives a child with slack $s+1$,
so their total contribution to $F_C(x)$ is $zxH_C(x,1)$.
Among them, only repetition of an admissible $M$ has an admissible
maximum afterward; it contributes $zxP_C(x)$ to $P_C(x)$.
These observations prove \eqref{eq:CF} and \eqref{eq:CP}.

For $H_C(x,y)$, repetition of $M$ leaves the one-site word $\Ssite$,
so contributes $zxP_C(x)$. A lower site $a$ leaves $\U^{k(a)}$,
whose polynomial is $(1-y^{k(a)})/(1-y)$. Summing over all sites gives
$zx(H_C(x,1)-H_C(x,y))/(1-y)$: the maximum site, when present,
has $k(M)=0$ and contributes zero to this expression.
For a new maximum, \eqref{eq:Ttransform} separates the child's
polynomial into the old-site term and the term independent of the
old sites. The old-site term, summed using the second identity in
\eqref{eq:sumidentities}, is
$zx(H_C(x,y)-H_C(y,y))/(x-y)$.
The independent term is $(2-y-y^{j-1})/(1-y)$; its sum is
\[
\frac{2-y}{1-y}J_C(x)
 -\frac{zx(F_C(x)-F_C(y))}{(1-y)(x-y)}.
\]
This accounts for every child exactly once and proves \eqref{eq:CH}.
Both identities in \eqref{eq:sumidentities} are finite geometric sums;
thus each displayed quotient, after its indicated combination, is a
polynomial at each length. Finally the factor $z$ on each child
contribution determines the coefficients recursively from the root.
\end{proof}

\section{The fixed-maximum expansion for 207(2)}\label{sec:fixed}
Let $c_{n,m}$ count the $C$-avoiding inversion sequences of maximum
exactly $m$, and write $R_m(z)=\sum_{n\ge1}c_{n,m}z^n$.
Fix the alphabet $\{0,\ldots,m\}$. A state $q=(M,S,D)$ records the
maximum, the seen values still admissible, and the forbidden values.
Start at $q_0=(0,\{0\},\varnothing)$. For $x\notin D$, set
\begin{equation}\label{eq:finiteupdate}
D'=\begin{cases}
D,&x>M,\\
D\cup\{0,\ldots,M-1\},&x=M,\\
D\cup\{0,\ldots,x\}\cup\{a\in S:a>x\},&x<M,
\end{cases}
\qquad
q_x=\bigl(\max(M,x),(S\cup\{x\})\setminus D',D'\bigr).
\end{equation}
Let $Q_m$ be the closure of $q_0$ under these children, and define
\[
(T_m)_{q,t}=\#\{x\in\{0,\ldots,m\}\setminus D(q):q_x=t\},
\qquad r_m(q)=\ind_{\{M(q)=m\}}.
\]
Let $p_m(q)$ count legal inversion prefixes of length $m+1$ at $q$:
start at $q_0$ at length one and, at length $n\le m$, allow only
$0\le x\le n$ in \eqref{eq:finiteupdate}. Regard $p_m,r_m$ as column vectors.

\begin{proposition}\label{prop:fixed}
For every $m\ge0$,
\begin{align}
\det(I-zT_m)&=(1-z)^{m+1},\notag\\
R_m(z)&=z^{m+1}p_m^{\mathsf T}(I-zT_m)^{-1}r_m
       =\frac{P_m(z)}{(1-z)^{m+1}},\label{eq:Rm}\\
G_C(z)&=1+\sum_{m\ge0}\frac{P_m(z)}{(1-z)^{m+1}}.\label{eq:expansion}
\end{align}
Here $P_m(z)=z^{m+1}p_m^{\mathsf T}\operatorname{adj}(I-zT_m)r_m$
is an integer polynomial, $z^{m+1}\mid P_m$, $\deg P_m\le4m+2$,
and $P_m(1)=1$. For fixed $m$, $c_{n,m}$ is polynomial in $n$ of
degree $m$ and leading coefficient $1/m!$ for $n\ge3m+2$.
\end{proposition}

\begin{proof}
The pair exclusions used in Section~\ref{sec:C} give
\eqref{eq:finiteupdate}. Seen forbidden values can be forgotten because
their subsequent exclusions are already forbidden. Every non-loop step
increases $M$ or strictly enlarges $D$. The only self-loops, each of
multiplicity one, are
$L_j=(j,\{j\},\{0,\ldots,j-1\})$, $0\le j\le m$.
Each is reachable by introducing and repeating $j$ (with $L_0=q_0$).
Ordering states by $(M,|D|)$ makes the matrix upper triangular and gives
the determinant. After length $m+1$, every value in the fixed alphabet
satisfies the inversion bound, so subsequent walks give \eqref{eq:Rm}.

Remove self-loops from a walk. Its remaining path has at most $m$ steps
increasing $M$ and at most $m$ other steps, since each of the latter
forbids a new value below $m$. Its contribution has the form
$z^{m+1+d}/(1-z)^\ell$, with $d\le2m$ and $\ell\le m+1$,
including prefix multiplicity. Bringing this finite sum to the common
denominator proves the degree bound. Exactly one contribution visits
all $m+1$ loops: the prefix is $0^{m+1}$, followed between $L_j$ and
$L_{j+1}$ by $j+1,j+1$. It contributes
$z^{3m+1}/(1-z)^{m+1}$, proving $P_m(1)=1$.
If $P_m(z)=\sum_i u_i z^i$, then for $n\ge3m+2$ the coefficient is the
polynomial $\sum_i u_i\binom{n-i+m}{m}$. Terms with $i>n$ vanish
in its polynomial interpretation by the degree bound. Its leading
coefficient is $P_m(1)/m!$. Finally $z^{m+1}\mid P_m$ makes
\eqref{eq:expansion} coefficientwise convergent.
\end{proof}

The first four summands are
\begin{align*}
R_0(z)&=\frac{z}{1-z},&
R_1(z)&=\frac{z^2+z^3-z^4}{(1-z)^2},\\
R_2(z)&=\frac{2z^3+4z^4-7z^5+z^6+z^7}{(1-z)^3},\\
R_3(z)&=\frac{6z^4+12z^5-34z^6+11z^7+11z^8-4z^9-z^{10}}
                    {(1-z)^4}.
\end{align*}
This is an all-length rational expansion for $G_C$, with a finite
definition of every numerator. Section~\ref{sec:class207-equality} gives
a Lambert series expression for the same generating function.

\section{207(3) and a common catalytic form}\label{sec:B}
Here is a direct four-type generating tree for $B$. Let $M,L$ be the
largest and second largest entries, counting multiplicity; for the
initial word $0$ put $L=0$. Define
\[
D(w)=\{w_i:w_i>w_j\text{ for some }j>i\},\qquad
s=n-M,\quad d=M-L.
\]
The types are $R_s$ if $M=L$, and, when $M>L$,
\[
\begin{array}{c|cc}
\text{type}&M\in D(w)&L\in D(w)\\\hline
U_{s,d}&\text{no}&\text{no}\\
V_{s,d}&\text{yes}&\text{no}\\
W_{s,d}&\text{no}&\text{yes}
\end{array}
\]
The root is $R_1$, and the complete succession rule is
\begin{align}
R_s&\longrightarrow R_{s+1},\quad U_{s-j+1,j}\ (1\le j\le s),\notag\\
U_{s,d}&\longrightarrow R_{s+1},\quad V_{s+1,k}\ (1\le k\le d),
                         \quad U_{s-j+1,j}\ (1\le j\le s),\notag\\
V_{s,d}&\longrightarrow V_{s+1,k}\ (1\le k\le d),
                         \quad W_{s-j+1,j}\ (1\le j\le s),\label{eq:Btree}\\
W_{s,d}&\longrightarrow R_{s+1},\quad V_{s+1,k}\ (1\le k<d),
                         \quad U_{s-j+1,j}\ (1\le j\le s).\notag
\end{align}
To prove the rule, an earlier pair $a\le b$ forbids an appended $x<a$
by $110$ or $120$; a descent $a>b$ forbids $x<b$ by $210$ and $x=a$
by $101$. To identify their intersection, the pair consisting of the two largest
entries forbids every $x<L$, while every earlier pair has its smaller
entry at most $L$ and forbids nothing larger by its inequality condition.
The equality exclusions are exactly $D(w)$. Thus the legal next entries
are $[L,n]\setminus D(w)$ (also for the single-entry root).
If $M>L$, the maximum occurs once. The two endpoints cannot both
belong to $D(w)$: if every $L$ occurs after $M$, a descent from $L$
gives $M,L,v$ with $v<L$, a $210$. If some $L$ occurs before $M$,
a descent from $M$ gives a later $v\le L$; then $L,M,v$ is $120$
when $v<L$, and $101$ when $v=L$.
If $M=L$ and the maximum is repeated, a descent after its second
occurrence gives $110$, and a descent before that occurrence followed
by the second maximum gives $101$. Thus the maximum remains admissible.
No earlier value lies strictly between $L$ and $M$, so only the two
endpoints can be excluded from $[L,M]$. Appending $M+j$ has parameters
$(s-j+1,j)$ and type $W$ precisely when the old $M$ was excluded.
Repeating a legal $M$ gives $R_{s+1}$. Appending $M-k<M$ gives
$V_{s+1,k}$ in the ranges displayed above. This also proves inductively
that the four types exhaust all prefixes. The tree in
\cite[Table~1, p.~28]{CM} gives another enumeration of $B$.

\subsection{The two explicit systems}
For $X=B,C$, use the admissible sites of $X$ and the definitions
of $p,h,T_p$ in Section~\ref{sec:C}. Define
\[
\begin{aligned}
F_X(x)&=\sum_{w\in\bigcup_{n\ge1}\In_n(X)}z^{|w|}x^{s(w)},&
P_X(x)&=\sum_{w\in\bigcup_{n\ge1}\In_n(X)}z^{|w|}x^{s(w)}h(w),\\
H_X(x,y)&=\sum_{w\in\bigcup_{n\ge1}\In_n(X)}
 z^{|w|}x^{s(w)}T_{p(w)}(y),&
Z_X(x,y)&=H_X(x,y)-P_X(x).
\end{aligned}
\]
Thus $Z_X(x,y)$ sums the site polynomials with the maximum site removed:
that site, if admissible, has weight $y^0=1$.
In both classes $x$ marks the slack $s=n-M$, and a site's exponent
of $y$ counts unused admissible sites strictly to its right. In
particular, $y$ is a catalytic variable recording the unused sites above a selected site; a prefix contributes a polynomial in $y$, rather than necessarily a single monomial.
The unknowns in each system are $P_X(x),F_X(x),Z_X(x,y)$;
$H_X(x,y)=P_X(x)+Z_X(x,y)$ is then recovered.

\begin{proposition}\label{prop:common}
For $X=B$, the complete system in ordinary functional-equation form is
\begin{subequations}\label{eq:Bsystem}
\begin{align}
P_B(x)&=zx+zxP_B(x)+zx\frac{F_B(x)-F_B(1)}{x-1},\label{eq:BP}\\
F_B(x)&=P_B(x)+zxZ_B(x,1),\label{eq:BF}\\
Z_B(x,y)&=zx\left[\frac{Z_B(x,1)-Z_B(x,y)}{1-y}+Z_B(x,y)
 +\frac{P_B(x)-P_B(y)}{x-y}\right.\notag\\
&\hspace{16mm}\left.
 +\frac{1}{1-y}\left(\frac{F_B(x)-F_B(1)}{x-1}
                    -\frac{F_B(x)-F_B(y)}{x-y}\right)\right].\label{eq:BZ}
\end{align}
\end{subequations}
For $X=C$, it is
\begin{subequations}\label{eq:Csystem}
\begin{align}
P_C(x)&=zx+zxP_C(x)+zx\frac{F_C(x)-F_C(1)}{x-1},\label{eq:CPexplicit}\\
F_C(x)&=P_C(x)+zxZ_C(x,1),\label{eq:CFexplicit}\\
Z_C(x,y)&=zx\left[\frac{Z_C(x,1)-Z_C(x,y)}{1-y}
 +\frac{Z_C(x,y)-Z_C(y,y)}{x-y}
 +\frac{P_C(x)-P_C(y)}{x-y}\right.\notag\\
&\hspace{16mm}\left.
 +\frac{1}{1-y}\left(\frac{F_C(x)-F_C(1)}{x-1}
                    -\frac{F_C(x)-F_C(y)}{x-y}\right)\right].\label{eq:CZ}
\end{align}
\end{subequations}
All constant terms in $z$ are zero. The root conditions are
\[
[z]P_X(x)=[z]F_X(x)=x,\qquad [z]Z_X(x,y)=0,\qquad [z]H_X(x,y)=x.
\]
Each system has a unique solution in the coefficient rings stated
above, and $G_X(z)=1+F_X(1)$.
\end{proposition}

\begin{proof}
For either class, a new maximum is admissible and contributes
$J_X(x)=zx(F_X(x)-F_X(1))/(x-1)$.
Repetition of an admissible maximum leaves a single seen site and
contributes $zxP_X(x)$. A lower choice makes the old maximum
inadmissible: it creates a descent ending below that maximum and
$101$ forbids repeating it. Including the root proves
\eqref{eq:BP} and \eqref{eq:CPexplicit}.
Every admissible site gives a child of slack $s+1$, and every new
maximum is already counted in $J_X(x)$, so
\[
F_X(x)=zx+zxH_X(x,1)+J_X(x).
\]
Subtracting the equation for $P_X(x)$ proves \eqref{eq:BF} and
\eqref{eq:CFexplicit}.

For $C$, use the proved equation \eqref{eq:CH} and substitute
$H_C(x,y)=P_C(x)+Z_C(x,y)$. Its old-site new-maximum term becomes
\[
\frac{zx(P_C(x)-P_C(y))}{x-y}
 +\frac{zx(Z_C(x,y)-Z_C(y,y))}{x-y}.
\]
The lower-site difference loses its $P_C(x)$ terms. Subtracting
\eqref{eq:CPexplicit} leaves the independent source
\[
\frac{zx}{1-y}\left(\frac{F_C(x)-F_C(1)}{x-1}
                         -\frac{F_C(x)-F_C(y)}{x-y}\right),
\]
which proves \eqref{eq:CZ}.

For $B$, the legal interval is $[L,n]\setminus D(w)$, as proved
above. A lower choice $a$ leaves one seen site at $a$ and the $k(a)$
unused sites to its right; the old maximum is forbidden. Its site
polynomial is therefore $(1-y^{k(a)+1})/(1-y)$.
Repetition of $M$ leaves one site, the same expression for $k(M)=0$.
A new maximum $M+j$ discards every old site below $M$: the old
maximum becomes the second largest entry. The old maximum itself
remains admissible exactly when $h=1$. Consequently the child's
polynomial is
\[
hy^{j-1}+\frac{1-y^{j-1}}{1-y}+1.
\]
The finite sums \eqref{eq:sumidentities} now give the full equation
\begin{align}
H_B(x,y)&=zx+\frac{zx}{1-y}\bigl(H_B(x,1)-yH_B(x,y)\bigr)
 +\frac{zx(P_B(x)-P_B(y))}{x-y}\notag\\
&\quad+\frac{2-y}{1-y}J_B(x)
 -\frac{zx(F_B(x)-F_B(y))}{(1-y)(x-y)}.\label{eq:BH}
\end{align}
Inserting $H_B(x,y)=P_B(x)+Z_B(x,y)$ and subtracting
\eqref{eq:BP}, the site term left is
\[
\frac{Z_B(x,1)-yZ_B(x,y)}{1-y}
 =\frac{Z_B(x,1)-Z_B(x,y)}{1-y}+Z_B(x,y).
\]
The independent source is the same fraction involving $F_B$ as in
\eqref{eq:BZ}. This proves that equation too.

For completeness, the last fraction in either system is polynomial:
for $f(x)=x^a$ with $a\ge2$ it is
\[
\frac{1}{1-y}\left(\frac{x^a-1}{x-1}
                      -\frac{x^a-y^a}{x-y}\right)
 =\sum_{\substack{i,j\ge0\\i+j\le a-2}}x^iy^j,
\]
and it is zero for $a=0,1$. The other divided differences are also
polynomial finite sums. Starting with the constant terms zero,
compute the coefficient of $P_X(x)$ and $Z_X(x,y)$ at length $n$
from length $n-1$, then obtain that of $F_X(x)$ from its displayed
equation. This proves uniqueness and the root conditions.
The definition of $F_X(x)$ proves $G_X(z)=1+F_X(1)$.
\end{proof}

The slack distributions differ at length $4$: for $s=1,2,3,4$ they
are $(6,9,5,1)$ for $B$ and $(6,10,4,1)$ for $C$.
Also $R_1^B(z)=(z^2+z^3)/(1-z)^2$, whereas $R_1(z)$ for $C$ is
as in Section~\ref{sec:fixed}. The desired identity concerns the
specialization $x=1$.

\section{Kernel solution for 207(3)}\label{sec:Bkernel}
The elimination below uses catalytic variables and iterated kernel roots,
as in the generating-tree applications of the kernel method surveyed and
developed by Bousquet-M\'elou~\cite{BM}. We give the specialized elimination
for the present succession rule explicitly.

Write $b_B=F_B(1)$, $\beta_B=b_B+(z-1)/z$, and put
\[
\phi(x)=\frac1{1-zx},\qquad
D(x)=zx^2-x+1,\quad A(x)=D(x)+zx,\quad M(x)=zx+z-1.
\]
Here $D,A,M$ are auxiliary polynomials, distinct from the forbidden
set $D(w)$ in the generating tree.
Clearing denominators in \eqref{eq:BZ} gives
\begin{align}
(1-y+zxy)Z_B(x,y)
={}&zxZ_B(x,1)+zx(1-y)\frac{P_B(x)-P_B(y)}{x-y}\notag\\
&+zx\left(\frac{F_B(x)-b_B}{x-1}
                    -\frac{F_B(x)-F_B(y)}{x-y}\right).\label{eq:Bcleared}
\end{align}
The kernel vanishes at $y=\phi(x)$; this substitution is valid
coefficientwise because the catalytic coefficients are polynomials.
Solving \eqref{eq:BP} first gives
\begin{equation}\label{eq:PXsolution}
P_X(x)=\frac{zx}{1-zx}
 \left(1+\frac{F_X(x)-F_X(1)}{x-1}\right),\qquad X=B,C.
\end{equation}
With $Y=\phi(x)$, the vanishing-kernel condition, multiplied by $x-Y$,
is
\begin{align*}
0={}&(x-Y)\frac{F_B(x)-P_B(x)}{zx}
 +(1-Y)(P_B(x)-P_B(Y))\\
&+(x-Y)\frac{F_B(x)-b_B}{x-1}-F_B(x)+F_B(Y).
\end{align*}
Inserting \eqref{eq:PXsolution} and $Y=1/(1-zx)$ and clearing
rational denominators yields the scalar equation
\begin{equation}\label{eq:Bscalar}
M(x)A(x)F_B(x)+zx(1-zx)^2F_B(\phi(x))
 =xz^2D(x)\beta_B.
\end{equation}
These operations are identities over $\mathbb Q(z,x)$ with the three
series values $F_B(x),F_B(Y),b_B$ treated as independent symbols.
More explicitly, after the substitutions the right-hand side of the
vanishing-kernel display is $-1/(zx(1-zx)M(x))$ times the difference
between the left and right sides of \eqref{eq:Bscalar}. This verifies
all factors and signs in the reduction.

Define
\[
\lambda=\frac{1-\sqrt{1-4z}}{2z},\qquad
\alpha=\frac{1-z-\sqrt{1-6z+z^2}}{2z},\qquad \rho=\phi(\alpha),
\]
where the square roots have constant term $1$. Then $D(\lambda)=0$,
$A(\alpha)=0$, and both $\lambda$ and $\alpha$ have constant term $1$.
For the rational functions
\[
f(x)=\frac{xz^2D(x)}{M(x)A(x)},\qquad
 g(x)=-\frac{zx(1-zx)^2}{M(x)A(x)},
\]
put $a_i=\phi^i(a)$ and define
\begin{equation}\label{eq:orbits}
\Pi(a)=\prod_{i\ge0}g(a_i),\qquad
\mathscr S(a)=\sum_{i\ge0}f(a_i)\prod_{0\le j<i}g(a_j).
\end{equation}
The empty product is $1$. Below these are coefficientwise products
and sums, for $a=1$, $a=\rho$, or an indeterminate $a=x$.

\begin{proposition}\label{prop:Bproduct}
Set
\[
\tau=-z^2\alpha\rho^2,\qquad
\mathscr Q=\mathscr S(1)+\frac{\Pi(1)}{\Pi(\rho)}
                     \bigl(\tau-\mathscr S(\rho)\bigr).
\]
The orbit expressions are well-defined formal series,
$\mathscr Q=-z^2+O(z^3)$, and
\begin{equation}\label{eq:GBproduct}
\boxed{G_B(z)=1+\frac{(1-z)\mathscr Q}{z(\mathscr Q-1)}.}
\end{equation}
With $b_B=G_B(z)-1$ and $\beta_B=b_B+(z-1)/z$, the full slack series is
\begin{equation}\label{eq:FBproduct}
F_B(x)=\beta_B\left[\mathscr S(x)+\frac{\Pi(x)}{\Pi(\rho)}
                          \bigl(\tau-\mathscr S(\rho)\bigr)\right].
\end{equation}
The series $P_B(x)$ is given by \eqref{eq:PXsolution}. The remaining
series is recovered explicitly as
\begin{align}
Z_B(x,y)=\frac{1}{1-y+zxy}\left[F_B(x)-P_B(x)
 +zx(1-y)\frac{P_B(x)-P_B(y)}{x-y}\right.\notag\\
\left.\hspace{17mm}
 +zx\left(\frac{F_B(x)-b_B}{x-1}
                    -\frac{F_B(x)-F_B(y)}{x-y}\right)\right],\label{eq:ZBsolution}
\end{align}
and $H_B(x,y)=P_B(x)+Z_B(x,y)$.
\end{proposition}

\begin{proof}
The fixed point $\phi(\lambda)=\lambda$ and the identity
\[
\phi(u)-\lambda=\frac{z(u-\lambda)}{(1-zu)(1-z\lambda)}
\]
show that each orbit step raises the $z$-order of $u-\lambda$ by one.
The initial orders of $1-\lambda$ and $\rho-\lambda$ are $1$ and $2$,
respectively. Along either orbit $M(a_i)$ is a unit and $A(a_i)$ has
order $1$. Also $g(\lambda)=1$ and $f(\lambda)=0$.
To see the orders explicitly, factor
\[
D(u)=(u-\lambda)\bigl(z(u+\lambda)-1\bigr).
\]
Then the factor $z^2/A(u)$ in $f(u)$ raises the order by one.
For $g(u)-1$, its numerator is
$-zu(1-zu)^2-M(u)A(u)$, which vanishes at $u=\lambda$;
division by $A(u)$ can lower the order by at most one.
Consequently, on the $1$ orbit,
\[
g(a_i)-1=O(z^i)\ (i\ge1),\qquad f(a_i)=O(z^{i+2}),
\]
and on the $\rho$ orbit the respective orders are at least $i+1$
and $i+3$. Direct expansion gives $g(1)=1/2+O(z)$ and
$g(\rho)=1+O(z)$. The products therefore converge and are units,
with constant terms $1/2$ and $1$, and the sums converge.
For an indeterminate $x$, the first iterate has constant term $1$;
$A(\phi(x))=z(2-x)+O(z^2)$ is nonzero over $\mathbb Q(x)$.
After the second iterate the same order estimates apply, proving
convergence over $\mathbb Q(x)[[z]]$ as well.

Equation \eqref{eq:Bscalar} is
\[
F_B(u)=\beta_B f(u)+g(u)F_B(\phi(u)).
\]
After $k$ iterations this is a finite sum plus
$\prod_{i<k}g(a_i)F_B(a_k)$.
Since each coefficient of $F_B(u)$ is polynomial in $u$ and
$a_k\longrightarrow\lambda$, its coefficientwise limit is
\begin{equation}\label{eq:orbitlimit}
F_B(a)=\beta_B\mathscr S(a)+\Pi(a)F_B(\lambda).
\end{equation}
At $u=\alpha$, the scalar equation is valid before dividing by
$A(u)$. As $D(\alpha)=-z\alpha$, it gives
\[
F_B(\rho)=-\frac{z^2\alpha}{(1-z\alpha)^2}\beta_B=\tau\beta_B.
\]
Use \eqref{eq:orbitlimit} at $a=\rho$ to eliminate $F_B(\lambda)$.
At $a=1$ the result is $b_B=\beta_B\mathscr Q$. Solving this
linear equation using $\beta_B=b_B-(1-z)/z$ proves
\eqref{eq:GBproduct}; using $a=x$ proves \eqref{eq:FBproduct}.
The first terms are $\mathscr S(1)=-z^2/2+O(z^3)$,
$\mathscr S(\rho)=O(z^3)$ and $\tau=-z^2+O(z^3)$, so
$\mathscr Q=-z^2+O(z^3)$ and $\mathscr Q-1$ is a unit.
Finally \eqref{eq:PXsolution} and \eqref{eq:Bcleared} prove
\eqref{eq:ZBsolution}. Its apparent catalytic denominators are
interpreted after cancellation; the recovered coefficients are the
polynomials of the unique solution to Proposition~\ref{prop:common}.
\end{proof}

\section{The kernel reduction for 207(2)}\label{sec:Ckernel}
Clearing denominators in \eqref{eq:CZ} gives
\begin{align}
\kappa(x,y)Z_C(x,y)
={}&zx\left[(x-y)Z_C(x,1)-(1-y)Z_C(y,y)
 +(1-y)(P_C(x)-P_C(y))\right.\notag\\
&\left.\hspace{12mm}
 +(x-y)\frac{F_C(x)-F_C(1)}{x-1}-F_C(x)+F_C(y)\right],\label{eq:Ccleared}\\
\kappa(x,y)&=(1-y)(x-y)+zx(x-1).\notag
\end{align}
Let $b_C=F_C(1)$ and $\beta_C=b_C+(z-1)/z$. Define
\[
W_C(x)=\frac{zx(1-x)}{1-zx}F_C(x),\qquad
R(y)=\frac{z(1-y)^2}{(2-y)\{2-y+z(y^2-y-1)\}}.
\]
The following gives a scalar reduction with both kernel roots kept
inside legitimate formal compositions.

\begin{proposition}\label{prop:Ckernel}
For an indeterminate $t$, put
\[
T(t)=\frac{t(zt-t+1)}{zt^2-t+1},\qquad
x_0=T(t),\quad x_1=-\frac{t-1}{zt^2-t+1}=T(\phi(t)),\quad y_0=T(t)/t.
\]
Then $x_0,x_1$ are the two roots of $\kappa(x,y_0)=0$ and
\begin{equation}\label{eq:Cadditive}
W_C(x_0)-W_C(x_1)=\beta_C(x_0-x_1)R(y_0).
\end{equation}
\end{proposition}

\begin{proof}
Substitution into the quadratic verifies $\kappa(x_i,y_0)=0$.
The constant terms of $x_0,x_1,y_0$ are $t,1,1$, respectively;
thus compositions with the original prefix series are defined over
$\mathbb Q(t)[[z]]$. Rational divisions below may be taken in its
Laurent-series field.
In \eqref{eq:CH}, set $\kappa(x,y)=0$ and use
\eqref{eq:CF}, \eqref{eq:CP}, and \eqref{eq:PXsolution}.
The common diagonal boundary is
\begin{equation}\label{eq:Cboundary}
H_C(y,y)-\frac{F_C(y)}{1-y}
 =a(x,y)F_C(x)+q(x,y)b_C+c(x,y),
\end{equation}
where
\begin{align*}
a(x,y)&=\frac{x-y}{(x-1)(1-zx)}
                      +\frac{x-y-zx}{zx(1-y)},\\
q(x,y)&=-\frac{x-y}{(x-1)(1-zx)},\qquad
c(x,y)=\frac{x-y}{1-zx}-\frac{x-y}{1-y}.
\end{align*}
For an explicit derivation, the cleared equation for $H_C(x,y)$
on the kernel is
\begin{align*}
H_C(y,y)-\frac{F_C(y)}{1-y}
={}&(x-y)(1+P_C(x))
 +\frac{x-y}{1-y}H_C(x,1)\\
&+\frac{(2-y)(x-y)}{zx(1-y)}J_C(x)-\frac{F_C(x)}{1-y}.
\end{align*}
Insert $zxH_C(x,1)=F_C(x)-zx-J_C(x)$ and
$1+P_C(x)=(1+J_C(x))/(1-zx)$; the displayed $a,q,c$ result.

In the following root calculations set $y=y_0$. On the kernel, $(x-y)/(x-1)=-zx/(1-y)$, which gives
\[
a(x,y)\frac{1-zx}{zx(1-x)}=-\frac{2-y}{(1-y)^3}.
\]
The root sum and product of
$zx^2+(1-y-z)x-y(1-y)$ give
\[
(1-zx_0)(1-zx_1)=2-y+z(y^2-y-1).
\]
Also, writing $\delta q=q(x_0,y)-q(x_1,y)$, the kernel relation gives
\[
\delta q=\frac{z(x_0-x_1)}{(1-y)(1-zx_0)(1-zx_1)},\qquad
c(x_0,y)-c(x_1,y)=\frac{z-1}{z}\delta q.
\]
Subtract \eqref{eq:Cboundary} for the two roots; its left sides
are identical. The result is
\[
\frac{2-y}{(1-y)^3}\bigl(W_C(x_0)-W_C(x_1)\bigr)
 =\beta_C\delta q.
\]
The preceding product formula proves \eqref{eq:Cadditive} at $y=y_0$.
\end{proof}

The $C$ scalar equation has an additive forcing term, whereas the
$B$ scalar equation has the convergent product solution above. An
unrestricted iteration of \eqref{eq:Cadditive} cannot be justified by
ordinary prefix-series composition. For example, direct iteration gives
\[
T(\phi^3(t))=\frac{zt+z-1}{z(zt^2-t+1)},
\]
which has order $-1$ in $z$.
There is a concrete obstruction at the critical value $y=2$. The
roots then are
\[
x_-=1+\alpha,\qquad x_+=1+\frac1{z\alpha},\qquad 1-zx_+=z\alpha.
\]
Using the rational summands from Section~\ref{sec:fixed}, a proposed
continuation of the slack series would require the sum
\[
\sum_{m\ge0}x_+^{-m}R_m(zx_+)
 =\sum_{m\ge0}\frac{x_+^{-m}P_m(zx_+)}{(1-zx_+)^{m+1}}.
\]
Since $x_+^{-1}=z+O(z^2)$, $zx_+=1+O(z)$, and $P_m(1)=1$, every
summand has order $-1$ and leading coefficient $1$. The sum therefore
does not exist coefficientwise in the formal Laurent-series topology.
This calculation prevents using that branch to read off $b_C$.
The coefficientwise rational expansion \eqref{eq:expansion}, whose
numerators are specified in \eqref{eq:Rm}, gives one all-length formula
for $G_C$. We next obtain a Lambert series expression and use it to
identify $G_C$ with the B product formula \eqref{eq:GBproduct}.

\section{The Class 207 equality and its common generating function}
\label{sec:class207-equality}
The C kernel equation gives a Lambert series quotient, whereas the
published B tree gives a difference of two product orbit sums. We evaluate
that difference by a cubic limit of the nonterminating
$q$-Saalsch\"utz summation. Together with the A--C bijection of
Callan--Mansour~\cite[Section~4.10]{CM}, this proves the Class~207 case
of their Conjecture~1.

Let $q=\lambda-1$, where $\lambda=(1-\sqrt{1-4z})/(2z)$ is the Catalan series, so $z=q/(1+q)^2$, and set
\[
 L=1-2q-2q^2+q^3,\qquad q+Lr+q^2r^2=0,\qquad r=-q+O(q^2).
\]
This specifies $r\in q\QQ[[q]]$ uniquely; explicitly,
$r=-2q/(L+\sqrt{L^2-4q^3})$, with square root $1+O(q)$.
Define
\begin{equation}
 T(q,r)=\sum_{k\ge1}\frac{k^2(qr)^k}{1-q^k},\qquad
 \theta=\frac{(1-q)^3(1+q)}{q^2},\qquad H=1+\theta T(q,r(q)).
 \label{count:eq:H}
\end{equation}
\begin{theorem}\label{count:thm:main}
For all $n\ge0$, the three classes $A,B,C$ have the same number of
inversion sequences of length $n$. Their common generating function satisfies
\begin{equation}
 \boxed{G_A\!\left(\frac q{(1+q)^2}\right)
 =G_B\!\left(\frac q{(1+q)^2}\right)
 =G_C\!\left(\frac q{(1+q)^2}\right)
 =\frac{(1+q)^2}{q}+\frac1H.}
 \label{count:eq:main}
\end{equation}
Here $H=-q+O(q^2)$; the two Laurent expressions on the right have
canceling negative powers and their sum belongs to $\QQ[[q]]$.
\end{theorem}
The proof first obtains the two counting formulas and then evaluates the
orbit difference that connects them. The summation used in the last step is
classical~\cite[Eq.~(II.24)]{GR}; see also~\cite[Eq.~17.7.5]{DLMF}.
We give the specialization and cubic limit explicitly.

\subsection{The C boundary and its Lambert series}
Use the C statistics of Section~\ref{sec:C} and, in this subsection,
write $P=P_C$, $F=F_C$ and $Z=Z_C$. In particular,
$G_C=1+F(1)$ and $\deg_x[z^n]P\le n$, since the slack is at most
the prefix length. Equations~\eqref{eq:CPexplicit}--\eqref{eq:CZ}
provide the complete system used below.

Write $\kappa(x,y)=(1-y)(x-y)+zx(x-1)$. At a kernel root, clearing
\eqref{eq:CZ} and substituting \eqref{eq:CPexplicit}--\eqref{eq:CFexplicit} gives
\begin{equation}
 (1-y)^2\left(Z(y,y)+\frac{P(y)}{zy}\right)
 =(1-x)(G_C-x)+\frac{(2-y)(1-x)^2}{1-y}P(x).
 \label{count:eq:kernel}
\end{equation}
For two roots $x,x'$, subtraction and $x+x'=1+(y-1)/z$ give
\begin{equation}
 (2-y)\bigl((1-x)^2P(x)-(1-x')^2P(x')\bigr)
 =(1-y)(x-x')\left(G_C+\frac{1-y}{z}\right).
 \label{count:eq:two}
\end{equation}
In particular, for an indeterminate $t$, put
\[
 d=t^2-t+z,\quad X=\frac{t(t-1)}d,\quad
 X'=\frac{z-t}d,\quad Y=\frac{(t-1)(t-z)}d.
\]
They satisfy $\kappa(X,Y)=\kappa(X',Y)=0$ and have nonnegative $z$-orders.
Substitution in \eqref{count:eq:two} gives
\begin{equation}
 (d+tz)\bigl(z^2P(X)-t^4P(X')\bigr)
 =zt(t^2-z)(dG_C+t).
 \label{count:eq:input}
\end{equation}
All substitutions so far are coefficientwise legitimate over $\QQ(t)[[z]]$.

Put $D(u)=qu^2+Lu+q^2$ and define
\begin{align*}
 j_G(u)&=\frac{(u-1)(u-q)(u^2-q)(1-q)^2}{uD(u)},&
 j_0(u)&=\frac{(u-1)^2(u-q)^2(u^2-q)(1+q)}{u^2D(u)},\\
 \cLam(j)&=\sum_{k\ge1}\frac{k^2q^k}{1-q^k}[u^k]j(u),&
 \Lambda_G&=\cLam(j_G),\qquad\Lambda_0=\cLam(j_0).
\end{align*}
Laurent expansions here are $q$-first: expand in $q$ with Laurent-polynomial
coefficients in $u$, and then extract $[u^k]$.
Take
\[
 t=\frac{u-q}{(1+q)(u-1)},\quad
 d=\frac{(1-q)^2u}{(1+q)^2(u-1)^2},\quad
 X=1-\frac q{(1-q)^2}(u+u^{-1}-2),\quad Q(u)=P(X)/d^2.
\]
Equation~\eqref{count:eq:input} becomes, with $\cGc=G_C(q/(1+q)^2)$ and $\cGb=G_B(q/(1+q)^2)$,
\begin{equation}
 Q(u)-Q(q/u)=\frac{(1+q)^4}{q(1-q)^3}
                 (j_G(u)\cGc+j_0(u)).\label{count:eq:add}
\end{equation}
The substitution deserves an order check. Writing $w=u+u^{-1}-2$,
$Q$ is a unit in $\QQ[[q]]$ times $w^2P(z,1-qw/(1-q)^2)$.
Since $\deg_x[z^n]P\le n$, a term $w^{j+2}$ from $P$ has $q$-order at
least $2j$. Consequently, if $Q_k=[u^k]Q$, then
$\operatorname{ord}_q Q_k\ge2|k|-4$ for $|k|\ge3$.
This proves that $Q(q/u)$ and the coefficient sums below are well defined
in the $q$-adic Laurent field. Also $Q(u)=Q(1/u)$ and $Q''(1)=0$, since
$d^{-2}$ supplies a fourth-order zero at $u=1$.
Taking positive coefficients in \eqref{count:eq:add} gives
\[
 Q_k=-\frac{(1+q)^4}{q(1-q)^3}\frac{q^k}{1-q^k}
           [u^k](j_G\cGc+j_0).
\]
Thus $0=Q''(1)=2\sum_{k>0}k^2Q_k$ implies
\begin{equation}
 \cGc=-\Lambda_0/\Lambda_G.\label{count:eq:Cquot}
\end{equation}
To evaluate it, use the two rational identities
\begin{align}
 j_0+\frac{(1+q)^2}{q}j_G
 &=\frac{1+q}{q}\frac{(u-1)(u-q)(u^2-q)}{u^2},\label{count:eq:red}\\
 j_G&=\frac{(1-q)^2}{q}\left(u-\frac q u\right)
 +\frac{(1-q)^4(1+q)}{q^2}
       \left(\frac1{1-ru}-\frac1{1-qr/u}\right).\label{count:eq:partial}
\end{align}
The second follows by $D(u)=-(q/r)(1-ru)(u-qr)$ and the quadratic equation
for $r$. Applying $\cLam$ yields
\[
 \Lambda_0+\frac{(1+q)^2}{q}\Lambda_G=-(1-q),\qquad
 \Lambda_G=(1-q)H.
\]
The first coefficients are $H=-q+q^2-q^3+O(q^4)$, so division is
legitimate in $\QQ((q))$. Therefore
\begin{equation}
 \cGc=\frac{(1+q)^2}{q}+\frac1H.\label{count:eq:Cfinal}
\end{equation}

\subsection{The published B tree and its orbit sums}
We use the tree for $B$ in Callan--Mansour~\cite[Table~1, p.~28]{CM}.
Let $a_m(z),b_{m,j}(z),c_{m,j}(z)$ count nodes of its three types by length,
with root $a_1$ at length one and $1\le j\le m$.
Writing $m'=m+1-j$, the rules are
\begin{align*}
 a_m&\longrightarrow a_{m+1},\ b_{m,1},\ldots,b_{m,m},\\
 b_{m,j}&\longrightarrow c_{m,j},c_{m-1,j-1},\ldots,c_{m',1},\
                    a_{m'+1},b_{m',1},\ldots,b_{m',m'},\\
 c_{m,j}&\longrightarrow c_{m+1,j},c_{m,j-1},\ldots,c_{m'+1,1},\
                    a_{m'+1},b_{m',1},\ldots,b_{m',m'}.
\end{align*}
There is exactly one $a$-child per node, and the empty object precedes
the root. Hence $\sum_m a_m=zG_B$. Define a separate boundary $f$ by
\[
 a(x)=\sum_{m\ge1}a_m x^{m-1}=z(1+f(x)),\quad
 \mathsf B(x,y)=\sum_{m,j}b_{m,j}x^{m-j}y^{j-1},\quad
 \mathsf C(x,y)=\sum_{m,j}c_{m,j}x^{m-j}y^{j-1}.
\]
In particular $f(1)=G_B-1$. This $f$ is a boundary of the published tree;
it is different from both slack series $F_B$ and $F_C$ of
Proposition~\ref{prop:common}.
The rules give $b_{m,j}=a_{m+1}$, since their parent sets coincide, and
therefore
\begin{align}
 \mathsf B(x,y)&=z\frac{f(x)-f(y)}{x-y},\notag\\
 (1-y+zxy)\mathsf C(x,y)
 &=z(\mathsf B(x,1)-y\mathsf B(x,y))+zx\mathsf C(x,1),\notag\\
 a(x)&=z+zx a(x)+zx(\mathsf B(x,1)+\mathsf C(x,1)).\label{count:eq:Btree}
\end{align}
At $y=(1-zx)^{-1}$, eliminate $\mathsf C(x,1)$ to obtain
$(1-zx)f(x)-zx=(x-1)\mathsf B(x,1)+y\mathsf B(x,y)$.
Inserting $\mathsf B$ and clearing denominators proves
\begin{equation}
 zf\!\left(\frac1{1-zx}\right)
 =(1-zx)(1-x+zx+zx^2)f(x)
       +z(1-x+zx^2)(f(1)-x).\label{count:eq:Bfun}
\end{equation}
Thus the scalar equation determines the total of the B class through
$f(1)=G_B-1$, without identifying $f$ with the slack series $F_B$.

Put
\[
 x(v)=\frac{(1+q)(1-v)}{1-qv},\quad E(v)=q+Lv+q^2v^2,\quad
 \mathfrak A(v)=\frac{(1-q^2v)E(v)}{q(1-qv)^3},\quad
 \mathfrak b(v)=\frac{(1-q)^2v}{(1-qv)^2}.
\]
For $\mathcal F(v)=f(x(v))$ and $g=G_B-1$, equation~\eqref{count:eq:Bfun} is
\begin{equation}
 \mathcal F(qv)=\mathfrak A(v)\mathcal F(v)+\mathfrak b(v)(g-x(v)).\label{count:eq:Bdiff}
\end{equation}
With $(a;q)_\infty=\prod_{j\ge0}(1-aq^j)$ and products of several
arguments denoting products of these factors, define
\begin{equation}
 U(v)=\frac{(qv/r,q^2rv;q)_\infty}{(qv;q)_\infty^2},\quad
 W(v)=\frac{(1-q)^2v}{1-q^2v}U(v),\quad
 J(v)=\frac{(v/r,qrv;q)_\infty}{(1-qv)(qv;q)_\infty^2}.\label{count:eq:products}
\end{equation}
Set
\[
 M_a=\sum_{k\ge0}W(aq^k),\quad N_a=\sum_{k\ge0}W(aq^k)x(aq^k),\quad
 \mathcal D=M_1-M_r,\quad\mathcal N=N_1-N_r.
\]
Cancellation of product factors, using $E(r)=0$, gives
\[
 J(qv)\mathfrak A(v)=J(v),\qquad J(qv)\mathfrak b(v)=W(v).
\]
Multiplying \eqref{count:eq:Bdiff} by $J(qv)$ and summing on an orbit therefore
telescopes. On the orbit from $1$, use $\mathcal F(1)=f(0)=0$;
on the orbit from $r$, use $J(r)=0$ without dividing by $\mathfrak A(r)=0$.
Since type $a_m$ first occurs at length $m$, $\deg_x[z^n]f\le n$.
Thus $\mathcal F(aq^k)\to\mathcal F(0)$ and $J(aq^k)\to1$ $q$-adically.
Since $J(0)=1$, the result in either case is
$\mathcal F(0)=gM_a-N_a$. Thus
\begin{equation}
 \cGb=1+\frac{\mathcal N}{\mathcal D}.\label{count:eq:Bquot}
\end{equation}
The formal limits exist: $W(q^k)=O(q^k)$,
$W(rq^k)=O(q^{k+1})$, and $\mathcal D=2+O(q)$.

For $K=(1+q+q^2)/q$, finite product cancellation also gives
\begin{equation}
 J(qv)-J(v)=W(v)(x(v)-K).\label{count:eq:tel}
\end{equation}
Indeed, after division by $U(v)$, the required polynomial identity is
\[
 (1-qv)^3-(1-q^2v)(1-v/r)(1-qrv)
 =-\frac{(1-q)^2}{q}v(1-q^3v).
\]
Since $J(0)=1$ and $J(r)=0$, summing \eqref{count:eq:tel} yields
$\mathcal N-K\mathcal D=-J(1)$. Equation~\eqref{count:eq:Bquot} becomes
\begin{equation}
 \cGb=\frac{(1+q)^2}{q}-\frac{J(1)}{\mathcal D}.\label{count:eq:Bfinal}
\end{equation}
It remains to prove $\mathcal D=-J(1)H$.

\subsection{A two-parameter summation}
\begin{lemma}\label{count:lem:sum}
For independent real parameters $0<q<1$ and $-1<r<0$, define
$U,W,J,M$ by \eqref{count:eq:products} and the orbit sums above. Then
\begin{equation}
 M_1-M_r=-J(1)(1-\beta T(q,r)),\qquad
 \beta=\frac{(1-q)(1-r)(1-qr)}{qr}.\label{count:eq:sum}
\end{equation}
The quadratic relation between $q$ and $r$ is not assumed here.
\end{lemma}
\begin{proof}
For $(a;q)_n=\prod_{j=0}^{n-1}(1-aq^j)$, use the convention
\[
 {}_3\phi_2\!\left(\begin{matrix}a,b,c\\ e,f\end{matrix};q,q\right)
 =\sum_{n\ge0}\frac{(a,b,c;q)_n}{(q,e,f;q)_n}q^n.
\]
The nonterminating $q$-Saalsch\"utz identity~\cite[Eq.~(II.24)]{GR},
in the form of~\cite[Eq.~17.7.5]{DLMF}, holds for $0<|q|<1$ and
$ef=abcq$, with denominator factors nonzero. It states
\begin{align}
 {}_3\phi_2\!\left(\begin{matrix}a,b,c\\ e,f\end{matrix};q,q\right)
 &+\frac{(q/e,a,b,c,qf/e;q)_\infty}
              {(e/q,aq/e,bq/e,cq/e,f;q)_\infty}
 {}_3\phi_2\!\left(\begin{matrix}aq/e,bq/e,cq/e\\ q^2/e,qf/e\end{matrix};q,q\right)\notag\\
 &=\frac{(q/e,f/a,f/b,f/c;q)_\infty}
               {(aq/e,bq/e,cq/e,f;q)_\infty}.\label{count:eq:saal}
\end{align}
Set
\[
 a=b=qr,\quad c=q^2r,\quad e=q^3rs,\quad f=q^2r^2/s,
 \qquad h=1-1/s,\quad t=1-h.
\]
The balance is exact. Denote the first series by $\Phi(s)$, the multiplier
of the second by $\mathcal A(s)$, and its terms by $b_n(s)$; explicitly,
\begin{align*}
 b_n(s)&=\frac{(t/q,t/q,t;q)_n}{(q,t/(qr),rt^2;q)_n}q^n,\\
 \mathcal A(s)&=
 \frac{(t/(q^2r),qr,qr,q^2r,rt^2;q)_\infty}
      {(q^2r/t,t/q,t/q,t,q^2r^2t;q)_\infty}.
\end{align*}
At $s=1$, the first series is $\Phi(1)=M_r/W(r)$, as follows by comparing
the $n=0$ terms and successive ratios. The multiplier has a triple pole:
\[
 (t/q;q)_\infty^2(t;q)_\infty=(1-t/q)^2h^3(qt;q)_\infty^3.
\]
Write $A_{-3}=\lim_{h\to0}h^3\mathcal A(s)$. Finite shifts give
\begin{align}
 A_{-3}&=\frac{q^2(q^{-2}/r;q)_\infty(qr;q)_\infty^2(r;q)_\infty}
 {(1-q)^2(q^2r^2;q)_\infty(q;q)_\infty^3},&
 W(r)A_{-3}&=\beta J(1).\label{count:eq:Aminus}
\end{align}
For clarity, the last cancellation reduces to
\[
 \frac{W(r)A_{-3}}{J(1)}
 =\frac{q^2r(1-q)(1-r)}{1-q^2r}
        (1-q^{-2}/r)(1-q^{-1}/r)=\beta.
\]
For every $n\ge2$,
\begin{equation}
 \lim_{h\to0}\mathcal A(s)b_n(s)
 =-\frac{W(q^{n-2})}{W(r)}.\label{count:eq:tail}
\end{equation}
A base case and ratio prove this for all indices. At $n=2$,
\[
 \lim_{h\to0}\frac{b_2}{h^3}
 =\frac{(1-q^{-1})^2q^2}
 {(1-q^2)(1-(qr)^{-1})(1-r^{-1})(1-r)(1-qr)};
\]
multiplication by $W(r)A_{-3}/W(1)$ gives $-1$ by finite product shifts.
The common successive ratio for \eqref{count:eq:tail} is
\[
 q\frac{(1-q^{n-1})^2(1-q^n)}
 {(1-q^{n+1})(1-q^{n-1}/r)(1-rq^n)}.
\]
In particular the $n=2$ term retains the orbit endpoint $W(1)$.

Divide the right side of \eqref{count:eq:saal} by $\mathcal A(s)$ and call it
\[
 P(h)=\frac{(qrt;q)_\infty^2(rt;q)_\infty(q^2r/t;q)_\infty}
 {(qr;q)_\infty^2(q^2r;q)_\infty(rt^2;q)_\infty}.
\]
The two exceptional terms are $b_0=1$ and
\[
 b_1=-\frac{rh(1-q-h)^2}
 {(1-q)(1-qr-h)(1-r(1-h)^2)}.
\]
Let
\[
 R(h)=\frac{(1-r+rh)(1-qr)(1-h)}
 {(1-qr-h)(1-r(1-h)^2)}.
\]
The following two identities isolate the finite part:
\begin{align}
 R-1-b_1&=\frac{rqh^3}{(1-q)(1-qr-h)(1-r(1-h)^2)},\label{count:eq:cubic}\\
 \frac{P}{R}&=\frac{(qrt;q)_\infty^3(qr/t;q)_\infty}
 {(qr;q)_\infty^3(qrt^2;q)_\infty}.\label{count:eq:PR}
\end{align}
The first is a polynomial identity after clearing denominators; the second
uses four single-step product shifts. Put $a=qr$, $\lambda=\log t$ and
$L_j(a)=(a\partial_a)^j\log(a;q)_\infty$. Taking logarithms of
\eqref{count:eq:PR} gives
\[
 \log(P/R)=\sum_{j\ge1}\frac{3+(-1)^j-2^j}{j!}L_j(a)\lambda^j.
\]
The linear and quadratic coefficients vanish, the cubic coefficient is
$-L_3(a)$, and
\[
 L_3(a)=-\sum_{k\ge1}\frac{k^2a^k}{1-q^k}=-T(q,r),\qquad
 \lambda^3=-h^3+O(h^4).
\]
Consequently $P-R=-Th^3+O(h^4)$, and \eqref{count:eq:cubic} proves
\begin{equation}
 [h^3](P-1-b_1)=\frac1\beta-T.\label{count:eq:finitepart}
\end{equation}
Now \eqref{count:eq:saal} reads
$\Phi(s)=\mathcal A(s)(P-1-b_1-\sum_{n\ge2}b_n)$.
Using \eqref{count:eq:tail} and \eqref{count:eq:finitepart} gives
\[
 \frac{M_r}{W(r)}=A_{-3}(1/\beta-T)+\frac{M_1}{W(r)}.
\]
Equation~\eqref{count:eq:Aminus} proves \eqref{count:eq:sum}.

To justify the limit through the infinite tail, for $n\ge2$ extract the
three zeros explicitly:
\[
 \frac{b_n}{h^3}=
 \frac{(1-t/q)^2(qt;q)_{n-2}^2(qt;q)_{n-1}q^n}
 {(q;q)_n(t/(qr);q)_n(rt^2;q)_n}.
\]
For fixed $0<q<1$, $-1<r<0$, on a sufficiently small closed complex
neighborhood of $h=0$ the denominator products are bounded away from
zero uniformly in $n$, and the numerator products are bounded. Therefore
$|b_n/h^3|\le Cq^n$. The two ${}_3\phi_2$ series have argument $q$,
so they converge, and this bound justifies termwise passage to the limit.
All product shifts and logarithmic derivatives above are then valid there.
\end{proof}

\subsection{Formal specialization and completion of the proof}
On the branch $q+Lr+q^2r^2=0$, elementary algebra gives
\[
 \beta=-\theta,\qquad J(1)=\frac{1-q^2}{q}U(1).
\]
The implicit function theorem makes $r(q)$ analytic near zero; moreover,
$q/r=-L-q^2r$ removes the apparent singularity in $U$.
The functions $U(1)$, $M_1$, $M_r$ and $qJ(1)$ are analytic germs,
since their nonconstant factors and orbit tails converge uniformly on
sufficiently small disks. In $qJ(1)$ the prefactor $q$ removes its simple pole.
In particular,
\[
 \mathcal D(0)=2,\qquad (qJ(1))(0)=2.
\]
For sufficiently small positive $q$, the branch has $-1<r(q)<0$.
Lemma~\ref{count:lem:sum} thus gives
\[
 q\mathcal D+(qJ(1))H=0
\]
on a real interval. Both sides are analytic at zero, so the identity theorem
proves the corresponding formal identity. It gives $H=-q+O(q^2)$ and
$\mathcal D=-J(1)H$. Equations~\eqref{count:eq:Cfinal} and \eqref{count:eq:Bfinal}
coincide. Since $z=q+O(q^2)$ is formally invertible, $G_B=G_C$ as length
series. The published bijection $G_A=G_C$ completes Theorem~\ref{count:thm:main}.

\section{Exact coefficients and verification}
\label{sec:checks}
The checks are implemented in Python~\cite{Python}; exact rational identities
use SymPy~\cite{SymPy}. The optional numerical audit in
\texttt{class207-proof/referee/audit.py} uses mpmath~\cite{mpmath} and is
separate from the exact-integer and symbolic checks below.
The initial counts are
\begin{center}
\begin{tabular}{r r r r r r}
\toprule
$n$&$207(2)$&$207(3)$&$n$&$207(2)$&$207(3)$\\
\midrule
0&1&1&7&1375&1375\\
1&1&1&8&6052&6052\\
2&2&2&9&27425&27425\\
3&6&6&10&127273&127273\\
4&21&21&11&602487&602487\\
5&80&80&12&2900432&2900432\\
6&324&324&&&\\
\bottomrule
\end{tabular}
\end{center}
The recurrence in Proposition~\ref{prop:recurrence}, an independent
implementation of the four-type $B$ tree, and an implementation of the
$C$ site-word tree reproduce these counts through length $12$.
Literal pattern checks of inversion sequences through length $8$ verify
the initial counts for both Class~152 sets and both Class~207 sets,
and every legal-site transition from Class~207 parents in that range.
The same enumeration checks the stated length-$4$ slack distributions
and the first four fixed-maximum rational functions through length $8$.
These finite checks support the enumerations. The all-length
Class~207 equality follows from Theorem~\ref{count:thm:main}.
In addition, direct exact-integer expansion of the two product orbits
at internal order $136$ gives a zero residual for
$q\mathcal D+(qJ(1))H$ through degree $129$. The common Lambert
formula agrees with the C legal-site recurrence and the four-type
B tree through length $128$. These computations check the formulas
independently of the limiting argument; they are not premises of the theorem.

\paragraph{Data and code availability}
The accompanying supplementary material contains \texttt{counts.csv},
the independent generators and literal pattern checker in
\texttt{check\_counts.py}, the recurrence program in Appendix~A, and exact checks of the scalar eliminations
and generating-function identities in \texttt{symbolic\_check.py}.
Its \texttt{README.md} specifies the finite ranges, software versions,
commands, and expected outputs. The files are supplied with the
manuscript in the \texttt{supplement} directory. The subdirectory
\texttt{supplement/class207-proof} contains the exact product and
Lambert-series checks, the denominator-cleared symbolic certificates,
and separate checks of the published B-tree boundary and literal avoidance.

Combining Theorem~\ref{count:thm:main} with the earlier enumerations gives
\[
1+z\sum_{r\ge0}b_r(1)z^r
 =1+\sum_{m\ge0}\frac{P_m(z)}{(1-z)^{m+1}}
 =1+\frac{(1-z)\mathscr Q}{z(\mathscr Q-1)}.
\]
The first equality comes from the two enumerations of C, the last
expression from the B product solution, and the equality between them
from the summation proof in Section~\ref{sec:class207-equality}.

\section{A kernel calculation for 207(3)}
By Proposition \ref{prop:common}, we have 
\begin{subequations}\label{tm:eq:Bsystem}
\begin{align}
P_B(x)&=zx+zxP_B(x)+zx\frac{F_B(x)-F_B(1)}{x-1},\label{tm:eq:BP}\\
F_B(x)&=P_B(x)+zxZ_B(x,1),\label{tm:eq:BF}\\
Z_B(x,y)&=zx\left[\frac{Z_B(x,1)-Z_B(x,y)}{1-y}+Z_B(x,y)
 +\frac{P_B(x)-P_B(y)}{x-y}\right.\notag\\
&\hspace{16mm}\left.
 +\frac{1}{1-y}\left(\frac{F_B(x)-F_B(1)}{x-1}
                    -\frac{F_B(x)-F_B(y)}{x-y}\right)\right].\label{tm:eq:BZ}
\end{align}
\end{subequations}
By solving \eqref{tm:eq:BP} and \eqref{tm:eq:BF} for $P_B(x)$ and $F_B(x)$, we obtain 

\begin{align}
F_B(x)&=\frac{zx((1-x)(1-xz)Z_B(x,1)+F_B(1)-x+1)}{x^2z-x+1},\label{tm:eq:t00}\\ P_B(x)&=\frac{zx(F_B(1)-zxZ_B(x,1)-x+1)}{x^2z-x+1}.\label{tm:eq:t01}
\end{align}

By substituting \eqref{tm:eq:t00}-\eqref{tm:eq:t01} into \eqref{tm:eq:BZ}, we have

\begin{align*}
0={}&\frac{xyz-y+1}{y-1}Z_B(x,y)-\frac{(x-1)(x^2z-xyz+xz-x+y)zx}
{(x-y)(x^2z-x+1)(y-1)}Z_B(x,1)\\
&\quad+\frac{yz^2x}{(y^2z-y+1)(x-y)}Z_B(y,1)-\frac{(xz+yz-1)zx}
{(x^2z-x+1)(y^2z-y+1)}F_B(1)\\
&\quad+\frac{(xy-x-y)z^2x}
{(x^2z-x+1)(y^2z-y+1)}.
\end{align*}

By taking $y=1/(1-xz)$, we obtain 
\begin{align}
Z_B(x,1)&=\frac{zx}{(x^2z+xz-x+1)(x-1)}Z_B\left(\frac{1}{1-xz},1\right)\nonumber\\
&+\frac{x(x^2z^2-2xz-z+1)}{(xz-1)(x^2z+xz-x+1)(x-1)}F_B(1)
+\frac{xz(x^2z-1)}{(xz-1)(x^2z+xz-x+1)(x-1)}.\label{tm:eq:t1}
\end{align}

Define $\rho_0(x,z)=x$ and $\rho_{j+1}(x,z)=\phi(\rho_j(x,z))$,
where $\phi(x)=1/(1-zx)$. By induction on $j\ge1$,
\[
\rho_j(x,z)=\frac{-x\sqrt z\,U_{j-2}(1/(2\sqrt z))+U_{j-1}(1/(2\sqrt z))}
 {\sqrt z\{-x\sqrt z\,U_{j-1}(1/(2\sqrt z))+U_j(1/(2\sqrt z))\}},
\]
where $U_{-1}=0$, $U_0=1$, and $U_j(v)=2vU_{j-1}(v)-U_{j-2}(v)$.
Indeed $K_j=z^{j/2}U_j(1/(2\sqrt z))$ satisfies
$K_{-1}=0$, $K_0=K_1=1$, $K_j=K_{j-1}-zK_{j-2}$, and the displayed
quotient is $(K_{j-1}-zxK_{j-2})/(K_j-zxK_{j-1})$.
Substitution into $\rho_{j+1}=1/(1-z\rho_j)$ proves the formula,
including its initial case $j=1$.

Write $U_B(x)=Z_B(x,1)$ and define
\[
 \alpha_T(x)=\frac{zx}{(x^2z+xz-x+1)(x-1)},
\]
\begin{align*}
 \beta_T(x)&=\frac{x(x^2z^2-2xz-z+1)}
 {(xz-1)(x^2z+xz-x+1)(x-1)},\\
 \gamma_T(x)&=\frac{xz(x^2z-1)}
 {(xz-1)(x^2z+xz-x+1)(x-1)}.
\end{align*}
Equation~\eqref{tm:eq:t1} is the rational affine equation
$U_B(x)=\alpha_T(x)U_B(\phi(x))+\beta_T(x)b_B+\gamma_T(x)$,
where $b_B=F_B(1)$. For every finite $N\ge1$ its exact iteration is
\begin{align}
 U_B(x)={}&\sum_{j=0}^{N-1}
  \bigl(\beta_T(\rho_j)b_B+\gamma_T(\rho_j)\bigr)
    \prod_{i=0}^{j-1}\alpha_T(\rho_i)\notag\\
 &+\prod_{i=0}^{N-1}\alpha_T(\rho_i)U_B(\rho_N).
 \label{tm:eq:finite}
\end{align}
The empty product is $1$. Substituting the one-step affine equation
for the last value proves this formula by induction on $N$.
The remainder must be retained: over $\QQ(x)((z))$,
$\operatorname{ord}_z\alpha_T(\rho_0)=1$ and
$\operatorname{ord}_z\alpha_T(\rho_i)=-1$ for $i\ge1$.
For $i=1$, $\rho_1=1+zx+O(z^2)$, and for $i\ge2$,
$\rho_i=1+z+O(z^2)$; substitution in $\alpha_T$ proves these orders.
Thus the product in the remainder has order $2-N$. The prefix $01$
gives $U_B(z;x)=z^2x+O(z^3)$, and $\rho_N=1+O(z)$ for $N\ge1$.
Consequently the full remainder in \eqref{tm:eq:finite} has
\[
 \operatorname{ord}_z\!\left(
   \prod_{i=0}^{N-1}\alpha_T(\rho_i)U_B(\rho_N)\right)=4-N.
\]
It does not tend to zero coefficientwise as $N$ grows.
The valid formal product formula for B is the different iteration in
Section~\ref{sec:Bkernel}.

At the Catalan fixed point $\lambda=(1-\sqrt{1-4z})/(2z)$,
the correct boundary identity is
\begin{equation}
 U_B(\lambda)=\frac{b_B+1-\lambda}{z\lambda}.
 \label{tm:eq:fixed}
\end{equation}
For a direct proof, eliminate $P_B$ from the common equations to obtain
\[
 (1-x+zx^2)F_B(x)=zx\{(1-x)(1-zx)U_B(x)+b_B+1-x\}.
\]
Since $1-\lambda+z\lambda^2=0$ and
$(1-\lambda)(1-z\lambda)=-z\lambda$, substitution proves
\eqref{tm:eq:fixed}. All these substitutions are legitimate formal
compositions because $\lambda$ has nonnegative $z$-order and each
coefficient of the original series is polynomial in the catalyst.

\section{A common counting formula for Class 152}\label{sec:class152}
An inversion sequence of length $n$ is $e=(e_0,\ldots,e_{n-1})$ with
$0\le e_i\le i$. A pattern occurs on three increasing indices when both
inequalities and equalities agree with the pattern. Write $\In_n(B_{152})$ for the
sequences avoiding every pattern in $B_{152}$, and put
\[
 A_{152}=\{010,100,102,210\},\qquad B_{152}=\{011,201,210\}.
\]
\begin{theorem}\label{count:class152}
For every $n\ge0$, $|\In_n(A_{152})|=|\In_n(B_{152})|$. The common count is $1$ at $n=0$;
for every $n\ge1$ it is
\[
 a_n=\frac{1+\sum_{j=0}^{n-1}\binom{2j}{j}}2.
\]
\end{theorem}
For $n\ge1$, the number $a_n$ is \textup{OEIS A024718}$(n-1)$~\cite{OEIS}.
That entry already records the partial central-binomial sum and its generating
function
\[
 \sum_{k\ge0}\textup{A024718}(k)z^k
 =\frac{1+(1-4z)^{-1/2}}{2(1-z)},
\]
the latter attributed to Michael Somos (14 February 2006). The new
contribution is the enumeration of the two specified inversion-sequence
avoidance classes and the proof of their Wilf equivalence.
This is Class 152 of Conjecture 1 in Callan--Mansour, \emph{Inversion sequences avoiding quadruple length-3 patterns}~\cite{CM}, Section 1.
We use ordinary length $n$ throughout, including the empty sequence. The triple-avoidance class $B_{152}$ appears as class~73 in Table~7 of Callan, Jel\'inek and Mansour~\cite{CJM}, where its initial counting sequence is recorded. The theorem above gives its closed counting formula and its equality with the quadruple-avoidance class $A_{152}$.

\subsection{The first-descent decomposition}
If a value $a$ follows a larger value, any second occurrence of $a$ gives
$010$ or $100$. Thus $a$ occurs exactly once. If an adjacent descent $M>a$
has a larger value elsewhere, that value forms $210$ when it is earlier and
$102$ when it is later. Consequently every adjacent descent starts at the
global maximum.

A nonmonotone avoider has a unique cut $p\,s$ at its first adjacent descent:
$p$ is a nondecreasing inversion sequence ending in the global maximum $M$,
and $s$ begins below $M$. Its entries other than $M$ are distinct and increasing
(a decrease between two of them gives $210$ with the preceding $M$).
They avoid all values of $p$. Moreover they lie in one gap $(u,v)$ between
successive distinct values of $p$: crossing a prefix value $v$ gives the triple
$v,a,b$ of type $102$. If $v<M$, no $M$ can occur in $s$, since $v,a,M$
would give $102$. Conversely these conditions exclude all four forbidden
patterns. The bounds $e_i\le i$ are automatic in $s$, since its entries are
at most $M\le |p|-1$.

Let $m=|p|$, $r=|s|\ge1$, and let $g_1,\ldots,g_k$ be the positive jumps
between consecutive distinct values of $p$ (starting from $0$).
An internal gap of size $g$ contributes $\binom{g-1}{r}$.
The top gap contributes
\begin{equation}
 \sum_{b\ge1}\binom{g_k-1}{b}\binom{r-1}{b-1}
 =\binom{g_k+r-2}{r}.
\label{c152:gap}
\end{equation}
Here $b$ is the number of smaller entries and the second factor chooses their
positions, with the first position required to be smaller. Empty sums contribute
zero; a gap of size $1$ contributes zero.

Encode $p$ as a Dyck path by requiring exactly $p_i$ down-steps in total before its $(i+1)$st
up-step and completing with down-steps. Nondecreasing inversion sequences
of length $m$ are thereby in bijection with Dyck paths of semilength $m$.
Their positive jumps are precisely the lengths of the nonfinal down-runs;
the top jump is the penultimate down-run. This translates \eqref{c152:gap} into the
weighted path series below.

\subsection{Weighted Dyck paths and the left generating function}
Work in formal power series over $\mathbb Q$. Set
\begin{equation}
 \mathcal C=\sum_{m\ge0}\frac1{m+1}\binom{2m}{m}z^m,
 \qquad t=z\mathcal C,\qquad \mathcal C=1+z\mathcal C^2=\frac1{1-t},\qquad z=t(1-t).
\label{c152:cat}
\end{equation}
For a weight $w$ with $w(0)=0$, let $D_w$ count paths weighted by their final
down-run, $H_w$ count the sum of weights of all nonfinal down-runs, and
$J_w$ count the weight of the penultimate down-run (zero if absent).
All series include the factor $z^m$ for semilength $m$. Then
\begin{equation}
 (2-\mathcal C)H_w=(\mathcal C-1)D_w,\qquad (1-z)J_w=tD_w.
\label{c152:runs}
\end{equation}
Here is a derivation that also checks the marked-run multiplicities.
Reflect a Dyck path to exchange final down-runs and initial up-runs.
A path whose initial up-run has length $k$ has a unique factorization
$U^k D P_1 D P_2\cdots D P_k$, with $k$ arbitrary Dyck paths $P_i$.
Thus paths with that initial run contribute $t^k$. Let $S_w$ count the sum
of weights of \emph{all} up-runs. The factorization gives
\[
 S_w=D_w+S_w\sum_{k\ge1}k z^k \mathcal C^{k-1}
     =D_w+t\mathcal C S_w.
\]
Since $t\mathcal C=\mathcal C-1$ and reflection gives $H_w=S_w-D_w$, the first identity in \eqref{c152:runs}
follows. For the second up-run, select the first nonempty $P_i$; factors
before it are empty, and factors after it are arbitrary. Hence
\[
 J_w=D_w\sum_{k\ge1}z^k\sum_{j=0}^{k-1}\mathcal C^j
     =\frac{z\mathcal C}{1-z}D_w.
\]
Reflection gives the asserted identity for the penultimate down-run.

For $q$ of constant term $1$, use $w_q(\ell)=q^{\ell-1}-1$ for $\ell\ge1$
and $w_q(0)=0$, and write $D_q,H_q,J_q$. The same run factorization yields
\begin{equation}
 D_q=\sum_{\ell\ge1}(q^{\ell-1}-1)t^\ell
     =\frac{\mathcal C t^2(q-1)}{1-qt}.
\label{c152:D}
\end{equation}
With $q_1=1+z$ and $q_2=(1-z)^{-1}$, binomial expansion gives
$[z^r]w_{q_1}(g)=\binom{g-1}{r}$ and
$[z^r]w_{q_2}(g)=\binom{g+r-2}{r}$ for $r\ge1$.
The first-descent decomposition therefore gives the left series
\begin{equation}
 L=\sum_{n\ge0}|\In_n(A_{152})|z^n=\mathcal C+H_{q_1}+J_{q_2}-J_{q_1}.
\label{c152:L}
\end{equation}
In \eqref{c152:L}, $H_{q_1}$ initially weights every nonfinal down-run as an internal
gap; the difference of the two $J$ terms replaces the top-gap weight.
The path's power of $z$ records prefix length, and the weight's power records
suffix length.

Equations \eqref{c152:cat}--\eqref{c152:D} give, directly,
\[
 D_{q_1}=\frac{t^3}{(1-t)(1-t^2)},\qquad
 D_{q_2}=\frac{t^3}{(1-t)^2}=(1+t)D_{q_1},
\]
\begin{equation}
 H_{q_1}+J_{q_2}-J_{q_1}
 =D_{q_1}\left(\frac{t}{1-2t}+\frac{t^2}{1-z}\right)
 =\frac{\mathcal C t^4}{(1-z)(1-2t)}.
\label{c152:simplify}
\end{equation}
All denominators used here have nonzero constant term, so these operations
are valid in formal power series.

\subsection{The last-zero decomposition and the right generating function}
Avoiding $011$ means that all positive entries are distinct, since $e_0=0$.
Cut a nonempty avoider uniquely as $p\,0\,s$ at its last zero, with $m=|p|$
and $M=\max(p)$ (take $M=0$ for empty $p$). Positive entries of $p$ are
increasing: a positive descent followed by the last zero gives $210$.
Every entry of $s$ exceeds $M$: equality is excluded by $011$, and a positive
entry below $M$ gives $M,0,a$ of type $201$. Thus $s$ consists of distinct
entries above $M$, avoids $201,210$, and satisfies $s_j\le m+1+j$.
Conversely these conditions suffice: a forbidden triple crossing the cut
would contradict either the increasing positive prefix or the separation
of prefix and suffix values.

The running-maximum map sends $p$ bijectively to a nondecreasing inversion
sequence: the inverse retains strict increases and replaces other entries
by zero (retaining the first entry $0$). It preserves $M$. Its Dyck path has final down-run length
$h=m-M$; the empty prefix corresponds to the empty path with $h=0$.

A distinct-entry word avoiding $201,210$ has a unique decomposition into
rotated increasing blocks: sort its entries, partition into nonempty
consecutive blocks, and move the minimum of each block to its end.
Indeed, at the position of the global minimum, the preceding entries must
increase (otherwise $210$), and must be below every following entry
(otherwise $201$). Remove this first block and repeat; this also proves
uniqueness and the converse. For example, sorted blocks
$(3,4)\mid(6)$ become the suffix $4,3,6$: the entry $4$ moves one
place to the left and $3$ moves to the end of its block, while the
singleton $6$ stays in place.

For the sorted entries $a_i$, put $\lambda_i=a_i-M-(i+1)$, with $i$ starting
at $0$. These labels are weakly increasing and nonnegative. The bounds on
the rotated word are exactly $0\le\lambda_i\le h$, with labels in block
tails strictly below $h$. Here the \emph{tail} of a block consists of
all entries other than its minimum in the unrotated increasing block;
each such entry moves one place earlier, whereas the minimum moves
to the end of its block. A singleton can have label $h$;
in a longer block its minimum is bounded by its tail labels. For the example above, take $M=1$, $m=3$, so $h=2$.
The labels of $3,4,6$ are $1,1,2$: the tail label $1$ is strictly
below $h$, and the singleton label is $h$. Adding back $M+i+1$
reverses this encoding.

For each label $0,\ldots,h-1$, an occurrence may start a block or continue
one, giving two choices; the first occurrence overall must start a block.
Every occurrence of $h$ must start a singleton. Words containing a lower
label therefore have half the unrestricted binary choices, while words
containing only $h$ have one choice. The suffix series (including the empty
suffix) is consequently
\begin{equation}
 S_h=\frac1{1-z}+\frac{(1-2z)^{-h}-1}{2(1-z)}
     =\frac{1+(1-2z)^{-h}}{2(1-z)}.
\label{c152:suffix}
\end{equation}
The last zero contributes $z$. The final-run factorization in the preceding subsection,
now with weight $r^h$, gives $\sum_P z^{|P|}r^{h(P)}=(1-rt)^{-1}$.
Using $r=(1-2z)^{-1}$ in \eqref{c152:suffix}, we obtain
\begin{equation}
 R=\sum_{n\ge0}|\In_n(B_{152})|z^n
 =1+\frac{z}{2(1-z)}\left(\mathcal C+\frac1{1-t/(1-2z)}\right)
 =1+\frac{z(1-t)}{(1-z)(1-2t)}=L,
\label{c152:R}
\end{equation}
where the last two simplifications use $z=t(1-t)$ and \eqref{c152:simplify}.
Finally $1-2t=\sqrt{1-4z}$, so
\[
 L=R=1+\frac{z}{2(1-z)}\left(1+\frac1{\sqrt{1-4z}}\right).
\]
Extracting coefficients proves the theorem. For $n=0,\ldots,10$, the counts are
\[1,1,2,5,15,50,176,638,2354,8789,33099.\]

\paragraph{Formal verification record}
The formal record uses Lean~4~\cite{Lean} and mathlib~\cite{Mathlib}.
The source theorem \texttt{InversionSeqClass152.result} states the all-length
equality of the finite avoider counts. Its pinned source is
\href{https://github.com/the-omega-institute/trureturing/blob/7e22ca771362e94db46dd4394637c98db0f942f1/D5/S3/Combinatorics/InversionSeq/InversionSeqClass152.lean}{\texttt{InversionSeqClass152.lean}}.
The closed coefficient formula above follows by the displayed generating-function
calculation. The mathematical proof here is complete independently of
this verification record. The supplementary material gives the pinned
Lean and mathlib versions and a command to check this source theorem.

\section{Algebraic and linear differential obstructions}\label{ana:sec:obstructions}

Retain the Class~207 site statistics from Section~\ref{sec:B} and write
$T_C(z;t)=Z_C(z;t,t)$, $U_X(z;t)=Z_X(z;t,1)$, and
$A_B(z;t)=Z_B(z;1,t)$. The notation $F_X(z;x)$ displays the length
variable when it needs emphasis. All derivatives below are with respect
to the stated catalytic argument.
\begin{theorem}\label{ana:thm:main}
Each of the pairs
\[
 \bigl(Z_B(z;x,y),Z_C(z;x,y)\bigr),\qquad
 \bigl(Z_B(z;x,y),Z_C(z;y,y)\bigr)
\]
is algebraically independent over $\CC(z,x,y)$. Explicitly, for
every $A\in\CC(z,x,y)[S,T]$, either identity
\begin{align*}
 A\bigl(Z_B(z;x,y),Z_C(z;x,y)\bigr)&=0,\\
 A\bigl(Z_B(z;x,y),Z_C(z;y,y)\bigr)&=0
\end{align*}
implies $A=0$.
\end{theorem}
There is no restriction on the polynomial degrees or on coefficients
in the stated field.
The proof first fixes $z=q/(1+q)^2$, where $0<q<10^{-4}$, and
establishes algebraic independence of suitable boundary functions over
$\CC(t)$. At this stage the coefficient field includes all complex
constants, including both total-count values. The full-series
statements then follow by two different specializations. The diagonal
assertion needs its own continuation and pole argument; it is not
inferred by specializing the first assertion.

Section~\ref{ana:sec:differential} establishes a further obstruction:
no finite linear relation over $\CC(z,x,y)$ exists among the
$y$-derivatives of $Z_B(z;x,y)$ and $T_C(z;y)$, including a rational
inhomogeneous term. The derivative of $Z_B$ holds $x$ fixed; the
derivative of $T_C$ differentiates the diagonal specialization.
This uses the genuine B poles, rather than algebraic independence
alone: for example, $\exp(1/t)$ has no meromorphic extension at zero
but satisfies $t^2f'+f=0$.

\section{The scalar equations}
Fix real $0<q<10^{-4}$ and put $z=q/(1+q)^2$. Define
\[
 c=\frac{q}{(1-q)^2},\qquad
 \XB(w)=\frac{(1+q)(1-w)}{1-qw},\qquad
 \XC(w)=1-c(w+w^{-1}-2).
\]

Suppress $z$ in the series arguments and write $b_X=F_X(1)$.
The site recurrences proved in Proposition~\ref{prop:common} give
\begin{align}
 P_X(x)&=zx+zxP_X(x)+zx\frac{F_X(x)-b_X}{x-1},\label{ana:eq:P}\\
 F_X(x)&=P_X(x)+zxU_X(x).\label{ana:eq:F}
\end{align}
Their cleared equations for the lower-site series are
\begin{align}
 (1-y+zxy)Z_B(x,y)
 ={}&F_B(x)-P_B(x)
       +zx(1-y)\frac{P_B(x)-P_B(y)}{x-y}\notag\\
 &+zx\left(\frac{F_B(x)-b_B}{x-1}
                  -\frac{F_B(x)-F_B(y)}{x-y}\right),\label{ana:eq:BZ}\\
 \kappa(x,y)Z_C(x,y)
 ={}&zx\left[(x-y)U_C(x)-(1-y)Z_C(y,y)
                   +(1-y)(P_C(x)-P_C(y))\right.\notag\\
 &\hspace{12mm}\left.
             +(x-y)\frac{F_C(x)-b_C}{x-1}-F_C(x)+F_C(y)\right],\label{ana:eq:CZ}\\
 \kappa(x,y)={}&(1-y)(x-y)+zx(x-1).\notag
\end{align}
Divided differences are interpreted after cancellation. A complete
derivation from the legal-site updates is given in Proposition~\ref{prop:common}.
We first give the scalar eliminations used in the proof.

Set
\begin{gather*}
 D(x)=1-x+zx^2,\quad A(x)=D(x)+zx,\quad M(x)=zx+z-1,\\
 \phi(x)=\frac1{1-zx},\qquad\beta_B=G_B-\frac1z.
\end{gather*}
Equation~\eqref{ana:eq:P} gives
\begin{equation}
 F_X(x)=G_X-x-\frac{(1-x)(1-zx)}{zx}P_X(x).\label{ana:eq:FP}
\end{equation}
At $y=\phi(x)$ the left side of~\eqref{ana:eq:BZ} vanishes.
Multiply the right side by $(x-y)/(zx)$ and substitute~\eqref{ana:eq:P}.
Before clearing denominators its vanishing expression is
\[
 (x-y)\frac{F_B(x)-P_B(x)}{zx}
 +(1-y)(P_B(x)-P_B(y))
 +(x-y)\frac{F_B(x)-b_B}{x-1}-F_B(x)+F_B(y).
\]
After substitution it equals $-1/(zx(1-zx)M(x))$ times the
left side minus the right side of
\begin{equation}
 M(x)A(x)F_B(x)+zx(1-zx)^2F_B(\phi(x))
       =xz^2D(x)\beta_B.\label{ana:eq:Bscalar}
\end{equation}
Equivalently, substituting~\eqref{ana:eq:FP} and solving for the other
boundary value gives
\[
 P_B(\phi(x))=\frac{(x-1)A(x)}{z^2x^3}P_B(x)
              +\frac{D(x)(G_B-x)}{zx^2(1-zx)}.
\]

Let $\zeta\in(-2,-1)$ be the root of
\begin{equation}
 q^3\zeta^2+(1-2q-2q^2+q^3)\zeta+1=0.\label{ana:eq:b}
\end{equation}
It is the branch $\zeta=-1+O(q)$. Put $a=q^2\zeta$ and $r=q\zeta$.
Using~\eqref{ana:eq:FP} in~\eqref{ana:eq:Bscalar}, with
$\phi(\XB(w))=\XB(qw)$, gives for $f(w)=P_B(\XB(w))$
\begin{align}
 f(qw)&=R(w)f(w)+B_0(w),\label{ana:eq:Baffine}\\
 R(w)&=\frac{(1-w/q)(1-w/r)(1-aw)}{(1-w)^3},\label{ana:eq:R}\\
 B_0(w)&=\frac{w(q-1)^2(q+1)
   (G_Bqw-G_B-qw+q-w+1)}{q(w-1)^2(q^2w-1)}.\label{ana:eq:B0}
\end{align}
The factorization of $R$ uses~\eqref{ana:eq:b}; all scalar series values
in these equations are constants with respect to $w$.

For C, let $x,x'$ be the two roots of $\kappa(\,\cdot\,,y)$.
At a root, equations~\eqref{ana:eq:P}--\eqref{ana:eq:CZ} imply
\begin{equation}
 (1-y)^2\left(Z_C(y,y)+\frac{P_C(y)}{zy}\right)
 =(1-x)(G_C-x)+\frac{(2-y)(1-x)^2}{1-y}P_C(x).
 \label{ana:eq:Ckernel}
\end{equation}
To see the coefficients explicitly, replace $zxU_C$ by $F_C-P_C$
in~\eqref{ana:eq:CZ}, and then use~\eqref{ana:eq:FP}. The coefficients of
$P_C(y)$ and $P_C(x)$ are respectively $-(1-y)^2/(zy)$ and
$(2-y)(1-x)^2/(1-y)$ after applying the kernel equation.
Subtract~\eqref{ana:eq:Ckernel} for the two roots and use
$x+x'=1+(y-1)/z$. The result is
\begin{equation}
 (2-y)\bigl((1-x)^2P_C(x)-(1-x')^2P_C(x')\bigr)
 =(1-y)(x-x')\left(G_C+\frac{1-y}{z}\right).\label{ana:eq:Ctwo}
\end{equation}
Now take $x=\XC(w)$, $x'=\XC(qw)$ and
\begin{equation}
 Y(w)=1+z(x+x'-1)
 =\frac{(q-w)(q^2w-1)}{w(q-1)^2(q+1)}.\label{ana:eq:Y}
\end{equation}
Direct substitution verifies both kernel equations. The connection with
Proposition~\ref{prop:Ckernel} is
\[
 t=\frac{\XC(w)}{Y(w)}=\XB(qw),\qquad
 T(\XB(qw))=\XC(w),\qquad
 T\!\left(\phi(\XB(qw))\right)=\XC(qw).
\]
Moreover, \eqref{ana:eq:FP} gives
\[
 W_C(x)=\frac{zx(1-x)(G_C-x)}{1-zx}-(1-x)^2P_C(x).
\]
These identities relate the earlier formal reduction to the additive
equation below. With
$h(w)=(1-\XC(w))^2P_C(\XC(w))$, we obtain
\begin{align}
 h(qw)-h(w)&=C_0(w),\label{ana:eq:Cadd}\\
 C_0(w)&=-\frac{(1-Y(w))(x-x')}{2-Y(w)}
                     \left(G_C+\frac{1-Y(w)}z\right).\label{ana:eq:C0}
\end{align}
The kernel calculations hold on nonempty open subsets where the
displayed denominators do not vanish, and the resulting identities
extend meromorphically. No identically zero factor is divided out.
This derivation needs no infinite Lambert-series formula.

\section{Convergence and continuation}
Repeated $w\mapsto qw$ need not act on a formal $q$-series with
coefficients in $\QQ(w)$. We instead use actual meromorphic functions
for a fixed small $q$.

\begin{lemma}\label{ana:lem:analytic}
For $0<q<10^{-4}$, $f$ and $h$ extend to single-valued meromorphic
functions on $\CC^*$ satisfying~\eqref{ana:eq:Baffine} and~\eqref{ana:eq:Cadd}.
Also $1\le G_X(z)<2$, and $G_X(z)\ne1/z$.
\end{lemma}
\begin{proof}
Every word here avoids a strictly decreasing subsequence of length
three. Color each entry by the length of a longest strictly decreasing
subsequence ending there. The colors are $1,2$, and each color subword
is weakly increasing. A weakly increasing subword of length at most $n$
with values in $\{0,\ldots,n-1\}$ has at most
$\binom{2n}{n}\le4^n$ possibilities. With at most $2^n$ colorings,
there are at most $32^n$ words. Overcounting is harmless.
Thus $F_X$ and $P_X$ converge absolutely if
$32|z|\max(1,|x|)<1$. For a lower site, $u\le M-1$, so
$s+u\le n-1$, and there are at most $n$ sites per word.
Hence $Z_X$ converges if $32|z|\max(1,|x|,|y|)<1$ as well.
The coefficient recurrences hold in these convergence domains.
In particular,
\begin{equation}
 1\le G_X(z)\le(1-32z)^{-1}<2<1/z.\label{ana:eq:Gbound}
\end{equation}

For B, both $\XB(w)$ and $\XB(qw)$ are bounded near $w=0$.
The convergent germ therefore satisfies~\eqref{ana:eq:Baffine} there.
For arbitrary $w\ne0$, choose $m$ so large that $q^mw$ is in this
disk and apply the affine equation backwards $m$ times.
Its coefficients are rational and $R\ne0$, so this gives a
meromorphic function locally. Different choices of $m$ agree by the
equation in the original disk, giving a single-valued continuation.

For C, put $\varepsilon=128q^2$ and
$\mathcal D=\{\varepsilon<|w|<\varepsilon^{-1}\}$.
Since $c<2q$ and $|w|+|w|^{-1}\le\varepsilon+\varepsilon^{-1}$, throughout
this annulus
\[
 z\max(1,|\XC(w)|)
 \le q+\frac1{64}+4q^2+256q^4<\frac1{32}.
\]
Thus $h$ is analytic on $\mathcal D$.
For $|w-1|<1/10$ and $q<10^{-4}$, the displayed rational formulas
give $|x|,|x'|,|Y|<4$, so every series used in
deriving~\eqref{ana:eq:Ctwo} converges. Equation~\eqref{ana:eq:Cadd}
holds on an open subset of the connected annulus
$\mathcal D\cap q^{-1}\mathcal D$; the meromorphic identity theorem
extends it to that overlap. The annuli $q^m\mathcal D$, $m\in\mathbb Z$,
cover $\CC^*$ and adjacent ones overlap. The additive equation
extends $h$ between them, consistently on overlaps, producing the
required single-valued meromorphic continuation. For consistency on
multiple overlaps, if $w$ and $q^k w$ lie in $\mathcal D$, every
intermediate radial shift does too. Iterating the identity on adjacent
overlaps therefore gives the same extension on every overlap.
\end{proof}

\section{The B singularity}
We first establish the property of B that will distinguish it from
both the C boundary and the C diagonal.
Write $K=\CC(w)$, $\s v(w)=v(qw)$, and
$[\alpha]=\{q^m\alpha:m\in\mathbb Z\}$.
\begin{lemma}\label{ana:lem:trace}
Let $u$ be a meromorphic function on $\CC^*$ satisfying
$\s u=A u+B$, with $A\in K^*$ and $B\in K$.
If that equation has no rational solution, then $u$ is transcendental
over $K$.
\end{lemma}
\begin{proof}
If $u$ were algebraic, $L=K(u)$ would be a finite extension of $K$.
The affine equation and its inverse
\[
 \s^{-1}u=\frac{u-\s^{-1}B}{\s^{-1}A}
\]
show that $L$ is stable under both $\s$ and $\s^{-1}$.
The normalized trace
\[
 v=\frac{\operatorname{Tr}_{L/K}(u)}{[L:K]}\in K
\]
satisfies $\s v=Av+B$, contradicting the absence of a rational
solution. Indeed, trace commutes with the semilinear automorphism
$\s$: apply $\s$ to the multiplication matrix in a basis, using
the transformed basis to compute the trace of $\s u$.
\end{proof}

\begin{lemma}\label{ana:lem:transcendence}
The function $f(w)=P_B(z;\XB(w))$ is transcendental over $\CC(w)$.
\end{lemma}
\begin{proof}
$B_0$ has a genuine simple pole at $w=q^{-2}$: its
parenthesized numerator there is
$(1-q)(G_B-1/z)/q\ne0$.
Its only other possible finite nonzero pole is $1$, of order at
most two. On the positive orbit $[1]$, $R$ has a simple zero
at $q$, a pole of order three at $1$, and is holomorphic and
nonzero (a unit) elsewhere on this orbit.
Suppose a rational $g$ solves~\eqref{ana:eq:Baffine}. It must have a
pole on $[1]$ to match the pole of $B_0$ at $q^{-2}$.
If $m$ is the least exponent at which $g$ has a pole, the pole
of $g(qw)$ at $q^{m-1}$ must come from $R$ or $B_0$, so
$m-1\in\{0,-2\}$. If $m=1$, both $g(qw)$ and $R(w)g(w)$
would be regular at $w=q^{-2}$, leaving the genuine forcing pole
unmatched. Hence $m=-1$.
At $w=q^{-2}$, $g$ is regular and $R$ is a unit, so the simple
forcing pole makes $g(qw)$ have order one there. Successive
applications of the equation then force
\[
 \begin{array}{c|cccc}
 \text{pole of }g & q^{-1} & 1 & q & q^k\ (k\ge2)\\ \hline
 \text{pole order} & 1 & 1 & 4 & 3
 \end{array}
\]
At $q^{-1}$ both $R$ and $B_0$ are regular and $R$ is a unit,
giving the pole at $1$. At $1$ the order-four product $Rg$
cannot cancel against the forcing pole of order at most two.
At $q$ the simple zero of $R$ reduces the order to three,
which then propagates indefinitely because $R$ is a unit and
$B_0$ is regular at every $q^k$ with $k\ge2$.
This contradicts rationality. Lemma~\ref{ana:lem:trace} now proves
that $f$ is transcendental.
\end{proof}

\begin{lemma}\label{ana:lem:local}
Let $f$ be a single-valued meromorphic function on a punctured
neighborhood of $t_*$ with no meromorphic extension to $t_*$.
Let $g$ be meromorphic at $t_*$ and transcendental over $\CC(t)$
as a germ. Then $f,g$ are algebraically independent over $\CC(t)$.
\end{lemma}
\begin{proof}
Suppose $\sum_{j=0}^d a_j(t,g(t))f(t)^j=0$, where
$a_j\in\CC(t)[T]$ and $a_d\ne0$. If $d=0$, this contradicts
the transcendence of $g$. Otherwise that transcendence makes
$a_d(t,g(t))$ a nonzero meromorphic germ. Divide by it to obtain
a monic polynomial in $f$ whose coefficients $c_j(t)$ are
meromorphic at $t_*$. On a sufficiently small punctured disk
they are holomorphic and bounded by $C|t-t_*|^{-m}$ for a common
integer $m\ge0$. The monic equation precludes poles of $f$
in this punctured disk: at such a pole its highest power would
have strictly greater pole order than every other term.
The polynomial root bound
$|f|\le1+\max_{j<d}|c_j|$ now gives
$|f(t)|\le C'|t-t_*|^{-m}$.
Thus $(t-t_*)^m f(t)$ has a removable singularity, and $f$
extends meromorphically to $t_*$, a contradiction.
\end{proof}

\begin{lemma}\label{ana:lem:Bpoint}
At $t_*=(1+q)/q$, the boundary function $P_B(z;t)$ has a single-valued
meromorphic continuation to a punctured neighborhood but has no
meromorphic extension to $t_*$.
\end{lemma}
\begin{proof}
The B continuation in Lemma~\ref{ana:lem:analytic} is analytic at
$w=0$, since it begins with the actual convergent series there.
It is meromorphic on $\CC^*$ and transcendental by
Lemma~\ref{ana:lem:transcendence}. It therefore cannot be meromorphic
at $w=\infty$: otherwise it would be meromorphic on the Riemann
sphere and hence rational. The M\"obius map $\XB$ has
\[
 \XB(\infty)=t_*:=\frac{1+q}{q},\qquad
 w=\frac{1+q-t}{1+q-qt}.
\]
It follows that $P_B(z;t)$ has a single-valued meromorphic
continuation to a punctured neighborhood of $t_*$ with no
meromorphic extension at $t_*$.

\end{proof}

\section{The exact C boundary poles}
For brevity put $P(t)=P_C(z;t)$ and $G=G_C(z)$.
\begin{lemma}\label{ana:lem:Cpoles}
The C boundary $P$ extends meromorphically to the whole finite
$t$-plane. Its poles are exactly the simple poles
\begin{equation}
 t_j=\XC(aq^j)=1+2c-ca q^j-\frac ca q^{-j},\qquad j\ge0.
 \label{ana:eq:tj}
\end{equation}
They are strictly increasing positive real numbers tending to infinity.
In particular, $P$ is transcendental over $\CC(t)$.
\end{lemma}
\begin{proof}
The polynomial defining $\zeta$ is positive at $-1$ and negative at $-2$,
and its derivative on $[-2,-1]$ is positive. Thus its small branch
lies in $(-2,-1)$. The negative $q$-orbits of $\zeta$ and $\zeta^{-1}$
are distinct: equality would imply $\zeta^2=q^m$ for an integer $m$,
whereas $1<\zeta^2<4<q^{-1}$.
The denominator in the C forcing is
\[
 2-Y(w)=\frac{q^2w^2+(q^3-2q^2-2q+1)w+q}
                  {w(q-1)^2(q+1)}.
\]
Its two simple zeros are $r=q\zeta=a/q$ and $a^{-1}$.
These are the only possible finite nonzero poles of $C_0$.
At either zero, the numerator contains the nonzero factor
$G-1/z$, by~\eqref{ana:eq:Gbound}.
Also $x-x'\ne0$, since the roots of $\kappa(x,2)$ have
discriminant $1-6z+z^2>0$. Both forcing poles are therefore genuine
and simple. This elementary pole propagation is consistent with the
orbit-residue viewpoint of rational $q$-summability%
~\cite[Section 2.3, Proposition 2.10]{CS};
we give the needed continuation argument explicitly.

The function $h$ is analytic on
$\mathcal D=\{128q^2<|w|<(128q^2)^{-1}\}$ by
Lemma~\ref{ana:lem:analytic}. In particular it is analytic at $r$ and
$r^{-1}$, since $q<|r|<2q$.
Put $R_0=\operatorname{Res}_{w=r}C_0(w)\ne0$.
The equation $h(qw)=h(w)+C_0(w)$ first creates a simple pole
at $a=qr$, with residue $qR_0$. At $a,aq,aq^2,\ldots$ the
forcing is regular, so the pole propagates without cancellation:
\begin{equation}
 \operatorname{Res}_{w=aq^j}h(w)=q^{j+1}R_0\ne0
 \qquad(j\ge0).
 \label{ana:eq:hresidue}
\end{equation}
There is no second forcing pole on this orbit.

On the original annulus, $h(1/w)=h(w)$; the meromorphic identity
theorem extends this equality to $\CC^*$. Hence the reciprocal
points $(aq^j)^{-1}$ are also simple poles.
There are no other poles. Indeed, a point on any other $q$-orbit
can be shifted into $\mathcal D$, and continuation along that orbit
encounters only regular forcing terms. On the orbit of $\zeta$, points
$q^m \zeta$ with $m\le1$ remain regular: start at $r$ and continue
outwards with $h(w)=h(qw)-C_0(w)$, encountering no forcing pole.
The points with $m\ge2$ are exactly the poles already listed.
Inversion gives the complementary statement on the orbit of $\zeta^{-1}$.
This proves that the complete pole set of $h$ is
\[
 \{aq^j,(aq^j)^{-1}:j\ge0\}.
\]

Now $p(w)=h(w)/(1-\XC(w))^2$ is meromorphic on $\CC^*$ and
invariant under $w\mapsto1/w$. Its possible additional singularity
at $w=1$ is removable, because it agrees there with the original
convergent $P_C(\XC(w))$. The map $\XC:\CC^*\to\CC$ is the
proper degree-two quotient by this involution, with fibers given by
\[
 w^2-\left(2+\frac{1-t}{c}\right)w+1=0.
\]
The invariant meromorphic function descends to a meromorphic function
of $t$. At the branch points $w=1,-1$, use respectively
$(w-1)/(w+1)$ and $(w+1)/(w-1)$ as local coordinates. Inversion
negates the coordinate, so an invariant Laurent series contains only
even powers and is meromorphic in $t$. Thus descent holds at both
branch values as well. The descended function is the continuation of $P$.
The poles of $h$ are away from these branch points and from $\XC=1$,
so they give exactly the simple poles~\eqref{ana:eq:tj}.
Their residues are
\[
 \operatorname{Res}_{t=t_j}P(t)
 =\frac{\XC'(aq^j)\,q^{j+1}R_0}{(1-t_j)^2}\ne0.
\]
Writing $d=-a\in(q^2,2q^2)$ gives
$t_j=1+2c+c d q^j+c d^{-1}q^{-j}$, and
\[
 t_{j+1}-t_j=c(1-q)
       \left(d^{-1}q^{-j-1}-d q^j\right)>0.
\]
The same formula shows $t_j\to+\infty$.
Finally, a meromorphic function on $\CC$ algebraic over $\CC(t)$
has only finitely many finite poles: away from the finitely many
poles of the rational coefficients of a monic equation, its highest
power cannot cancel at a pole. Thus $P$ is transcendental.
\end{proof}

\section{Meromorphic continuation and poles of the diagonal}
Define $t_{-1}=\XC(a/q)$ and $y_j=Y(aq^j)$ for $j\ge-1$.
\begin{lemma}\label{ana:lem:interlace}
We have $y_{-1}=2$, $t_0>2$, and
\[
 t_{j-1}<y_j<t_j\qquad(j\ge0).
\]
For each boundary pole $x=t_j$, the two roots in $y$ of
$\kappa(t_j,y)=0$ are $y_j$ and $y_{j-1}$.
Consequently none of these root values is any $t_k$, $k\ge0$.
\end{lemma}
\begin{proof}
The two kernel identities at $x=\XC(w),\XC(qw)$ and $y=Y(w)$
were checked in deriving~\eqref{ana:eq:Cadd}.
Since $\XC(w)=\XC(w^{-1})$, the same identities with $w/q$
show that the two $y$ roots at $x=\XC(w)$ are $Y(w),Y(w/q)$.
They are distinct at $w=aq^j$: subtracting their explicit formulas
gives
\[
 Y(w)-Y(w/q)=\frac{q(w^2-1)}{w(1-q)(1+q)}>0
 \quad\text{for } -2q^2<w<0.
\]
Also $2-Y(r)=0$, so $y_{-1}=2$, and
$t_0>c/(2q^2)>1/(2q)>2$.

For the inequalities write $w=aq^j=-v$, with $0<v<2q^2$.
Direct subtraction gives
\[
 \XC(w)-Y(w)=\frac{q(q+v)(v+1)}{v(1-q)^2(1+q)}>0,
\]
and $Y(w)-\XC(w/q)=-N/[v(1-q)^2(1+q)]$, where
\[
 N=-q+q^2+q^3+(q+q^2)v+(1+q-q^2)v^2
 \le -q+q^2+3q^3+6q^4+4q^5<0.
\]
The last inequality follows at once from $0<q<10^{-4}$.
Together with the strict increase of the $t_j$, these gaps place
each $y_j$ strictly between consecutive $t$ values. The remaining
root value $y_{-1}=2$ is below $t_0$, proving the assertion.
\end{proof}

\begin{lemma}\label{ana:lem:Tmeromorphic}
The actual diagonal $T_C(z;t)$ extends meromorphically to every
finite $t$. It has a simple pole at every $t_j$, $j\ge0$, and is
transcendental over $\CC(t)$.
\end{lemma}
\begin{proof}
At a C kernel root,~\eqref{ana:eq:Ckernel} gives
\begin{equation}
 T_C(z;y)= -\frac{P(y)}{zy}
 +\frac{(1-x)(G-x)}{(1-y)^2}
 +\frac{(2-y)(1-x)^2P(x)}{(1-y)^3}.
 \label{ana:eq:translator}
\end{equation}
We justify its global single-valued meaning, including possible
branch points. Parameterize the kernel by $x=\XC(w)$ and
$y=Y(w)$. The proper degree-two map $Y:\CC^*\to\CC$ is the
quotient by $w\mapsto1/(qw)$. This involution exchanges
$x=\XC(w)$ with $x'=\XC(qw)$.
The right side of~\eqref{ana:eq:translator} is meromorphic on $\CC^*$,
by Lemma~\ref{ana:lem:Cpoles}, and~\eqref{ana:eq:Ctwo} says it is
invariant under this involution. The identity holds initially on an
open convergent set, and extends to $\CC^*$ by the identity theorem.
It therefore descends meromorphically to the finite $y$-plane.
At either branch point $w_0=\pm q^{-1/2}$, the local coordinate
\[
 u=\frac{w-w_0}{w+w_0}
\]
is negated by the involution. The invariant Laurent expansion
contains only even powers, which proves meromorphic descent there.
The descended function agrees with the original
diagonal on the initial convergent set, hence is its continuation.

Fix $y=t_k$. Both kernel roots $x$ are finite because the
coefficient of $x^2$ is the fixed nonzero number $z$.
Neither is a pole $t_j$ of $P$: such a pole would force
$y=y_j$ or $y_{j-1}$, excluded by Lemma~\ref{ana:lem:interlace}.
Thus the last two terms of~\eqref{ana:eq:translator} are regular
at $t_k$. To include a possible colliding root explicitly, their
invariant holomorphic pullback descends holomorphically; no pole
is introduced at a branch point.
The denominators $1-y$ are nonzero since $t_k>2$.
The first term has a genuine simple pole, with residue
\[
 \mathop{\rm Res}_{y=t_k}T_C(z;y)
 =-\frac1{z t_k}\mathop{\rm Res}_{y=t_k}P(y)\ne0.
\]
There are infinitely many such poles, so the finite-pole argument
at the end of Lemma~\ref{ana:lem:Cpoles} proves transcendence.
Other poles of $T_C$ are not excluded by this argument.
\end{proof}

\section{Algebraic independence of both pairs}
\begin{proof}[Proof of Theorem~\ref{ana:thm:main}]
By Lemma~\ref{ana:lem:Bpoint}, $P_B$ has no meromorphic extension at
$t_*=(1+q)/q$. Lemmas~\ref{ana:lem:Cpoles} and~\ref{ana:lem:Tmeromorphic}
make both $P_C$ and $T_C$ meromorphic there and transcendental
over $\CC(t)$. Their germs there are also transcendental: any
polynomial relation on one such germ would extend over the connected
plane by the identity theorem. Lemma~\ref{ana:lem:local} proves that
each of the pairs
\[
 (P_B(z;t),P_C(z;t)),\qquad (P_B(z;t),T_C(z;t))
\]
is algebraically independent over $\CC(t)$.

Two rational affine changes give the other boundary statements.
Equations~\eqref{ana:eq:F} and~\eqref{ana:eq:FP} imply
\begin{equation}
 z^2t^2U_X(z;t)=zt(G_X-t)-(1-t+zt^2)P_X(z;t),
 \qquad X=B,C.
 \label{ana:eq:transfer}
\end{equation}
Specializing~\eqref{ana:eq:BZ} at $x=1$ and using~\eqref{ana:eq:FP} gives
\begin{equation}
 (1-t+zt)A_B(z;t)=K_B-\frac{P_B(z;t)}t,
 \quad K_B=F_B(1)+(z-1)P_B(1)+zF'_B(1)+z.
 \label{ana:eq:Atransfer}
\end{equation}
Here the prime differentiates the catalytic variable and all scalar series
are evaluated at the fixed $z$.
For clarity, the right side of~\eqref{ana:eq:BZ} at $x=1$ before
substitution is
\[
 F_B(1)-P_B(1)+z(P_B(1)-P_B(t))
 +zF'_B(1)-z\frac{F_B(1)-F_B(t)}{1-t}.
\]
This becomes~\eqref{ana:eq:Atransfer}. Both changes are invertible
over $\CC(t)$, since their multipliers are nonzero rational
functions. They preserve algebraic independence and prove the
following fixed-$z$ independence statements:
\[
 (U_B(z;t),U_C(z;t)),\qquad
 (U_B(z;t),T_C(z;t)),\qquad
 (A_B(z;t),T_C(z;t)).
\]

For the first full-series assertion, suppose a nonzero polynomial
relation exists between $Z_B(z;x,y)$ and $Z_C(z;x,y)$.
Clear rational denominators to obtain
\[
 A=\sum_{i,j}a_{ij}(z,x,y)S^iT^j\in\CC[z,x,y][S,T].
\]
Divide all coefficients $a_{ij}$ by their largest common power of
$y-1$. This preserves the relation, since cancellation is valid in
the integral domain $\CC[x,y][[z]]$. At least one resulting
coefficient is not divisible by $y-1$, so the specialization
$\overline A=A|_{y=1}\in\CC[z,x][S,T]$ is nonzero and gives
\[
 \overline A\bigl(z,x,U_B(z;x),U_C(z;x)\bigr)=0.
\]
At all but finitely many positive $z=q/(1+q)^2$ with $0<q<10^{-4}$,
this remains a nonzero polynomial in $x,S,T$: one of its nonzero
coefficients, viewed as a polynomial in $z$, has only finitely many
roots. Choose such a $z$. Absolute convergence in
Lemma~\ref{ana:lem:analytic} turns the formal identity into an analytic
identity near $x=1$. Continue it to a punctured neighborhood of $t_*$
using the continuations already constructed. This contradicts the
fixed-$z$ independence of $U_B,U_C$.

For the second assertion, start with a nonzero relation between
$Z_B(z;x,y)$ and $T_C(z;y)$. Clear denominators as above, but now
divide all coefficients by their largest common power of $x-1$.
Cancellation again preserves the identity. Specializing at $x=1$
leaves a nonzero polynomial $\widehat A\in\CC[z,y][S,T]$ and gives
\[
 \widehat A\bigl(z,y,A_B(z;y),T_C(z;y)\bigr)=0.
\]
Choose an admissible $z$ outside the finite exceptional set of one
nonzero coefficient. Convergence near $y=1$ and continuation then
contradict the fixed-$z$ independence of $A_B,T_C$.
This proves both assertions.
\end{proof}

\section{Linear differential relations in the catalytic variable}
\label{ana:sec:differential}
\begin{theorem}\label{ana:thm:differential}
Let $r,s\ge0$. If $a_i,c_j,d\in\CC(z,x,y)$ satisfy
\begin{equation}
 \sum_{i=0}^{r}a_i(z,x,y)\,\partial_y^i Z_B(z;x,y)
 +\sum_{j=0}^{s}c_j(z,x,y)\,\frac{d^j}{dy^j}T_C(z;y)
 +d(z,x,y)=0
 \label{ana:eq:differential}
\end{equation}
as a formal identity in $z$, then all $a_i,c_j,d$ are zero.
Here $\partial_y$ holds $x$ fixed, whereas $d/dy$ differentiates
the complete diagonal $T_C(z;y)=Z_C(z;y,y)$.
The same assertion holds over $\CC(z,t)$ when $Z_B$ is replaced
by any of $P_B(z;t)$, $U_B(z;t)$, or $A_B(z;t)$, and $T_C(z;y)$
is replaced by $T_C(z;t)$.
\end{theorem}
There is no bound on the finite derivative orders or on the rational
coefficient degrees. The proof strengthens the B singularity argument
by locating simple poles of the actual continued function.

\subsection{Simple B poles accumulating at a finite point}
\begin{lemma}\label{ana:lem:Bpoles}
For each fixed $0<q<10^{-4}$, the actual B boundary $P_B(z;t)$
has a genuine simple pole at every
\begin{equation}
 \tau_n=\frac{1+q}{q}\frac{1-q^n}{1-q^{n-1}},\qquad n\ge2.
 \label{ana:eq:Bpoles}
\end{equation}
These are distinct, strictly decreasing, and converge to
$t_*=(1+q)/q$.
\end{lemma}
\begin{proof}
Use the actual continuation $f(w)=P_B(z;\XB(w))$ from
Lemma~\ref{ana:lem:analytic}. The inverse M\"obius coordinate in
Lemma~\ref{ana:lem:Bpoint}, together with the original germ at $w=0$,
shows that $P_B(z;t)$ is meromorphic on $\CC\setminus\{t_*\}$.
The function $f$ is analytic at $w=1$ because
$\XB(1)=0$ lies in the original convergence domain. The rational
function $R$ is a unit at every positive $q^{-n}$, $n\ge1$:
its only positive zero is $q$, and its only positive pole is $1$.
The forcing $B_0$ is regular at $q^{-1}$. Therefore
\[
 f(w)=\frac{f(qw)-B_0(w)}{R(w)}
\]
shows that $f$ is analytic at $q^{-1}$ as well.

At $w=q^{-2}$, $B_0$ has a genuine simple pole. Indeed its
parenthesized numerator in~\eqref{ana:eq:B0} is
\[
 \frac{1-q}{q}\left(G_B-\frac1z\right)\ne0
\]
there, by~\eqref{ana:eq:Gbound}. The other numerator factors are
nonzero, and $q^2w-1$ has a simple zero.
As $f(qw)$ is analytic at this point and $R$ is a unit, $f$
has a genuine simple pole at $q^{-2}$.
At each $q^{-n}$, $n\ge3$, the forcing is regular and $R$ is
a unit. A genuine simple pole of $f$ at $q^{-(n-1)}$
therefore propagates to a genuine simple pole at $q^{-n}$.
More explicitly, writing $\eta_n=\operatorname{Res}_{w=q^{-n}}f(w)$,
\[
 \eta_2=-\frac{\operatorname{Res}_{w=q^{-2}}B_0(w)}{R(q^{-2})}\ne0,
 \qquad
 \eta_n=\frac{\eta_{n-1}}{qR(q^{-n})}\ne0\quad(n\ge3).
\]

The map $\XB$ is locally invertible at these points, since
\[
 \XB'(w)=\frac{(1+q)(q-1)}{(1-qw)^2}\ne0.
\]
Its images are exactly~\eqref{ana:eq:Bpoles}, with residues
$\XB'(q^{-n})\eta_n\ne0$. Finally,
\[
 \tau_n-t_*=\frac{(1+q)(1-q)q^{n-2}}{1-q^{n-1}}>0,
 \qquad
 \frac{\tau_{n+1}-t_*}{\tau_n-t_*}
 =q\frac{1-q^{n-1}}{1-q^n}<1.
\]
This proves the assertions.
\end{proof}

\subsection{Separation by pole order}
\begin{lemma}\label{ana:lem:pole-separation}
Let $f$ be meromorphic on $\CC\setminus\{t_*\}$, with distinct
simple poles $\tau_n\to t_*\in\CC$.
Suppose $g$ is meromorphic on the finite plane and has infinitely
many distinct simple poles there. Then, for arbitrary finite $r,s$,
the functions
\[
 1,f,f',\ldots,f^{(r)},g,g',\ldots,g^{(s)}
\]
are linearly independent over $\CC(t)$.
\end{lemma}
\begin{proof}
Consider $\sum_i a_i f^{(i)}+\sum_j c_j g^{(j)}+d=0$ with
rational coefficients. If some $a_i$ is nonzero, choose the largest
such index $i=r_0$. A sufficiently small punctured disk about $t_*$
contains no poles of $g$ or of any coefficient: meromorphicity at
$t_*$ and rationality give this after excluding the center itself.
Also a nonzero rational $a_{r_0}$ has finitely many zeros.
Choose $\tau_n$ in this disk where $a_{r_0}(\tau_n)\ne0$.
If the residue of $f$ there is $\eta\ne0$, then
\[
 a_{r_0}(t)f^{(r_0)}(t)
 =\frac{a_{r_0}(\tau_n)(-1)^{r_0}r_0!\eta}
        {(t-\tau_n)^{r_0+1}}
   +O\bigl((t-\tau_n)^{-r_0}\bigr).
\]
Every smaller derivative has lower pole order, and every $g$ term
and $d$ is regular. The leading principal part cannot cancel.
Hence all $a_i$ vanish.

If some $c_j$ is nonzero, take the largest index $j=s_0$.
Among the infinitely many simple poles of $g$, choose one at which
all coefficients are regular and $c_{s_0}$ is nonzero. The same
highest-pole-order argument gives a contradiction. Thus all $c_j$
vanish, and then $d=0$.
\end{proof}

\begin{proof}[Proof of Theorem~\ref{ana:thm:differential}]
First fix $z=q/(1+q)^2$ with $0<q<10^{-4}$.
Lemma~\ref{ana:lem:Bpoles} supplies the B simple poles accumulating
at $t_*$. Lemma~\ref{ana:lem:Tmeromorphic} makes $T_C$ meromorphic
at every finite point, including $t_*$, and supplies infinitely
many genuine simple poles tending to infinity.
Lemma~\ref{ana:lem:pole-separation} proves the assertion for
$P_B(z;t),T_C(z;t)$ over $\CC(t)$.

The invertible rational affine changes
\eqref{ana:eq:transfer} and~\eqref{ana:eq:Atransfer} give the other B series.
If $V=\alpha P_B+\beta$, with rational $\alpha\ne0,\beta$,
then $V^{(i)}=\alpha P_B^{(i)}$ plus lower derivatives and a
rational term. Thus this change is triangular and invertible
on each finite derivative span with the constant function included.
It transfers linear independence to $U_B,T_C$ and $A_B,T_C$.

Now suppose~\eqref{ana:eq:differential} has some nonzero coefficient.
Clear rational denominators to make all coefficients polynomials
in $z,x,y$. Divide them by their largest common power of $x-1$.
Cancellation preserves the identity in the integral domain
$\CC[x,y][[z]]$, and at least one resulting coefficient is
nonzero at $x=1$. Specialize $x=1$. Differentiation in $y$
commutes with this specialization, so a nonzero polynomial
coefficient relation between derivatives of $A_B(z;y)$ and
$T_C(z;y)$ results.

Choose $z=q/(1+q)^2$, $0<q<10^{-4}$, outside the finitely many
zeros of a nonzero polynomial in $z$ obtained as the coefficient
of a power of $y$ in one of the surviving coefficient polynomials.
The specialized relation remains nonzero. Absolute convergence near $y=1$,
including termwise differentiation to each fixed finite order,
turns the formal identity into an analytic identity. It continues
to the domains used above, contradicting the fixed-$z$ derivative
independence of $A_B,T_C$. The same argument without the $x=1$
step proves the stated assertions for functions of one catalytic variable over $\CC(z,t)$.
\end{proof}

\section{Polynomial-exponential expressions}
\label{exp:sec:exponential}
Fix $0<q<10^{-4}$, $z=q/(1+q)^2$, and $t_*=(1+q)/q$ as in the
preceding analytic arguments. The next result concerns each such fixed $z$.
Let $K$ be the field of germs of meromorphic functions at $t_*$.
An element of $K$ is holomorphic on a sufficiently small punctured disk
and has at most a pole at its center. The letter $F$ below is an
independent indeterminate. The actual B input $b(t)$ is any one of
\begin{equation}
 P_B(z;t),\qquad U_B(z;t)=Z_B(z;t,1),\qquad
 Z_B(z;x_0,t),\quad x_0\in(1/2,3/2).
 \label{exp:eq:inputs}
\end{equation}
Here $x_0$ is fixed independently of $t$.

\begin{theorem}\label{exp:thm:main}
Let $m\ge1$, $R_0,\ldots,R_m\in K(F)$, and
$P_1,\ldots,P_m\in K[F]$. Suppose
\[
 \deg_F P_j\ge1,\qquad
 \deg_F(P_i-P_j)\ge1\quad(i\ne j).
\]
For any input $b$ in \eqref{exp:eq:inputs}, an identity
\begin{equation}
 T_C(z;t)=R_0(t,b(t))+
       \sum_{j=1}^m R_j(t,b(t))\exp P_j(t,b(t))
 \label{exp:eq:excluded}
\end{equation}
on a nonempty open part of a sufficiently small punctured neighborhood
of $t_*$, after meromorphic continuation of the original series, forces
\[
 R_1=\cdots=R_m=0,\qquad R_0\in K,
 \qquad R_0(t)=T_C(z;t).
\]
Thus no surviving exponential term in \eqref{exp:eq:excluded} can depend
nontrivially and polynomially on the selected B input.
\end{theorem}

In particular, for every nonzero $R\in K(F)$ and every polynomial
$P\in K[F]$ of positive $F$-degree,
\[
 T_C(z;t)\ne R(t,b(t))\exp P(t,b(t)).
\]
There is no bound on the degree of $P$, on the rational prefactor,
or on the number of separated exponential terms. Finite Laurent
polynomials in $\exp P(t,b(t))$ are also covered by taking the
distinct exponents $jP$, including negative integers $j$.
The theorem permits arbitrary coefficients meromorphic in the catalytic variable,
which is stronger than restricting them to rational functions of $t$.
Its conclusion explicitly identifies the B-independent remainder.

Identical polynomials in an exponent can first be merged by adding their
rational prefactors. The separation hypothesis is literal: two exponents
whose difference is independent of $F$ are not separate terms under this
statement. A constant shift holomorphic at $t_*$ can be absorbed into
a prefactor in $K(F)$; an arbitrary meromorphic shift cannot always be
absorbed in this way.
\section{Exponential independence at accumulating poles}
We give the local analytic argument in full.
\begin{lemma}\label{exp:lem:localexp}
Let $h_1,\ldots,h_m$ be meromorphic germs at a point $\tau$.
If every difference $h_i-h_j$, $i\ne j$, has a pole there,
then $e^{h_1},\ldots,e^{h_m}$ are linearly independent over
the field of meromorphic germs at $\tau$.
\end{lemma}
\begin{proof}
A relation with one nonzero term is impossible. Otherwise take
a relation with the fewest nonzero terms, and divide by its last term
to obtain
\[
 1+\sum_{j=1}^{r-1}b_j e^{h_j-h_r}=0,
 \qquad b_j\ne0.
\]
Differentiating gives a relation among at most $r-1$ of the
original exponentials, after multiplication by $e^{h_r}$.
Minimality implies
$b'_j+b_j(h_j-h_r)'=0$ for every $j$.
A nonzero meromorphic germ is $(t-\tau)^k$ times an analytic
unit. Its logarithmic derivative has at most a simple pole.
But the derivative of a genuine pole has order at least two.
Thus $b'_j/b_j=-(h_j-h_r)'$ is impossible.
\end{proof}

\begin{lemma}\label{exp:lem:independence}
Let $f(t)=P_B(z;t)$. It is transcendental over $K$.
If $Q_0,\ldots,Q_m\in K[F]$ have pairwise differences of
positive $F$-degree, then
$\exp Q_0(t,f(t)),\ldots,\exp Q_m(t,f(t))$
are linearly independent over $K(f)$.
\end{lemma}
\begin{proof}
Suppose $\sum_{j=0}^d a_j(t)f(t)^j=0$ is a nonzero polynomial
relation over $K$, with $a_d\ne0$. For all sufficiently large
$n$, each $a_j$ is regular at $\tau_n$, and $a_d(\tau_n)\ne0$:
a nonzero meromorphic germ has no zeros or poles in a sufficiently
small punctured disk. At the simple pole $\tau_n$ of $f$, the
highest power has strictly higher pole order than every lower one.
It cannot cancel, a contradiction. This proves transcendence and
injectivity of rational evaluation $K(F)\to K(f)$.

For each of the finitely many differences $Q_i-Q_j$, its leading
coefficient is regular and nonzero at every sufficiently late
$\tau_n$. Thus $Q_i(t,f)-Q_j(t,f)$ has a genuine pole there.
Every rational coefficient in $K(f)$ is a meromorphic germ at that
point, and a nonzero one is not identically zero there by the
injectivity just proved and the identity theorem. Apply
Lemma~\ref{exp:lem:localexp} at any such late pole to exclude a
nonzero relation. The analytic identity theorem legitimizes
passing from an identity on an open subset of the punctured disk
to the local punctured neighborhood of this pole; all other
singularities are discrete away from the center.
\end{proof}

\begin{proof}[Proof of Theorem~\ref{exp:thm:main}]
For $b=P_B$, include $Q_0=0$ and $Q_j=P_j$, $j\ge1$,
in Lemma~\ref{exp:lem:independence}. By Lemma~\ref{ana:lem:Tmeromorphic},
$T_C\in K$. In \eqref{exp:eq:excluded} the coefficient of $e^{Q_0}$
is $R_0(t,f)-T_C(z;t)$. Independence makes every coefficient
zero. Injectivity of rational evaluation then gives
$R_0(t,F)=T_C(z;t)\in K$, as claimed.

The other inputs are invertible rational affine changes of $f$.
Equations~\eqref{ana:eq:F} and \eqref{ana:eq:FP} give
\begin{equation}
 z^2t^2U_B(z;t)=zt(G_B-t)-(1-t+zt^2)f(t).
 \label{exp:eq:Utransfer}
\end{equation}
For fixed $x_0\in(1/2,3/2)$ with $x_0\ne1$,
\eqref{ana:eq:BZ} and \eqref{ana:eq:FP} give
\begin{equation}
 Z_B(z;x_0,t)=A_{x_0}(t)
       -\frac{x_0(1-t)}{t(x_0-t)(1-t+zx_0t)}f(t),
 \label{exp:eq:fulltransfer}
\end{equation}
where, writing $F_0=F_B(z;x_0)$ and $P_0=P_B(z;x_0)$,
\[
 A_{x_0}(t)=\frac{1}{1-t+zx_0t}
 \left[F_0-P_0+\frac{zx_0(1-t)P_0}{x_0-t}
 +\frac{zx_0(F_0-G_B+1)}{x_0-1}
 +\frac{zx_0(G_B-t-F_0)}{x_0-t}\right].
\]
These constants are defined by convergent series. At $x_0=1$,
the quotient $(F_0-G_B+1)/(x_0-1)$ has the removable value
$\partial_x F_B(z;1)$. Taking this limit in \eqref{exp:eq:fulltransfer}
gives a rational affine expression for $Z_B(z;1,t)$ whose
multiplier of $f$ is $-1/[t(1-t+zt)]\ne0$.
Both affine
multipliers are nonzero elements of $K$, even if they vanish at
$t_*$ itself. Substitution therefore induces an automorphism of
$K(F)$, preserves positive polynomial degree and the degree of
each exponent difference, and transfers the result to both inputs.
\end{proof}

\section{Signed catalytic specializations}\label{sec:signed}
The suggestion to compare $Z_C(x,-x)$ with $-Z_C(-x,x)$ leads to an
exact finite identity, although that two-term equality does not hold.
For arbitrary catalysts $a,b$, write
\[
 \mathcal W(a,b)=\frac{\kappa(a,b)}{za}Z_C(a,b),\qquad
 \kappa(a,b)=(1-b)(a-b)+za(a-1).
\]
This is interpreted in $\mathbb Q(a,b)((z))$, with removable factors
canceled at any indicated specialization.
\begin{proposition}\label{prop:signed}
For any $a,c,u,v,w$, the six-term identity
\begin{align}
 &(v-w)\bigl(\mathcal W(a,u)-\mathcal W(c,u)\bigr)\notag\\
 &\quad +(w-u)\bigl(\mathcal W(a,v)-\mathcal W(c,v)\bigr)
 +(u-v)\bigl(\mathcal W(a,w)-\mathcal W(c,w)\bigr)=0
 \label{eq:sixterm}
\end{align}
holds. In contrast,
\[
 Z_C(z;x,-x)+Z_C(z;-x,x)=2x^2z^5+O(z^6).
\]
\end{proposition}
\begin{proof}
Put $U(a)=Z_C(a,1)$, $D(b)=Z_C(b,b)$ and
$V(a)=(F_C(a)-F_C(1))/(a-1)$. The cleared C equation gives
\[
 \mathcal W(a,b)=B_0(a)-bA_0(a)+C_0(b),
\]
where
\[
 \begin{aligned}
 A_0(a)&=U(a)+V(a)+P_C(a),\\
 B_0(a)&=a(U(a)+V(a))+P_C(a)-F_C(a),\\
 C_0(b)&=(b-1)(D(b)+P_C(b))+F_C(b).
 \end{aligned}
\]
Consequently $\mathcal W(a,b)-\mathcal W(c,b)$ is affine in $b$;
its second divided difference is zero, proving \eqref{eq:sixterm}.
In particular, taking $a=x,c=-x,u=x,v=-x,w=ix$, with $i^2=-1$,
gives an identity among signed and fourth-root specializations.
For the two-term proposal, direct coefficient extraction from the
unique C system gives
\[
 \begin{aligned}
 [z^2]Z_C(x,-x)&=x,& [z^3]Z_C(x,-x)&=3x,\\
 [z^4]Z_C(x,-x)&=9x+x^3,&
 [z^5]Z_C(x,-x)&=30x+x^2+6x^3.
 \end{aligned}
\]
Replacing $x$ by $-x$ proves the claimed first discrepancy.
\end{proof}
Independent literal reconstruction from admissible sites verifies the
full bivariate series $Z_C$ through length $8$. The supplementary program
\texttt{signed\_check.py} also verifies the signed and fourth-root
specialization of~\eqref{eq:sixterm} through $z^{60}$ using exact
Gaussian-integer arithmetic. These finite checks corroborate the
coefficient calculation and the identity proved above.

Equation~\eqref{eq:sixterm} follows from the C system alone. The equality
of total counts in Theorem~\ref{count:thm:main} additionally proves the
slack-series relation in the conclusion. A further transformation between
the full refined B and C series is a separate question.

\section{Conclusion and further relations}
The two enumeration problems considered here have explicit solutions.
For Class~152, the common count is a partial central-binomial sum.
For Class~207, Theorem~\ref{count:thm:main} proves equality of all three
counting series and gives their common Lambert series expression.
The C kernel and the published B tree lead to different scalar formulas;
the cubic limit of the nonterminating $q$-Saalsch\"utz summation is the
identity that connects them. The A--C step is the published bijection
of Callan--Mansour~\cite[Section~4.10]{CM}.

Equality of the totals also gives a relation between the slack series.
Put $\Delta F=F_C-F_B$ and $\Delta P=P_C-P_B$.
Subtracting the common $P$ equations and using $F_C(1)=F_B(1)$ gives
\[
 (x-1)(1-zx)\Delta P(x)=zx\Delta F(x),\qquad
 \Delta F(x)=-\frac{(1-x)(1-zx)}{zx}\Delta P(x).
\]
This proves the proposed relation between the two boundary systems.
The equality of the totals coexists with the different analytic behavior
of the refined series. An additional relation between the full $Z_B$ and
$Z_C$, or between the C diagonal and a transformed B series, remains open,
including the search for a relation valid through $z^{50}$. The six-term
C identity of Section~\ref{sec:signed} does not resolve this question.

Theorems~\ref{ana:thm:main} and~\ref{ana:thm:differential} exclude exact
polynomial relations for the two stated pairs and finite linear
relations among the specified catalytic derivatives. They impose no
fixed degree or derivative-order bound. Nonlinear differential relations
are a different question. Theorem~\ref{exp:thm:main} excludes the stated
separated polynomial-exponential expressions for each sufficiently small
fixed positive $z$. Proper rational exponents, quotients of exponential
sums, other catalytic substitutions, length derivatives, coefficient
transforms, and infinite constructions have separate scopes.
An explicit independent evaluation of the boundary $Z_C(x,0)$ and a
description of nonlinear differential relations remain further questions.

\section*{Declarations}
\paragraph{Funding}
This research received no specific grant from any funding agency in the public, commercial, or not-for-profit sectors.

\paragraph{Declaration of competing interest}
The authors declare that they have no competing interests.

\paragraph{Data availability}
The data and verification code supporting the finite computations are
available with reproduction instructions in the
\href{https://github.com/the-omega-institute/inversion-sequence-classes-152-207}{companion code and data repository}.
Section~\ref{sec:checks} describes the finite checks, and
\ref{sec:program} gives the complete recurrence program.
The ordinary proofs are contained in the main article and this expanded supplement.

\paragraph{Declaration of generative AI and AI-assisted technologies in the manuscript preparation process}
Haobo Ma and Wenlin Zhang used AI tools, namely ChatGPT Pro (GPT-6 Pro), Codex (gpt-6.1-sol) and Claude Opus 5.5 (via Claude Code), for literature search, computation and double-checking, including formal verification of the results in Lean. All proofs can be checked by hand, and the computations can be reproduced from the verification package. H. Ma and W. Zhang take full responsibility for their contributions.

\appendix
\section{Python program for 207(2)}\label{sec:program}
The file \texttt{class207\_recurrence.py} in the accompanying supplement
uses only the Python standard library. From the supplement directory, run
\texttt{python3 class207\_recurrence.py 12} to reproduce the $C$ column
of the table, or add \verb|--check counts.csv| to compare both archived
Class~207 columns. Run \texttt{python3 check\_counts.py} to regenerate
those columns independently and perform the literal pattern checks.
In the code, \texttt{row[s][k]} is $\sum_{i<k}a_r(s,i)$,
\texttt{diag[m][s]} is $\sum_{t=1}^s a_r(t,m-t)$, and
\texttt{correction[s]} is $\sum_{t=1}^s(b_r(t)-a_r(t,0))$.
Thus \eqref{eq:a} uses diagonal $m=s+k$ and \eqref{eq:b} uses $m=s+1$.
\begin{lstlisting}
#!/usr/bin/env python3
"""Class 207(2): exact descendant recurrence, using integer partial sums."""
import argparse
import csv


def counts(N):
    """Return the counts at lengths 0,...,N; O(N^3) integer operations."""
    if N < 0:
        raise ValueError("N must be nonnegative")
    if N == 0:
        return [1]
    # A[s][k] = a_r(s,k), B[s] = b_r(s), on s+k <= N-r.
    A = [[]] + [[1] * (N - s + 1) for s in range(1, N + 1)]
    B = [0] + [1] * N
    answer = [1, 1]
    for L in range(N, 1, -1):
        # Row sums, diagonal sums, and sums of b_r(t)-a_r(t,0).
        row = [[]]
        for s in range(1, L + 1):
            p = [0]
            for v in A[s]:
                p.append(p[-1] + v)
            row.append(p)
        diag = [[0] for _ in range(L + 1)]
        for m in range(1, L + 1):
            for s in range(1, m + 1):
                diag[m].append(diag[m][-1] + A[s][m - s])
        correction = [0]
        for s in range(1, L + 1):
            correction.append(correction[-1] + B[s] - A[s][0])
        A = [[]] + [
            [row[s + 1][k] + diag[s + k][s] + correction[s]
             for k in range(L - s)] for s in range(1, L)]
        B = [0] + [B[s + 1] + diag[s + 1][s] + correction[s]
                   for s in range(1, L)]
        answer.append(B[1])
    return answer


if __name__ == "__main__":
    parser = argparse.ArgumentParser(description=__doc__)
    parser.add_argument("N", nargs="?", type=int, default=12)
    parser.add_argument("--check", metavar="CSV", help="compare n,B207,C207 CSV")
    args = parser.parse_args()
    values = counts(args.N)
    if args.check:
        with open(args.check, newline="") as stream:
            rows = {int(r["n"]): r for r in csv.DictReader(stream)}
        for n, value in enumerate(values):
            if n not in rows or any(int(rows[n][c]) != value
                                    for c in ("B207", "C207")):
                raise SystemExit(f"Reference mismatch at n={n}")
        print(f"Verified all {len(values)} terms (0 <= n <= {args.N}).")
    else:
        print("n,count_207_2")
        for n, value in enumerate(values):
            print(f"{n},{value}")
\end{lstlisting}